\documentclass[12 pt]{article}
\usepackage{amssymb,epsfig,amsmath}
\usepackage[dvips]{color}
\usepackage[T1]{fontenc}
\usepackage[french]{babel}
\usepackage[utf8]{inputenc}
\usepackage{comment}
\usepackage{hyperref}
\usepackage{graphicx}
\usepackage{multicol}
\usepackage{caption}
\newtheorem{defi}{D\'efinition}[section]
\newtheorem{prop}[defi]{Proposition}
\newtheorem{theo}[defi]{Th\'eor\`eme}
\newtheorem{conj}[defi]{Conjecture}
\newtheorem{lemm}[defi]{Lemme}
\newtheorem{coro}[defi]{Corollaire}
\newtheorem{rema}[defi]{Remarque}
\newtheorem{exem}[defi]{Exemple}
\newtheorem{exems}[defi]{Exemples}

\newcommand{\bdefi}{\begin{defi}}
\newcommand{\edefi}{\end{defi}}
\newcommand{\bprop}{\begin{prop}}
\newcommand{\eprop}{\end{prop}}
\newcommand{\btheo}{\begin{theo}}
\newcommand{\etheo}{\end{theo}}
\newcommand{\blemm}{\begin{lemm}}
\newcommand{\brema}{\begin{rema}}
\newcommand{\erema}{\end{rema}}
\newcommand{\bexer}{\begin{exem}}
\newcommand{\eexer}{\end{exem}}
\newcommand{\bexems}{\begin{exems}}
\newcommand{\eexems}{\end{exems}}
\newcommand{\bconj}{\begin{conj}}
\newcommand{\econj}{\end{conj}}
\newcommand{\elemm}{\end{lemm}}
\newcommand{\bcoro}{\begin{coro}}
\newcommand{\ecoro}{\end{coro}}

\newcommand{\rem}{\noindent{\bf Remarque. }}

\usepackage{mathrsfs}
\renewcommand\mathcal{\mathscr}

\newcommand{\T}{{\cal T}}

\newcommand{\M}{{\cal M}}

\newcommand{\E}{{\cal E}}
\newcommand{\F}{{\cal F}}
\newcommand{\V}{{\cal V}}

\newcommand{\SSSS}{{\cal S}}

\newcommand{\PPPP}{{\cal P}}
\newcommand{\C}{{\cal C}}
\newcommand{\I}{{\cal I}}
\newcommand{\Q}{{\cal Q}}

\newcommand{\maths}[1]{{\mathbb #1}}  

\newcommand{\RR}{\maths{R}}
\newcommand{\NN}{\maths{N}}
\newcommand{\CC}{\maths{C}}

\newcommand{\DD}{\maths{D}}

\newcommand{\ZZ}{\maths{Z}}
\newcommand{\PP}{\maths{P}}

\newcommand{\CCC}{{\mathfrak C}}

\newcommand{\ga}{\gamma}
\newcommand{\Ga}{\Gamma}

\newcommand{\cqfd}{\hfill$\Box$}

\newcommand{\CAT}{\operatorname{CAT}}

\usepackage{geometry}
\usepackage{graphicx}

\newcommand{\dddp}{\partial_{\infty}^2\widetilde{\Sigma}}

\newcommand{\Flat}{\operatorname{Flat}}

\newcommand{\Mix}{\operatorname{Mix}}
\newcommand{\R}{\mathcal{R}}

\newcommand{\revet}{(\widetilde{\Sigma},[\widetilde{q}])}

\newcommand{\srfce}{({\Sigma},{[q]})}

\newcommand{\srfcem}{({\Sigma},m)}
\newcommand{\Id}{\operatorname{Id}}

\newcommand{\Area}{\operatorname{Area}}

\newcommand{\peri}{\operatorname{peri}}
\newcommand{\revetn}{(\widetilde{\Sigma},[\widetilde{q}_n])}
\newcommand{\revetm}{(\widetilde{\Sigma},\widetilde{m})}
\newcommand{\Ax}{\operatorname{Ax}}
\newcommand{\Sing}{\operatorname{Sing}}

\usepackage{fullpage}

\renewcommand{\bdefi}{\begin{defi}}
	\renewcommand{\edefi}{\end{defi}}
\renewcommand{\bprop}{\begin{prop}}
	\renewcommand{\eprop}{\end{prop}}
\renewcommand{\btheo}{\begin{theo}}
	\renewcommand{\etheo}{\end{theo}}
\renewcommand{\blemm}{\begin{lemm}}
	\renewcommand{\brema}{\begin{rema}}
		\renewcommand{\erema}{\end{rema}}
	\renewcommand{\bexer}{\begin{exem}}
		\renewcommand{\eexer}{\end{exem}}
	\renewcommand{\bexems}{\begin{exems}}
		\renewcommand{\eexems}{\end{exems}}
	\renewcommand{\bconj}{\begin{conj}}
		\renewcommand{\econj}{\end{conj}}
	\renewcommand{\elemm}{\end{lemm}}
\renewcommand{\bcoro}{\begin{coro}}
	\renewcommand{\ecoro}{\end{coro}}

\renewcommand{\rem}{\noindent{\bf Remark. }}

\newtheorem{defiA}{Definition}[section]
\newtheorem{propA}[defi]{Proposition}
\newtheorem{theoA}[defi]{Theorem}
\newtheorem{conjA}[defi]{Conjecture}
\newtheorem{lemmA}[defi]{Lemma}
\newtheorem{coroA}[defi]{Corollary}
\newtheorem{remaA}[defi]{Remark}
\newtheorem{exemA}[defi]{Example}
\newtheorem{exemsA}[defi]{Examples}

\newcommand{\bdefiA}{\begin{defiA}}
	\newcommand{\edefiA}{\end{defiA}}
\newcommand{\bpropA}{\begin{propA}}
	\newcommand{\epropA}{\end{propA}}
\newcommand{\btheoA}{\begin{theoA}}
	\newcommand{\etheoA}{\end{theoA}}
\newcommand{\blemmA}{\begin{lemmA}}
		\newcommand{\elemmA}{\end{lemmA}}
	\newcommand{\bremaA}{\begin{remaA}}
		\newcommand{\eremaA}{\end{remaA}}
	\newcommand{\bexerA}{\begin{exemA}}
		\newcommand{\eexerA}{\end{exemA}}
	\newcommand{\bexemsA}{\begin{exemsA}}
		\newcommand{\eexemsA}{\end{exemsA}}
	\newcommand{\bconjA}{\begin{conjA}}
		\newcommand{\econjA}{\end{conjA}}
\newcommand{\bcoroA}{\begin{coroA}}
	\newcommand{\ecoroA}{\end{coroA}}
\newcommand{\demA}{\noindent{\bf Proof. }}
\newcommand{\remA}{\noindent{\bf Remark. }}

\title{Geometric compactification of moduli spaces of half-translation structures on surfaces.}
\author{Thomas~Morzadec}

\begin{document}
\maketitle 

\sloppy
{\bf Résumé:~ } Dans cet article, on construit une compactification équivariante de l'espace $\PP\Flat(\Sigma)$ des classes d'homothétie de structures de demi-translation
sur une surface $\Sigma$ compacte, connexe, orientable. On définit l'espace $\PP\Mix(\Sigma)$ des classes d'homothétie de structures mixtes sur $\Sigma$, qui sont des 
structures arborescentes,
au sens de Drutu et Sapir, $\CAT(0)$, dont les pièces sont des arbres réels ou des complétés de surfaces munies de structures de demi-translation. 

En munissant $\Mix(\Sigma)$ 
de la topologie de Gromov équivariante, et en utilisant des techniques de cônes asymptotiques à la Gromov, on montre que $\PP\Mix(\Sigma)$ est une compactification
équivariante
de $\PP\Flat(\Sigma)$, ce qui nous permet de comprendre géométriquement les dégénérescences de structures de demi-translation sur $\Sigma$. On compare ensuite cette compactification à celle 
de Duchin-Leininger-Rafi, qui utilise des courants géodésiques, en passant par les distances de translation des éléments du groupe de revêtement de $\Sigma$. 

\medskip

{\bf Abstract:~ } In this paper, we give an equivariant compactification of the space $\PP\Flat(\Sigma)$ of homothety classes of half-translation structures on a compact, 
connected, orientable surface $\Sigma$. We introduce the space $\PP\Mix(\Sigma)$ of homothety classes of mixed structures on $\Sigma$, that are $\CAT(0)$ tree-graded spaces
in the sense of Drutu and Sapir, with pieces which are $\RR$-trees and completions of surfaces endowed with  half-translation structures. 

Endowing $\Mix(\Sigma)$ with the  equivariant Gromov topology, and using asymptotic cone techniques, we prove that $\PP\Mix(\Sigma)$ is an equivariant compactification of 
$\PP\Flat(\Sigma)$, thus allowing us to understand in a geometric way the degenerations of half-translation structures on $\Sigma$. We finally compare our compactification to the
one of Duchin-Leininger-Rafi, based on geodesic currents on $\Sigma$, by the mean of the translation distances of the elements of the covering group of $\Sigma$. 
\footnote{keywords: half-translation surface, holomorphic quadratic differential, tree-graded space, mixed structure on surfaces, asymptotic cone, flat surface with singularities,
geodesic lamination, compactification. AMS code 57M50, 30F60, 30F30, , 32G15, 51M05, 53C23, 53C45}

\section{Introduction.}

The goal of this paper is to construct and to describe a geometric  compactification, natural under the action of the mapping class group,  of the space
of homothety classes of half-translation stuctures on a compact surface,
endowed with the  equivariant Gromov topology.
  It is part of the wide field of study of deformations of geometric 
structures on surfaces (see for instance \cite{Goldman88}). Let $\Sigma$ 
be a compact, connected, orientable surface of genus $g > 2$, without 
boundary for simplicity (see \cite{Morzy4}) for the general case). After the founding fathers Gauss and 
Riemann who have studied conformal geometry on surfaces, the Teichmüller 
spaces $\T(\Sigma)$ of isotopy classes of hyperbolic metrics on 
$\Sigma$, have been studied by for instance Fricke, Klein, Fenchel, 
Nielsen, and the moduli spaces of real projective structures by for 
instance Goldman and Choi \cite{ChoGol97}. The analysis of the space 
$\operatorname{Flat}(\Sigma)$ of half-translation structures on a 
$\Sigma$ is currently blooming, with the works notably of Calta, Eskin, 
Hubert, Lanneau, Masur, McMullen, Möller, Myrzakhani, Schmidt, Smillie, 
Veech, Weiss, Yoccoz and Zorich. When these deformations spaces are non 
compact,
it is important and useful to consider the asymptotic behavior of the 
sequences of geometric structures that leave all compact subsets. Only 
few results are known about the compactifications of spaces of geometric 
structures, except for the Teichmüller space,
for which several compactifications have been built, notably by Thurston 
(see [FLP]), and
also by Bestvina, Morgan, Paulin, Shalen. A compactification of 
$\operatorname{Flat}(\Sigma)$ has been recently proposed in \cite{DucLeiRaf10}. 
This article aims at proposing a new one.

\medskip

Let $\Sigma$ be a connected, orientable surface.  
A {\it half-translation structure} (or flat structure with conical singularities and holonomies in $\{\pm\Id\}$) on $\Sigma$ is the data consisting of a 
(possibly empty) discrete subset $Z$ of $\Sigma$ and of a
Euclidean metric on $\Sigma-Z$ with conical  singular points of angles of the form $k\pi$, with $k\in\NN$ and $k\geqslant 3$, at each point of $Z$, such
that the holonomy of every
piecewise $\C^1$ loop of $\Sigma-Z$ is contained in $\{\pm\Id\}$. As the set $\Flat(\Sigma)$ of isotopy classes
of half-translation structures on $\Sigma$ identifies with 
the quotient by $\operatorname{SO}(2)$ of the set of isotopy classes of holomorphic quadratic differentials on $\Sigma$, we will denote by $[q]$, with $q$ 
a holomorphic quadratic differential on $\Sigma$, a half-translation structure on $\Sigma$. We refer for instance to \cite{Wright,Zorich06} and 
Section \ref{structuresplates} for background on half-translation structures.
 
\medskip

Assume that $\Sigma$ is compact and that $\chi(\Sigma)<0$. Let us define the mixed structures on $\Sigma$, which are $\CAT(0)$ tree-graded spaces in the sense of Drutu-Sapir,
arising, as we shall see, as geometric degenerations of elements of $\Flat(\Sigma)$.
Recall that a {\it tree-graded space} (see \cite[Def.~1.10]{Drutu05}) is the data of a complete geodesic metric space $X$ and of a covering $\cal{P}$ of $X$ by
closed geodesic subsets of $X$, whose elements are called {\it pieces}, such that:

\medskip
\noindent
$\bullet$~ two distinct pieces have at most one common point;

\medskip
\noindent
$\bullet$~ any simple geodesic triangle of $X$ is contained in a single piece.

\medskip

We will need the following definitions. Let $p:\widetilde{\Sigma}\to\Sigma$ be a universal cover of $\Sigma$ with covering group $\Gamma$. Let $\Sigma_0$ be 
a {\it tight} subsurface of $\Sigma$, that is  a proper closed subsurface with smooth boundary, such that no connected component of $\Sigma_0$ is a disk or a pair of pants,  
no connected component of $\Sigma-\Sigma_0$ is a disk, and no cylinder component of $\Sigma_0$
can be homotoped in another connected component of $\Sigma_0$. 
Let $W$ be a connected component of $\Sigma_0$ or of $\Sigma-\Sigma_0$ and  $\widetilde{W}$ be a connected component of the preimage of $W$ in $\widetilde{\Sigma}$.
 We associate to $\widetilde{W}$ a complete geodesic metric space $X_{\widetilde{W}}$ of the following type:

\medskip
\noindent
$\bullet$~ if $\widetilde{W}$ is a strip, then $X_{\widetilde{W}}$ is empty;

\medskip
\noindent
$\bullet$~ if $\widetilde{W}$ is a connected component of $\widetilde{\Sigma}_0$ which is not a strip or if $\widetilde{W}$ is a connected component of 
$\widetilde{\Sigma}-\widetilde{\Sigma}_0$ and $W$ is a pair of 
pants, then $X_{\widetilde{W}}$ is a point;

\medskip
\noindent
$\bullet$~ if $\widetilde{W}$ is a connected component of $\widetilde{\Sigma}-\widetilde{\Sigma}_0$ and $W$ is neither an annulus nor a pair of pants, 

\medskip 

- either $X_{\widetilde{W}}$ is the completion of the lift to $\widetilde{W}$ of a half-translation structure $[q_W]$ on $W$ which extends to the filled in surface, possibly
with singularities of angle $\pi$.

\medskip

- or $X_{\widetilde{W}}$ is the dual $\RR$-tree to a filling measured hyperbolic lamination on $W$ (endowed with any complete hyperbolic metric, see for instance
\cite{Bonahon97} and \cite[§~2.3]{Ota96} and Subsection \ref{areazero} for the definitions).

\medskip

Finally, to every proper homotopy class $\widetilde{c}$ of boundary components of $\widetilde{\Sigma}_0$, we associate a compact interval $X_{\widetilde{c}}$ of $\RR$,
called an {\it edge}. We refer to the picture of Section \ref{mixedstructure} for a graphical understanding of the following definition.

\bdefiA A mixed structure on $\Sigma$ is a metric space $X$ such that there exists a (possibly empty) tight, proper subsurface
$\Sigma_0$, as above, such that $X$ is obtained by gluing some complete metric spaces $X_{\widetilde{W}}$ and 
$X_{\widetilde{c}}$, where $\widetilde{W}$ is a connected components of $\widetilde{\Sigma}_0$ or of $\widetilde{\Sigma}-\widetilde{\Sigma}_0$, 
and $\widetilde{c}$ is a proper homotopy class of boundary components of $\widetilde{\Sigma}_0$, as above (see Section \ref{mixedstructure} for the precisions), such that:

\medskip
\noindent
$\bullet$~ if $\widetilde{W}$ is a strip of $\widetilde{\Sigma}-\widetilde{\Sigma}_0$, and if $\widetilde{c}$ is the proper homotopy class of the boundary components 
of $\widetilde{W}$, then the length of $X_{\widetilde{c}}$ is positive; 

\medskip
\noindent
$\bullet$~ if $\widetilde{W}$ is a connected component of $\widetilde{\Sigma}-\widetilde{\Sigma}_0$ whose image in $\Sigma$ is a pair of pants, there exists at least
a proper homotopy class $\widetilde{c}$ of boundary components of $\widetilde{W}$ such that the length of $X_{\widetilde{c}}$ is positive.
\edefiA

We will see (see Section \ref{mixedstructure}) that a mixed structure is a $\CAT(0)$ tree-graded space, and that the actions of the stabilizers 
of $\widetilde{W}$ in $\Gamma$ on the pieces $X_{\widetilde{W}}$ and $X_{\widetilde{c}}$ as above glue together to give an isometric action of $\Gamma$ on $X$. We will see that the action of $\Gamma$ on $X$ defines the subsurface $\Sigma_0$, up to isotopy.
Moreover, a mixed structure is defined (up to $\Gamma$-equivariant isometry) by the half-translation structures and the measured hyperbolic laminations on the
connected components of $\Sigma-\Sigma_0$, and by the lengths of the edges, that define its pieces (see Section \ref{mixedstructure}). 

\medskip

Let $\Mix(\Sigma)$ be the space of $\Gamma$-equivariant isometry classes of mixed structures on $\Sigma$. We endow it with the equivariant Gromov topology 
(see \cite{Paulin89a,Paulin09a} and Section \ref{mixedstructure}) which is the topology such that two elements are close if they contain large finite subsets which are almost isometric in an 
equivariant way under a large finite subset of $\Gamma$. 
Let $\PP\Mix(\Sigma)$ be the space of homothety classes of mixed structures on $\Sigma$, endowed with the quotient  topology. We identify $\Flat(\Sigma)$ with  
the space of $\Gamma$-equivariant isometry classes of  $\Gamma$-invariant half-translation structures on $\widetilde{\Sigma}$, which is a subset of $\Mix(\Sigma)$.
Let $\PP\Flat(\Sigma)$ be its image in $\PP\Mix(\Sigma)$. The two main results of this paper are the following ones.

\btheoA The space $\PP\Flat(\Sigma)$ is an open and dense subset of $\PP\Mix(\Sigma)$, and $\PP\Mix(\Sigma)$  is compact. The action of the mapping class group of
$\Sigma$
on $\PP\Flat(\Sigma)$ extends continuously to $\PP\Mix(\Sigma)$ 
\etheoA

Let $(X,d)$ be a metric space with an isometric action of a group $\Gamma$. For all $\gamma\in\Gamma$, the {\it translation distance} of $\gamma$ in $X$ is
$\ell_{X}(\gamma)=\inf_{x\in X}d(x,\gamma x)$. Let $\PP\RR^\Gamma=((\RR^{+})^\Gamma-\{0\})/\RR^{+*}$. We denote by $[X]$ the image of $X\in\Mix(\Sigma)$ in
$\PP\Mix(\Sigma)$ and by $[x_\gamma]_{\gamma\in\Gamma}$ the image of $(x_\gamma)_{\gamma\in\Gamma}\in(\RR^+)^\Gamma-\{0\}$ in $\PP\RR^\Gamma$. 

\btheoA\label{theorem3} The map $[X]\mapsto[\ell_X(\gamma)]_{\gamma\in\Gamma}$ is an embedding  of $\PP\Mix(\Sigma)$ onto its image.
\etheoA

In Section \ref{notation}, we recall some definition and basic facts about half-translation structures on a surface.  In Section \ref{action}, 
we give some results about the action of the covering 
group on the universal cover of a surface endowed with a half-translation structure. In Section \ref{tools}, we recall two tools that we will need in this paper,
which are the subsurfaces filled by a family of homotopy classes of simple closed curves and the generalized subsurfaces with geodesic boundary, for a half-translation 
structure. In Section \ref{ultralimite}, we recall what are the ultralimits of sequences of metric spaces in Gromov's sense. In Section \ref{sequenceofhalftranslations}, we study the ultralimits of universal covers of sequences of 
half-translation structures on a surface. Finally, in Section \ref{mixedstructure}, we define the mixed structures, we recall the definition of the
equivariant Gromov topology, and
we prove the above  two theorems. Using Theorem
\ref{theorem3}, we construct a mapping class group-equivariant homeomorphism between our compactification of $\PP\Flat(\Sigma)$ and the one constructed by Duchin-Leininger-Rafi in \cite{DucLeiRaf10}.

Our approach is fundamental in understanding in a geometric way the degenerations of half-translation structures, and allow an extension (work in progress) to flat structures
with conical singularities and finite holonomy groups (for instance the quotient by $\operatorname{SO}(2)$ of the set of isotopy classes of cubic holomorphic differentials on $\Sigma$,
whose moduli spaces have started to be studied for instance by Benoist-Hulin \cite{BenHulin13}, by Labourie \cite{Labourie07} and Loftin \cite{Loftin07}).

\medskip\noindent
{\small{\it Acknoledgement: I want to thank Frédéric Paulin for many advices and corrections that have deeply improved the writting of this paper.}}
\section{Notation and background.}\label{notation}

\subsection{Half-translation structures.}\label{structuresplates}

In the whole paper, we will use the definitions and notation of \cite{BriHae99} for a surface endowed with a distance
$(\Sigma,d)$: (locally) $\CAT(0)$, $\delta$-hyperbolic,... 
Notably, a {\it geodesic} (resp. a {\it local geodesic}) 
of $(\Sigma,d)$ is an isometric (resp. locally isometric) map $\ell:I\to \Sigma$, where $I$ is an interval of $\RR$. It will be called a 
  {\it segment}, a {\it ray} or a {\it  geodesic line} of $(\Sigma,d)$ if $I$ is respectively a compact interval, a closed half line 
  (generally $[0,+\infty[$)
  or $\RR$.
  If there is no precision, a {\it geodesic} is a geodesic line. Let $\Sigma$ be a connected, orientable surface, with empty boundary. 

\bdefiA
A half-translation structure (or flat structure with conical singularities and holonomies in $\{\pm\Id\}$) on a surface $\Sigma$ is the data of a (possibly empty) discrete
subset $Z$ of $\Sigma$
and a Euclidean metric on $\Sigma-Z$ with conical singularity of angle $k_z\pi$ at each $z\in Z$, with $k_z\in\NN$, $k_z\geqslant 3$, 
     such that the holonomy of every piecewise $\C^1$ loop in 
 $\Sigma-Z$ is contained in $\{\pm\Id\}$.
\edefiA

 As the set $\Flat(\Sigma)$ of isotopy classes
of half-translation structures on $\Sigma$ identifies with 
the quotient by $\operatorname{SO}(2)$ of the set of isotopy classes of holomorphic quadratic differentials on $\Sigma$, we will denote by $[q]$, with $q$ 
a holomorphic quadratic differential on $\Sigma$, a half-translation structure on $\Sigma$ (see \cite[§~2.5]{Morzy1} and \cite[Def.~1.2~p.~2]{Strebel84} for the definition
of a holomorphic quadratic differential).
 A half-translation structure defines a geodesic distance $d$ on $\Sigma$ that is locally $\CAT(0)$. We will call {\it local flat geodesic} a local geodesic of a half-translation 
 structure. A continuous map $\ell:\RR\to\Sigma$ is a local flat geodesic if and only if it satisfies (see \cite[Th.~5.4~p.24]{Strebel84} and
 \cite[Th.~8.1~p.~35]{Strebel84}): for every $t\in\RR$,

\medskip
\noindent
$\bullet$~  if $\ell(t)$ does not belong to $Z$, there exists a neighborhood $V$ of $t$ in $\RR$ such that $\ell_{|V}$ is a Euclidean segment
(hence, $\ell_{|V}$ has constant direction);
 
\noindent
$\bullet$~  if $\ell(t)$ belongs to $Z$, then the two angles defined by $\ell([t,t+\varepsilon[)$ and $\ell(]t-\varepsilon,t])$, 
with $\varepsilon>0$ small enough, measured in the two connected components of $U-\ell(]t-\varepsilon,t+\varepsilon[)$, 
   with $U$  a small enough neighborhood of $\ell(t)$, are at least $\pi$.

\begin{center}
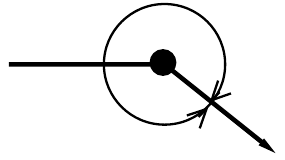
\end{center}

Let $\C(\Sigma)$ and $\SSSS(\Sigma)$ be the sets of free homotopy classes of essential closed and simple closed curves on $\Sigma$. Let $\alpha\in\C(\Sigma)$. 
If $m$ is a hyperbolic metric on $\Sigma$, then $\alpha$ has a unique $m$-geodesic representative. However, if $[q]$ is a half-translation structure on $\Sigma$, then 
$\alpha$ has
at least a $[q]$-geodesic representative, but it may not be unique. In that case, the set of geodesic representatives foliates a maximal flat cylinder (with singular points 
on its boundary), whose interior is embedded into $\srfce$. We denote by $\ell_m(\alpha)$ and $\ell_{[q]}(\alpha)$ the lengths of the geodesic representatives of $\alpha$.

In this article, we will consider (complete) hyperbolic metrics and half-translation structures on $\Sigma$. Whereas a hyperbolic (local) geodesic is
uniquely determined by its image (up to changing the origin), a flat (local) geodesic is not. However, we will sometimes still denote by $\ell$ the image of a flat (local)
geodesic $\ell$, 
if there is no confusion. 
Let $d$ be the distance defined by a half-translation structure
or a hyperbolic metric on $\Sigma$, and let $p:(\widetilde{\Sigma},\widetilde{d})\to(\Sigma,d)$ be a universal cover of covering group $\Gamma$. Since $(\Sigma,d)$ is complete and locally $\CAT(0)$, 
according to the theorem of Cartan-Hadamard, the space $(\widetilde{\Sigma},\widetilde{d})$ is complete and $\CAT(0)$. We will denote by 
$\partial_\infty\widetilde{\Sigma}$ the boundary at infinity of $(\widetilde{\Sigma},\widetilde{d})$ and $\partial^2_\infty\widetilde{\Sigma}=
\partial_\infty\widetilde{\Sigma}\times\partial_\infty\widetilde{\Sigma}-\Delta$ (with $\Delta=\{(x,x),x\in\partial_\infty\widetilde{\Sigma}\}$). If $\Sigma$ is a (possibly trivial) cover of 
a compact surface, which will always be the case in this article, the boundary at infinity $\partial_\infty\widetilde{\Sigma}$ does not depend on the complete locally
$\CAT(0)$ distance on $\Sigma$, up to a unique $\Gamma$-equivariant homeomorphism. If ${g}$ is a geodesic of $(\widetilde{\Sigma},\widetilde{d})$, we denote by
$E(g)=(g(-\infty),g(+\infty))\in\partial_\infty^2\widetilde{\Sigma}$ its pair of points at infinity. 

\section{Isometric action of the covering group on the universal cover of a surface endowed with a half-translation structure.}\label{action}

Let $(X,d)$ be a complete $\CAT(0)$ metric space. If $(X,d)$ is endowed with an isometric action of
a group $\Gamma$, for every $\gamma\in\Gamma$,
the {\it translation distance} of $\gamma$ is 
$\ell_X(\gamma)=\underset{x\in X}\inf d(x,\gamma x)$. The element $\gamma$ is said to be 
{\it elliptic} if it has a fixed point in $X$, {\it parabolic} if $\ell_X(\gamma)=0$ but
$\gamma$
has no fixed point in $X$, and {\it hyperbolic} if $\ell_X(\gamma)>0$. Then, if $(X,d)$ is complete (that will always be the case), there exists at least one geodesic $\Ax(\gamma)$ called a 
{\it translation axis} of $\gamma$ in $(X,d)$, which is invariant under $\gamma$ and such that $d(x,\gamma x)=\ell_X(\gamma)$ if 
$x\in\Ax(\ga)$ (see \cite[Chap.~2.6]{BriHae99}). Since the translation axes of a hyperbolic element are pairwise at finite Hausdorff distance,
according to \cite[Th.~2.13~p.182]{BriHae99}, if $X$ is a surface,
the union of all the translation axes of a hyperbolic element is a flat strip, possibly reduced to a single geodesic, or is isometric to $\RR\times[0,+\infty[$ or to $\RR^2$.

\medskip

Let $\srfce$ be a compact, connected surface endowed with a half-translation structure, and let $p:\revet\to\srfce$ be a universal cover with covering group $\Gamma$, as in Section \ref{structuresplates}. 
We will need the  Gauss-Bonnet formula for a half-translation structure on a surface. Let $P$ be a compact subsurface of $\srfce$ with piecewise geodesic boundary.
For every point $x\in\partial P$, the interior angle $\theta(x)>0$ is the flat angle between the two germs at $x$ of geodesic segments contained in $\partial P$, 
measured in the
angular sector inside of $P$. For every point $x\in\stackrel{\circ}{P}$, the angle $\theta(x)$ is the total flat angle at $x$ (i.e.
$n\pi$ with $n\in\NN$
and $n\geqslant 3$ 
if $x$ is a singular point, and $2\pi$ otherwise). 

\blemmA\cite[Prop.~3.6]{Dank10}\label{GaussBonnet} We have $2\pi\chi(P)=\underset{x\in\stackrel{\circ}{P}}\sum (2\pi-\theta(x)) + \underset{x\in\partial P}\sum (\pi-\theta(x))$.
\elemmA

Since the boundary of $P$ is piecewise geodesic, the angles $\theta(x)$ at the points $x\in\partial P$ are at least
$\pi$ except at the points where $\partial P$ is not locally geodesic. If we denote by $\theta_1,\dots,\theta_n$ the angles 
at these points,
the Gauss-Bonnet formula implies the following result.

\blemmA\label{GaussBonnet} We have $2\pi\chi(P)\leqslant\underset{x\in\stackrel{\circ}{P}}\sum (2\pi-\theta(x)) + \overset{n}{\underset{i=1}{\sum}} (\pi-\theta_i)$.
\elemmA

Let $\widetilde{\ell}$ be a geodesic of $\revet$.

\bdefiA The side $+$ and $-$ of $\widetilde{\ell}$ are the connected components of $\widetilde{\Sigma}-\widetilde{\ell}$.
\edefiA

Let $\ell$ be a periodic geodesic of $\srfce$ and let $\widetilde{\ell}$ be a lift of ${\ell}$ in $\widetilde{\Sigma}$.
      Let $\gamma\in\Gamma-\{e\}$ be an element of the stabilizer 
of $\widetilde{\ell}$, and let $x\in\widetilde{\ell}$. Then, the {\it total curvature of $[x,\gamma x[$ at side $+$} is the sum
$\sum_{t\in[x,\gamma x[}\;(\pi-\theta^+(t))$, where $\theta^+(t)$ is the angle of $\widetilde{\ell}$ at $t$ measured in the side $+$. 

\blemmA\label{Rafi}(see \cite[Rem.~3.2]{Raf05}) The total curvature of $[x,\gamma x[$ at side $+$ 
is an integral multiple of $\pi$. It is zero, or it is at most $-\pi$.
\elemmA

Let $\gamma\in\Ga$ be a hyperbolic element of translation length $\ell(\ga)$ and let $F(\gamma)$ be the (possibly degenerated) flat strip, union of all 
the translation axes of $\gamma$.
Let $x\in\widetilde{\Sigma}$ and let $x_\perp$ be the orthogonal projection of $x$ onto $F(\gamma)$. Let $\operatorname{Ax}(\gamma)$ be the translation axis of $\gamma$
in $\widetilde{\Sigma}$ that bounds $F(\gamma)$ and contains $x_\perp$.

\blemmA\label{pointfixe} The segment $[x,\gamma x]$ meets  $\operatorname{Ax}(\gamma)$ and   
$d(x,\ga\,x)^2\geqslant\ell(\gamma)^2 +2\,d(x,x_\perp)^2$.  
\elemmA 

\demA Let us consider the  geodesic quadrilateral with vertices $x$, $\ga x$, $\ga x_\perp$ and $x_\perp$. 
The two sides at the vertices $x$ and $\gamma x$ may 
share an initial segment. Let $[x,x_1]$ (resp. $[\gamma x,x_2]$)
be the intersection of the geodesic segments $[x,\gamma x]$ and $[x,x_\perp]$ (resp. $[\gamma x,x]$ and $[\gamma x,\gamma x_\perp]$).
Assume for a contradiction that the geodesic segment $[x_1,x_2]$ does not meet $[x_\perp,\gamma x_\perp]$. 
Then, the closed curve 
$[x_1,x_2]\cdot[x_2,\gamma x_\perp]\cdot[\gamma x_\perp,x_\perp]\cdot[x_\perp,x_1]$ is simple and bounds a topological disk $P$.  
For all $t\in [x_\perp,\gamma x_\perp]$, we denote by $\theta(t)$ the angle made by the rays of $\Ax(\gamma)$ starting at $t$, measured in the side of $\Ax(\gamma)$ containing $x$, and we denote by
$\theta_1,\,\theta_2$ the interior angles of $P$ at $x_1$ and $x_2$, and by $\theta_-$ and $\theta_+$ the interior angles of $P$ at $x_\perp$ 
and $\gamma x_\perp$.
The segment $[\gamma x,\gamma x_\perp]$ is the image of $[x,x_\perp]$ by $\gamma$, hence $\theta_-+\theta_+=\theta(x_\perp)$.
Since $\chi(P)=1$, according to the Gauss-Bonnet formula  (Lemma \ref{GaussBonnet}), we should have 
\begin{align*}
2\,\pi&\leqslant\sum_{z\in\stackrel{\circ}{P}}(2\pi-\theta(z))+\pi-\theta_1+\pi-\theta_2+\pi-\theta_-+\pi-\theta_++
\sum_{t\in]x_\perp,\gamma x_\perp[}(\pi-\theta(t))\\
&\leqslant 3\pi-(\theta_1+\theta_2)+\pi-(\theta_-+\theta_+)+
\sum_{t\in]x_\perp,\gamma x_\perp[}(\pi-\theta(t))\\
&<3\pi + \sum_{t\in[x_\perp,\gamma x_\perp[}(\pi-\theta(t)),
\end{align*}
\noindent
whereas, according to Lemma \ref{Rafi}, we have $\sum_{t\in[x_\perp,\gamma x_\perp[}(\pi-\theta(t))\leqslant-\pi$.
Hence, the segment $[x,\gamma x]$ meets $[x_\perp,\gamma x_\perp]$ at least at a point $y$. 
\begin{center}
	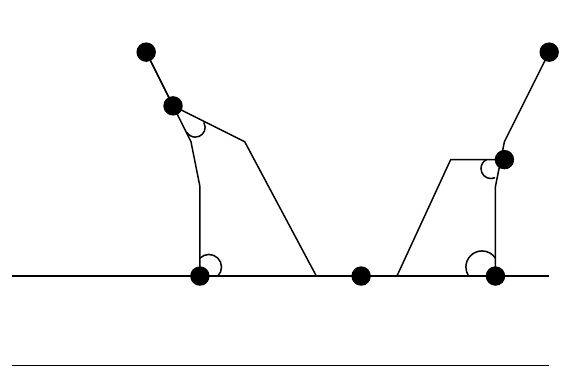
\end{center}
Since $\revet$ is $\CAT(0)$, we have $$d(x,y)^2\geqslant d(x,x_\perp)^2+d(x_\perp,y)^2 \mbox{ and } d(y,\gamma x)^2\geqslant d(\gamma x,\gamma x_\perp)^2+d(y,\gamma x_\perp)^2,$$
\noindent
hence $d(x,\ga\,x)^2\geqslant\ell(\gamma)^2 +2\,d(x,x_\perp)^2$.\cqfd

\medskip

Let $\Ax(\alpha)$ and $\Ax(\beta)$ be two 
translation axes of two hyperbolic elements $\alpha,\beta\in\Gamma$ that are not some powers of a common element. Then, 
their translation axes have no common point at infinity. Let $w$ be a point of $\Ax(\alpha)$ and $z$ be a point of $\Ax(\beta)$. We assume that 
$\Ax(\alpha)$ and $\Ax(\beta)$ are in the boundaries of the (possibly degenerated) flat strips union of all the translation axes of $\alpha$ and $\beta$, 
and that the segment $[w,z]$ does not meet a translation axis of $\alpha$ or $\beta$ other than $\Ax(\alpha)$ and $\Ax(\beta)$. Since $\revet$ is proper,
there exists $x\in\Ax(\alpha)$ and $y\in\Ax(\beta)$ such that $d(x,y)$ minimizes the distance between
$\Ax(\alpha)$ and
$\Ax(\beta)$. It may happen that $x=y$, then $[x,y]=\{x\}$. 

\blemmA\label{doublepointfixes} We have $d(z,w)\geqslant d(w,x)+d(z,y)+d(x,y)-2(\ell(\alpha)+\ell(\beta))$.
\elemmA

\demA We denote by $w'$ and $z'$ the endpoints of the intersection $[w,z]\cap\Ax(\alpha)$ and $[z,w]\cap\Ax(\beta)$ (possibly equal to $w$ and $z$). 
\begin{center}
	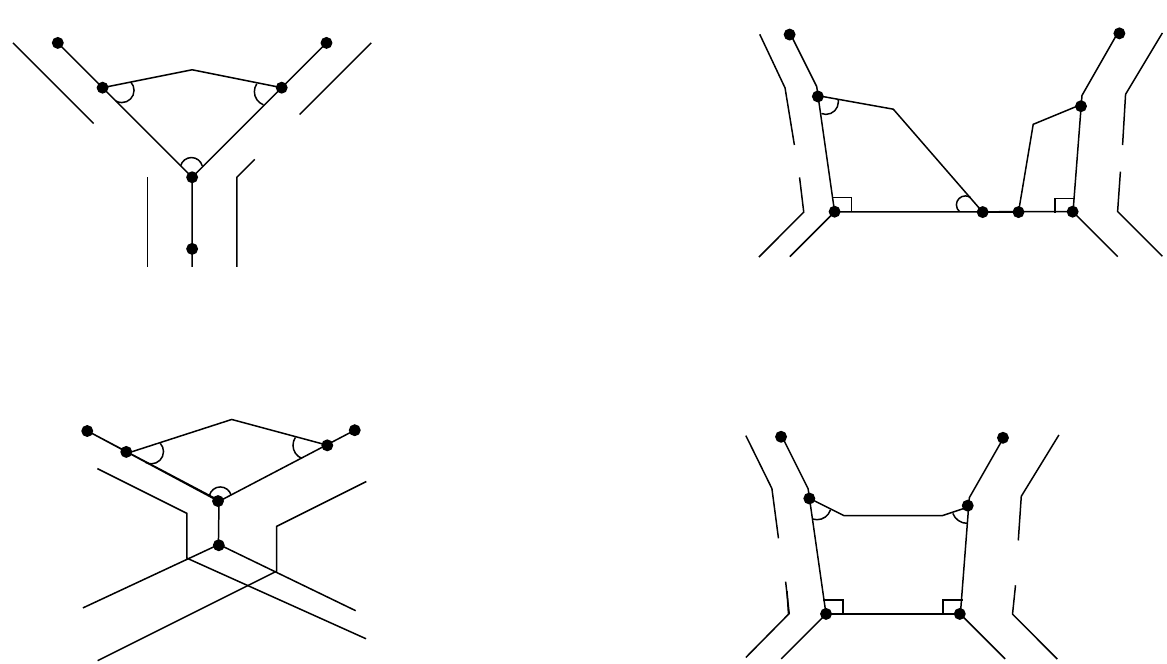
\end{center}
Assume first that $x=y$. Since $\Ax(\alpha)$ and $\Ax(\beta)$ are distinct and respectively $\alpha$ and $\beta$-invariant, they can share a segment of length at
most $\max\{\ell(\alpha),\ell(\beta)\}$. Let $t$ be the endpoint of the intersection $[x,w]\cap[x,z]$ (possibly equal to $x$). Note that $t\in[w',x]$.
If $w'=t$, then $z'=t$ and 
$d(z,w)=d(z,t)+d(t,w)\geqslant d(w,x)+d(z,x)-2\max\{\ell(\alpha),\ell(\beta)\}$, which proves the result. 
Hence, we assume for now on that $w'\not=z'$. Assume for a contradiction that $d(w,t)>\ell(\alpha)$ or
$d(z,t)>\ell(\beta)$. Let $T$ be the triangle with vertices $t$, $w'$ and $z'$.
According to Lemma \ref{Rafi}, by the assumption on $\Ax(\alpha)$, $\Ax(\beta)$, and on $[w,z]$, 
the sum of the total curvatures of $[t,w]$ and of $[t,z]$ at theirs sides containing $T$ is at most $-\pi$, and
according to the
Gauss-Bonnet formula applied to $T$, if we denote by $\theta_1$, $\theta_2$ and $\theta_3$ the interior angles of $T$ at $t$, $w'$ and 
$z'$, we should have $2\pi\leqslant 3\pi-(\theta_1+\theta_2+\theta_3)-\pi<2\pi$, a contradiction. Hence, we have $d(w',t)\leqslant\ell(\alpha)$ and
$d(z',t)\leqslant\ell(\beta)$, and
\begin{align*}
	d(z,w)&\geqslant d(z,z')+d(w,w')\geqslant d(w,x)-d(w',x)+d(z,x)-d(z',x)\\
	&\geqslant d(w,x)-(\ell(\alpha)+\max\{\ell(\alpha),\ell(\beta)\})+d(z,x)-(\ell(\beta)+\max\{\ell(\alpha),\ell(\beta)\}),
\end{align*}

\noindent
which proves the result. 

Assume next that $x\not=y$ and that $[w,z]$ meets $[x,y]$. Let $t_1$ and $t_2$ be the (possibly equal) endpoints of $[w,z]\cap[x,y]$, with $t_1$ the 
closest to $x$. If $t_1$ is distinct from $x$,
let $T$
be the triangle with vertices $w'$, $t_1$ and $x$. Assume for a contradiction that
$d(w',x)>\ell(\alpha)$. Then, by the assumption on $\Ax(\alpha)$, $\Ax(\beta)$, and since $]x,y[$ meets no translation axis of $\alpha$ nor $\beta$, and according to Lemma \ref{Rafi}, the total curvature of 
$[w',x]$ at the side of $\Ax(\alpha)$ containing $T$ is at most $-\pi$, and, if we denote by $\theta_1$ and $\theta_2$ the 
interior angles of $T$ at $w'$ and $t_1$, according to the Gauss-Bonnet formula applied to $T$, 
we should have $$2\pi\leqslant \pi-\theta_1+\pi-\theta_2+\pi-\frac{\pi}{2}-\pi<\frac{3\pi}{2},$$
a contradiction. Hence, we have $d(w',x)\leqslant\ell(\alpha)$ and since $x$ is the orthogonal projection of $t_1$ on $\Ax(\alpha)$, we have $d(w',t_1)\geqslant d(x,t_1)$. If 
$t_1=x$, then $w'=x=t_1$. Similarly, we have $d(z',t_2)\geqslant d(y,t_2)$ and $d(z',y)\leqslant\ell(\beta)$. 
Hence, we have, as wanted, 
\begin{align*}
	d(z,w)&=d(z,z')+d(z',t_2)+d(t_2,t_1)+d(t_1,w')+d(w',w)\\
	&\geqslant d(z,y)-\ell(\beta)+d(w,x)-\ell(\alpha)+d(x,y).
\end{align*}

Finally, assume that $[w,z]$ is disjoint from $[x,y]$. Then, by an argument similar to the one above, we have 
$d(z',y)\leqslant\ell(\alpha)$ and $d(w',x)\leqslant\ell(\beta)$. Hence, we have, as wanted
\begin{align*}
	d(z,w)&=d(z,z')+d(z',w')+d(w',w)\\
	&\geqslant d(z,y)-\ell(\beta)+d(w,x)-\ell(\alpha)+d(x,y).                                                                 
\end{align*}\cqfd

\section{Subsurface with geodesic boundary or filled up by simple closed curves}\label{tools} 

In this Section \ref{tools}, we assume that $\Sigma$ is a (possibly trivial) cover of a compact surface whose Euler characteristic is negative.
Let $p:\widetilde{\Sigma}\to\Sigma$ be a universal cover with covering group $\Gamma$. The results of Subsections \ref{filled} and \ref{subsurface} are 
folklore. We refer to \cite{Morzy4} for complete proofs.

\subsection{Subsurface filled up by a set of (free) homotopy classes of simple closed curves or by a measured hyperbolic lamination.}\label{filled}
Assume in this Subsection \ref{filled} that $\Sigma$ is compact. 
A {\it tight subsurface} of $\Sigma$ is a closed subsurface $\Sigma'$ of $\Sigma$ with smooth boundary, such that:

\medskip
\noindent
$\bullet$~ no connected component of $\Sigma'$ is a disk or a pair of pants;

\medskip
\noindent
$\bullet$~ no connected component of $\Sigma-\Sigma'$ is a disk.

\medskip
\noindent
$\bullet$~ no cylinder connected component of $\Sigma'$ can be homotoped in another connected component of $\Sigma'$.

If $C$ is a closed subset of $\Sigma$, 
an essential closed curve $\alpha$ of $\Sigma$ {\it topologically cuts} $C$ if it is not freely 
homotopic to a curve disjoint from $C$. Let $S$ be a subset of $\SSSS(\Sigma)$.

\blemmA (see \cite[Lem.~4.1]{Morzy4})\label{fills} There exists a (non unique) tight subsurface $W$ of $\Sigma$ such that:
\begin{itemize}
\item[$\bullet$] every essential closed curve of $\Sigma$ that topologically cuts $W$ has a positive intersection number with at least one element of $S$;
\item[$\bullet$] $W$ contains a  representative of each element of $S$.
\end{itemize}
\elemmA

\blemmA  \label{subfamily} (see \cite[Lem.~4.2]{Morzy4}) Let $S$ and $S'$ be two subsets of $\SSSS(\Sigma)$ such that $S'\subseteq S$. If $S$ and $S'$ fill up respectively $W$ and $W'$ then, 
up to isotopy, we have $W'\subseteq W$. 
\elemmA

\bcoroA(see \cite[Coro.~4.3]{Morzy4})\label{isotopy} Assume that $S$ fills up two subsurfaces $W$ and $W'$. Then $W$ and $W'$ are isotopic.   
\ecoroA

The subsurface $W$, which is unique up to isotopy, will be called the {\it subsurface filled up by $S$}.  It may  be non connected. Let $[q]$ be a half-translation structure and let $m$ be a (complete)
Riemannian
metric on $\Sigma$. If $S$ is finite and
$\{c[s]\}_{[s]\in S}$ is a set of $[q]$ or $m$-geodesic representatives of the elements of $S$, contained in $W$, then the complement of the union 
$\bigcup_{[s]\in S}c[s]$ in $W$ is a finite union of disks and (possibly) cylinders that can be homotoped to a boundary component of 
$W$. Otherwise, there would exist an essential closed curve $\alpha$ contained in a connected component of $\Sigma-(\bigcup_{[s]\in S}c[s])$,
and we would have $i(s,\alpha)=0$ for all $s\in
S$.


Let ${g}$ and $b$ be the genus and the number of boundary components of $W$. Let $K=(g+1)(2g+b)(3g+b)$ (this constant is not optimal).

\blemmA(see \cite[Lem.~4.4]{Morzy4})\label{FFini} Let $W$ be a tight subsurface filled up by a subset $S$ of $\SSSS(\Sigma)$. Then there exists a subset $S'\subseteq S$ with at most $K$ elements
that fills up $W$.  
\elemmA

Let $m$ be a hyperbolic metric on $\Sigma$.
A {\it measured hyperbolic lamination} $(\Lambda,\mu)$  of $\srfcem$ is a closed subset of $\Sigma$ which 
is a disjoint union of simple geodesics, called {\it leaves}, endowed with a transverse measure. We refer for instance to \cite{Bonahon97}
for the definition of a transverse measure and the definition of the intersection number with $(\Lambda,\mu)$.

\blemmA\label{remplieparunelamination}(see \cite[Lem.~4.5]{Morzy4}) There exists a unique tight subsurface $W(\Lambda)$ of $\Sigma$ (up to isotopy)
such that $\Lambda$ is contained in $W$ and for every essential closed curve $\beta$ of $\Sigma$ that 
topologically cuts $W$, we have $i((\Lambda,\mu),\beta)>0$.
\elemmA

The subsurface $W(\Lambda)$ is said to be {\it filled up} by $(\Lambda,\mu)$. It is connected if and only if $\Lambda$ is minimal.

\subsection{Subsurfaces with flat geodesic boundary.}\label{subsurface}

We now discuss some geometric properties of subsurfaces filled up by some set of (homotopy classes of) simple closed curves.  
Let $[q]$ be a half-translation structure on $\Sigma$. Let $D$ be
an open disk embedded in $\Sigma$ with piecewise smooth boundary. The {\it perimeter} $\peri(D)$ of $D$ is the length
of the boundary of the completion of $D$ for the induced metric. 

\blemmA\label{disk}(see \cite[Lem.~4.6]{Morzy4}) The diameter of $D$ is at most $\peri(D)$ and the area of $D$ 
is at most $\frac{1}{4\pi}\peri(D)^2$. 
\elemmA

If $\alpha$ is an isotopy class of simple closed curves, let $C(\alpha)$ be the 
(possibly degenerated) flat cylinder union of all the geodesic representatives of $\alpha$. 
If $c$ is a piecewise smooth simple closed curve in $\alpha$, disjoint from ${C}(\alpha)$, 
let $b(\alpha,c)$ be the boundary component of $C(\alpha)$ such that $b(\alpha,c)$ and $c$ bound a cylinder disjoint from the interior of  ${C}(\alpha)$
($b(\alpha,c)$ is unique since $\Sigma$ is not a torus). 

\blemmA\label{cylinder}(see \cite[Lem.~4.7]{Morzy4}) We have $d(c,b(\alpha,c))\leqslant\operatorname{length}(c)$. 
\elemmA

The following notion of geometric realization of surfaces is due to \cite[§2]{Raf05}. 
Let $W$ be a non trivial $\pi_1$-injective connected subsurface of $\Sigma$. 
The fundamental group $\pi_1(W)$ (with respect to any choice of basepoint in $W$)
is a subgroup of $\Gamma$, that is determined by the homotopy class of $W$, up to conjugation. 

\medskip

Let $\widehat{p}:\Sigma_W=\widetilde{\Sigma}/\pi_1(W)\to\Sigma$ be the $W$-cover of $\Sigma$, and let
$p':\widetilde{\Sigma}\to\Sigma_W$ be the unique universal cover such that $p=\widehat{p}\circ p'$. Let $[\widehat{q}]$ and $[\widetilde{q}]$ be the pullbacks of $[q]$ on $\Sigma_W$
and on $\widetilde{\Sigma}$. 
The surface $\Sigma_W$ is homotopically equivalent to $W$, and the preimage $\widehat{p}^{-1}(W)$ has a unique connected component
$\widehat{W}$ that is not simply connected. The map $\widehat{p}_{|\widehat{W}}:\widehat{W}\to W$ is a homeomorphism and the complement 
of $\widehat{W}$
in $\Sigma_W$ is a finite union of open annuli that can be homotoped to boundary components of $\widehat{W}$. A {\it generalized subsurface} is a closed
connected union of a (possibly empty) subsurface with some finite connected (metric) graphs, glued at some points of the boundary of the subsurface. Its {\it boundary} is 
the union of the boundary of this subsurface with these graphs. 
If $\widehat{W}$ is not a cylinder, 
the {\it $[\widehat{q}]$-geometric realization} of $\widehat{W}$ is the unique generalized subsurface $\widehat{W}_{[\widehat{q}]}$ of ${\Sigma}_W$
homotopic to $\widehat{W}$ within $\Sigma_W$,
whose boundary is the union of
some $[\widehat{q}]$-geodesic representatives of the boundary components of $\widehat{W}$ (see \cite[p.~188]{Raf05}),
that contains a unique $[\widehat{q}]$-geodesic representative of each boundary component of $\widehat{W}$. We call 
{\it boundary components of $\widehat{W}_{[\widehat{q}]}$} the $[\widehat{q}]$-geodesic representatives of the boundary components of $\widehat{W}$ contained in
$\widehat{W}_{[\widehat{q}]}$, and we say that a boundary component of $\widehat{W}_{[\widehat{q}]}$ {\it corresponds to} a boundary component of $\widehat{W}$
if they are freely homotopic. If $\widehat{W}$ is a cylinder, the {\it $[\widehat{q}]$-geometric realization} of $\widehat{W}$ is the (possibly degenerated)
flat cylinder, union of all the $[\widehat{q}]$-geodesic representatives of the boundary components of $\widehat{W}$.
For every  essential homotopy class of closed 
curves $\alpha$ of 
$\widehat{W}$, the $[\widehat{q}]$-geodesic representatives of $\alpha$ are contained in $\widehat{W}_{[\widehat{q}]}$. 

\medskip

Similarly, if $\widetilde{W}$ is a connected component of the preimage of ${W}$ in $\widetilde{\Sigma}$, the $[\widetilde{q}]$-geodesics having the same pair of
points at infinity than the boundary components of $\widetilde{W}$
may not be pairwise disjoint. However, there exists a (possibly non unique) generalized subsurface properly homotopic to $\widetilde{W}$ within 
$\widetilde{\Sigma}$ whose boundary is the union of some geodesics having the same pairs of points at infinity than the boundary components of $\widetilde{W}$. 
If $\widetilde{W}$ is not a strip, we call {\it $[\widetilde{q}]$-geometric realization} of $\widetilde{W}$ the unique such generalized subsurface 
$\widetilde{W}_{[\widetilde{q}]}$ that contains a unique $[\widetilde{q}]$-geodesic
representative of each boundary component of $\widetilde{W}$, and we call 
{\it boundary components of $\widetilde{W}_{[\widetilde{q}]}$} the $[\widetilde{q}]$-geodesic representatives of the boundary components of $\widetilde{W}$ contained in
$\widetilde{W}_{[\widetilde{q}]}$. We say   that a boundary component of $\widetilde{W}_{[\widetilde{q}]}$ {\it corresponds to} a boundary component of $\widetilde{W}$
if they have the same ordered pair of points at infinity. 
If $\widetilde{W}$ is a strip, we call {\it $[\widetilde{q}]$-geometric realization} of $\widetilde{W}$ the (possibly degenerated) flat strip, union of all the
$[\widetilde{q}]$-geodesics having the same (unordered) pair of points at infinity than the boundary components of $\widetilde{W}$.  

\medskip

We denote by $\Sing(\widetilde{W}_{[\widetilde{q}]})$ the set of singular points of $[\widetilde{q}]$ contained in $\widetilde{W}_{[\widetilde{q}]}$.

\blemmA\label{hyperbolicity}
The space $\widetilde{W}_{[\widetilde{q}]}$ is convex. Moreover, for all $\varepsilon>0$, if the union $\Sing(\widetilde{W}_{[\widetilde{q}]})\cup\partial\widetilde{W}_{[\widetilde{q}]}$
is $\varepsilon$-dense into $\widetilde{W}_{[\widetilde{q}]}$, then $\widetilde{W}_{[\widetilde{q}]}$ is $2\varepsilon$-hyperbolic. 
\elemmA

\demA The proof is 
essentialy the same as the one of \cite[Prop.~3.7]{Dank10}. It suffices to replace $\rho=\sup_{x\in \Sigma}\, d(x,\Sing(\Sigma))$ by
$\rho=
\sup_{x\in \widetilde{W}_{[\widetilde{q}]}}
d(x,\Sing(\widetilde{W}_{[\widetilde{q}]})\cup\partial \widetilde{W}_{[\widetilde{q}]})$.\cqfd

\medskip

\remA If the diameter $\operatorname{Diam}\srfce$ is finite, the universal cover $\revet$ is $2\operatorname{Diam}\srfce$-hyperbolic. 

\medskip


Assume that $\widehat{W}$ is neither a cylinder nor a pair of pants, and that $\widehat{W}$ is filled up by a finite set of isotopy classes of simple closed curves 
$S(\widehat{W})$ of cardinality at most $K$. We can always assume 
that $S(\widehat{W})$ does not contain the homotopy class
of any boundary component of $\widehat{W}$. Let $\varepsilon=
\max_{\alpha\in S(\widehat{W})}\ell_{[\widehat{q}]}(\alpha)$ and let $b$ be the number of boundary components of $\widehat{W}$.  

\blemmA\label{diameterandarea}(see \cite[Lem.~4.9]{Morzy4})The length of any boundary component of $\widehat{W}_{[\widehat{q}]}$ is at most $K\varepsilon$, the
diameter of $\widehat{W}_{[\widehat{q}]}$ is at most $11K\varepsilon$ and its area is at most
$\frac{1}{\pi}(1+2b)^2(K\varepsilon)^2$. 
\elemmA

\section{Ultralimits of  sequences of metric spaces.}\label{ultralimite}

The ultralimit of a sequence of metric spaces is a notion introduced by M.~Gromov (see \cite{Gromov93}, we refer for example to \cite{Drutu01} for background 
and precisions on the content of this Section \ref{ultralimite}). It uses ultrafilters (introduced by H. Cartan) that are a
way of picking an accumulation value of a sequence in a compact metrizable space, avoiding extraction arguments.
Let $\omega$ be a non principal ultrafilter on $\NN$
(see \cite[§~6.4]{Bourbaki71a}). We say that a sentence $A(n)$ is true for {\it $\omega$-almost all $n\in\NN$} if there exists
$I\in\omega$ such that $A(n)$ is true for all 
$n\in I$. If $(a_n)_{n\in\NN}$ is a sequence of a topological space 
$E$ and $a\in E$, we say that $(a_n)_{n\in\NN}$ {\it $\omega$-converges to $a$}, and write $\lim_\omega a_n=a$,   
if for every neighborhood $V$ of $a$, the element $a_n$ belongs to $V$ for $\omega$-almost all $n\in\NN$.
The
$\lim_\omega$ satisfies the usual properties of limits. If $E$ is compact, for example if $E=[-\infty,+\infty]$, for every non principal ultrafilter $\omega$ and 
for every sequence $(a_n)_{n\in\NN}$ in $E$, there exists a unique  accumulation value $a$ of $(a_n)_{n\in\NN}$ in $E$ such that $\lim_\omega a_n=a$.

Let $(X_n,d_n,\star_n)_{n\in\NN}$ be a
sequence of pointed metric spaces endowed with an isometric action of a group $\Gamma$.
Let  $$X'_\omega=\{(x_n)_{n\in\NN}\in\underset{n\in\NN}\prod X_n\,: \lim_\omega d_n(x_n,\star_n)<+\infty\}.$$
Then, the function 
$d'_\omega:X'_\omega\times X'_\omega\to\RR^+$ defined by
$$d'_\omega((x_n)_{n\in\NN},(y_n)_{n\in\NN})=\lim_\omega d_n(x_n,y_n)$$ is a pseudo-distance on
$X'_\omega$. 
We denote by $X_\omega$ the quotient of $X'_\omega$ by the equivalence relation 
$(x_n)_{n\in\NN}\sim(y_n)_{n\in\NN}$ 
if $d'_\omega((x_n)_{n\in\NN},(y_n)_{n\in\NN})=0$,
and the equivalence class of an element $(x_n)_{n\in\NN}$ of $X'_\omega$ is denoted by 
$[x_n]_{n\in\NN}$.
Then, the pseudo-distance $d'_\omega$ induces a distance $d_\omega$ on $X_\omega$. If for $\omega$-almost all $n\in\NN$, $F_n$ is a subset of 
$X_n$, let $[F_n]_{n\in\NN}=\{[x_n]_{n\in\NN}\,:\,x_n\in F_n\,\mbox{ for } \omega\mbox{-almost all }\,n\in\NN\}$.
The diagonal action of $\Gamma$ on $\underset{n\in\NN}\prod X_n$ is said to be 
{\it admissible} if $\lim_\omega d_n(\star_n,\gamma\star_n)<+\infty$ for every $\gamma\in\Gamma$. Then, the action $\Gamma\times X_\omega\to X_\omega$ defined by 
$\gamma [x_n]_{n\in\NN}=[\gamma x_n]_{n\in\NN}$ is an
isometric action on $X_\omega$. Let $\star_\omega=[\star_n]_{n\in\NN}$. If the action of $\Gamma$ is admissible, the pointed metric space 
$(X_\omega,d_\omega,\star_\omega)$ is called the {\it ultralimit
	of $(X_n,d_n,\star_n)_{n\in\NN}$ for
	$\omega$} and is denoted by $\lim_\omega(X_n,d_n,\star_n)$. We recall a few properties of ultralimits of metric spaces: 

\begin{itemize}
	\item[$\bullet$] If for all $n\in\NN$, $(X_n,d_n,\star_n)$ is a geodesic metric space, so is $(X_\omega,d_\omega,\star_\omega)$;
	\item[$\bullet$] If for all $n\in\NN$, $(X_n,d_n,\star_n)$ is $\CAT(0)$, so is $(X_\omega,d_\omega,\star_\omega)$;
	\item[$\bullet$] If for all $n\in\NN$, $(X_n,d_n,\star_n)$ is $\delta_n$-hyperbolic in the sense of Gromov, and $\lim_\omega\delta_n=\delta$,
	then $(X_\omega,d_\omega,\star_\omega)$ is $\delta$-hyperbolic.
\end{itemize}

Let $(X_n,d_n,\star_n)_{n\in\NN}$ be a sequence of $\CAT(0)$ pointed metric spaces with an admissible isometric
action of $\Gamma$, and let 
$\gamma\in\Gamma$. If $\gamma$ is hyperbolic in $X_n$, 
and $\Ax_{X_n}(\gamma)$ is a translation axis of $\gamma$ in $X_n$ for $\omega$-almost all $n\in\NN$, and if $\lim_\omega d_n(\star_n,\Ax_{X_n}(\gamma))<+\infty$,
then
$[\Ax_{X_n}(\gamma)]_{n\in\NN}$ 
is a translation axis of $\gamma$ in $X_\omega$ if 
$\ell_{X_\omega}(\gamma)>0$, and it is contained in the set of fixed points of $\gamma$ in $X_\omega$ if $\ell_{X_\omega}(\gamma)=0$.
If $\gamma$ is elliptic in $X_n$, 
and $x_n(\gamma)$ is a fixed point of $\gamma$ in $X_n$ for $\omega$-almost all $n\in\NN$, and if $\lim_\omega d_n(\star_n,x_n(\gamma))<+\infty$,
then
$[x_{n}(\gamma)]_{n\in\NN}$ is a fixed point of $\gamma$ in $X_\omega$. In the two cases, we have $\ell_{X_\omega}(\gamma)=\lim_\omega\ell_{X_n}(\gamma)$. 

\medskip

\remA The ultralimit $X_\omega$ does not change if we choose another sequence of base points $(\star_n')_{n\in\NN}$, as long
as $\lim_\omega d_n(\star_n,\star_n')<+\infty$.
The choice of the base points will not always be specified.

\section{Ultralimits of sequences of half-translation structures.}\label{sequenceofhalftranslations}

Let $\Sigma$ be a compact, connected, orientable surface such that $\chi(\Sigma)<0$, and let $p:\widetilde{\Sigma}\to\Sigma$ be a universal cover with covering group
$\Gamma$. Recall that the boundary of $\Sigma$ is empty, to simplify the writting of the article, but the results can be extended very easily to a 
surface with nonempty boundary (see \cite{Morzy4}). Let 
$\Flat(\Sigma)$ be the set of isotopy classes of half-translation structures on $\Sigma$.
   Let $([q'_n])_{n\in\NN}$ be a sequence in $\Flat(\Sigma)$ and,  for all $n\in\NN$, let $[\widetilde{q}'_n]$ be
the pullback of $[q'_n]$ on $\widetilde{\Sigma}$. Let $\omega$ be a non principal ultrafilter on $\NN$ as 
in Section \ref{ultralimite}. 

\medskip

Let $S$ be a finite generating set of $\Gamma$. For all $n\in\NN$,
define $f_n:\widetilde{\Sigma}\to\RR^+$ by 
$f_n(x)=\max_{s\in S} d'_n(x,s x)$. Let
$\lambda_n=\inf_{x\in\widetilde{\Sigma}} f_n(x)$ and let us choose a point  $\star_n\in\widetilde{\Sigma}$  such
that 
$f_n(\star_n)\leqslant\lambda_n+1$. Let $[\widetilde{q}_n]=\frac{1}{\lambda_n}[\widetilde{q}'_n]$, and let $d_n$  be the distance defined by $[\widetilde{q}_n]$ on $\widetilde{\Sigma}$.

For all $\gamma\in\Gamma$, we have  $\lim_\omega d_n(\star_n,\gamma\star_n)<+\infty$, hence
$\lim_\omega(\widetilde{\Sigma},[\widetilde{q}_n],\star_n)$ is endowed with an
isometric action of $\Gamma$. Moreover, if $[x_n]_{n\in\NN}\in\lim_\omega(\widetilde{\Sigma},[\widetilde{q}_n],\star_n)$, 
there exists $s\in S$ such that $d_n(x_n,s x_n)\geqslant 1$ for $\omega$-almost all $n\in\NN$, hence this action has no global fixed point.

\medskip

\bremaA\label{convergencedistancetranslation} For all $\gamma\in\Gamma-\{e\}$ and $n\in\NN$, according to Lemma \ref{pointfixe}, if $F_n(\gamma)$ is the flat strip, union of all the translation axes of $\gamma$ in $\revetn$, we have 
$d_n(\star_n,F_n(\gamma))\leqslant d_n(\star_n,\gamma\star_n)$. Since the action of   $\Gamma$ on $\revetn_{n\in\NN}$ is admissible, 
we have $\lim_\omega d_n(\star_n,F_n(\gamma))<+\infty$. 
Hence, the ultralimit $[F_n(\gamma)]_{n\in\NN}$ exists, and $\ell_{\lim_\omega\revetn}(\gamma)=\lim_\omega\ell_{\revetn}(\gamma)$.
Moreover, according to Lemma \ref{pointfixe}, the set $[F_n(\gamma)]_{n\in\NN}$ is exactly the union of all the translation axes or the set of fixed
points of $\gamma$ in $\lim_\omega\revetn$. 
Moreover, according to Lemma \ref{doublepointfixes}, if $\gamma_1$ and $\gamma_2$ are not powers of a common element, 
then they have at most one common fixed point in $\lim_\omega\revetn$.
\eremaA

We will see (Section \ref{nondegeneration}) that the ultralimit
$\lim_\omega\revetn$ is a surface endowed with a half-translation structure if and only if there exists $\varepsilon>0$ such that for every $\alpha\in\SSSS(\Sigma)$
and for $\omega$-almost all $n\in\NN$, we have $\ell_{[q_n]}(\alpha)\geqslant\varepsilon$. In  Section \ref{areazero}, we consider the case where there does not exist such a 
uniform lower bound on the $[q_n]$-lengths of the elements of $\SSSS(\Sigma)$. We will notably consider the case where there exist subsurfaces with fixed homotopy type
whose boundary components (if any) have their lengths that $\omega$-converge to $0$, and whose areas $\omega$-converge to $0$.
Let us introduce some general definition. 
Let $m$ be a hyperbolic metric on $\Sigma$, and let $\widetilde{m}$ be its pullback on $\widetilde{\Sigma}$.

For $i\in\{1,2,3\}$, let $(x_i,y_i)$ be an element of $\dddp$ and let $\widetilde{\lambda}_i$ be the geodesic of $\revetm$ whose ordered pair of points at infinity is $(x_i,y_i)$. We assume 
that the corresponding unordered pairs of points are pairwise distinct.
The pairs $(x_1,y_1)$ and $(x_2,y_2)$ are {\it interlaced} if the geodesics $\widetilde{\lambda}_1$ and $\widetilde{\lambda}_2$ intersect each other.
If  $(x_1,y_1)$ and $(x_2,y_2)$ are not interlaced, the pair $(x_3,y_3)$ is {\it contained} between $(x_1,y_1)$ and $(x_2,y_2)$ if 
$\widetilde{\lambda}_3$ is contained in the connected component of $\widetilde{\Sigma}-\widetilde{\lambda}_1\cup\widetilde{\lambda}_2$ bounded by $\widetilde{\lambda}_1$ and $\widetilde{\lambda}_2$, and the pair $(x_3,y_3)$ is
{\it caught} between $(x_1,y_1)$ and $(x_2,y_2)$ if it is contained between $(x_1,y_1)$ and $(x_2,y_2)$ and if $\widetilde{\lambda}_3$ intersects any geodesic segment 
joining $\widetilde{\lambda}_1$ and $\widetilde{\lambda}_2$. These definitions do not depend on the choice of $m$.  

If $X$ is an oriented geodesic of $\widetilde{\Sigma}$ (for any complete, $\CAT(0)$ metric on $\widetilde{\Sigma}$) or an element of $\Gamma-\{e\}$, we denote by 
$E(X)\in\dddp$ its ordered pair
of points at infinity or of fixed points at infinity. If $X,\,Y$ and $Z$ are three such elements, we say that $X$ and $Y$ are {\it interlaced} if $E(X)$ and $E(Y)$ are interlaced and that $Z$ is
{\it caught} or {\it contained} between $X$ and $Y$ if $E(Z)$ is caught or contained between $E(X)$ and $E(Y)$.

\medskip

In the remainder of this section, we will consider a tight, connected subsurface $W$ of $\Sigma$
(possibly equal to $\Sigma$).
Let $\widehat{p}:{\Sigma}_W\to\Sigma$ be a 
$W$-cover of $\Sigma$ and let $\widehat{W}$ be the unique connected component of the preimage of $W$ in $\Sigma_W$ which is not simply connected. 
Let $\widetilde{W}$ be a connected component of the preimage of $W$ in $\widetilde{\Sigma}$.
Let $([\widehat{q}_n])_{n\in\NN}$ and $([\widetilde{q}_n])_{n\in\NN}$ be the pullbacks of $([q_n])_{n\in\NN}$ on 
${\Sigma}_W$ and on $\widetilde{\Sigma}$. For all $n\in\NN$, we denote by $\widehat{W}_n$ and $\widetilde{W}_n$ the $[\widehat{q}_n]$ and $[\widetilde{q}_n]$-geometric 
realizations of $\widehat{W}$ and of $\widetilde{W}$.

\subsection{Typical degenerations.}\label{areazero} 

In this subsection, we assume that:

\medskip
\noindent
$\bullet$ $\lim_\omega\ell_{[q_n]}(\alpha)=\epsilon_\alpha<+\infty$, for every $\alpha\in\C(\Sigma)$;

\medskip
\noindent
$\bullet$ $\lim_\omega\ell_{[{q}_n]}({b})=0$, for every boundary component ${b}$ of ${W}$. 
\medskip 

Let $\widetilde{c}$ be a boundary component of $\widetilde{W}$.
We denote by $\Gamma_{\widetilde{c}}$ the stabilizer of $\widetilde{c}$ in $\Gamma$ and by $\gamma_{\widetilde{c}}$ a primitive generator of
$\Gamma_{\widetilde{c}}$. For all $n$,
let $F_n(\widetilde{c})$ be the (possibly degenerated) flat strip, union of all the geodesic representatives of $\widetilde{c}$.
Let $\gamma$ be a primitive element of $\Gamma-\{e\}$ whose translation axes in any $\revetn$ are interlaced with
$\widetilde{c}$, and let $(\Ax_n(\gamma))_{n\in\NN}$ be a sequence of translation axes of $\gamma$ in $\revetn$. 

\blemmA\label{star} The geodesic $\Ax_\omega(\gamma)=[\Ax_n(\gamma)]_{n\in\NN}$ is the only translation axis or set of fixed points of $\gamma$ 
in $\lim_\omega\revetn$, and the intersection
$\Ax_\omega(\gamma)\cap[F_n(\widetilde{c})]_{n\in\NN}$ is a geodesic segment orthogonal to the boundary of $[F_n(\widetilde{c})]_{n\in\NN}$ if $[F_n(\widetilde{c})]_{n\in\NN}$
is not reduced to a geodesic line, and is reduced to a point  otherwise. Moreover, if $\gamma'$ is another primitive element of $\Gamma-\{e\}$ whose translation axes in any 
$\revetn$ are interlaced with 
$\widetilde{c}$, and if $\Ax_\omega(\gamma')$ is the set of fixed points or the translation axis of $\gamma'$ in $\lim_\omega\revetn$, then
$\Ax_\omega(\gamma')\cap[F_n(\widetilde{c})]_{n\in\NN}=\Ax_\omega(\gamma)\cap[F_n(\widetilde{c})]_{n\in\NN}$.
\elemmA

We denote by $\star_{\widetilde{c}}$ the unique point of $\Ax_\omega(\gamma)\cap[F_n(\widetilde{c})]_{n\in\NN}$ that belongs to $[\widetilde{W}_n]_{n\in\NN}$.
We will need the following lemma (which is not surprising for specialist) in order to prove Lemma \ref{star}. Let $c$ be the image of $\widetilde{c}$ in $\Sigma$. Let $\alpha_\gamma$ be the free homotopy class of closed curves defined by $\gamma$ in $\Sigma$. 
Let $S$ be a filling finite set of simple closed curves in $\Sigma$, such that $c$ does not belong to $S$. For all $n\in\NN$, we denote by $I_n$ and $J_n$ the geodesic
segments in $F_n(\widetilde{c})$ perpendicular to its boundary components $\widetilde{c}_n$ and $\widetilde{c}_n'$
(reduced to two points if $F_n(\widetilde{c})$ is reduced to a single geodesic) such that one endpoint of $I_n$ (resp. $J_n$) is the first (resp. last) intersection point between 
$\Ax_n(\gamma)$ and $F_n(\widetilde{c})$.  
\begin{center}
	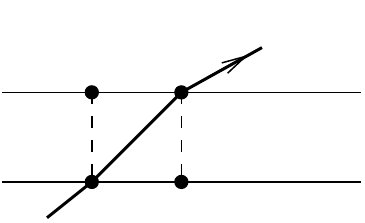
\end{center}

\blemmA\label{confonduborne}(see \cite[Lem.~6.3]{Morzy4}) We have $d_n(I_n,J_n)\leqslant\ell_{[{q}_n]}({c})
\sum_{s\in S}i(\alpha_\gamma,s)$.
\elemmA

{\noindent{\bf Proof of Lemma \ref{star}. }} According to Remark \ref{convergencedistancetranslation} the ultralimit
$[\Ax_n(\gamma)]_{n\in\NN}$ exists. 

Let us first prove that the intersection  $\Ax_\omega(\gamma)\cap[F_n(\widetilde{c})]_{n\in\NN}$ is a geodesic segment orthogonal to the boundary components
of $[F_n(\widetilde{c})]_{n\in\NN}$, if $[F_n(\widetilde{c})]_{n\in\NN}$ is a flat strip, and is reduced to a point otherwise.  

Let $I_n$ and $J_n$ be the segments (possibly reduced to a point) orthogonal to $F_n(\widetilde{c})$, having as an endpoint the first and the last intersection point between
$\Ax_n(\gamma)$ and $F_n(\widetilde{c})$, as in Lemma \ref{confonduborne}. Since $\lim_\omega\ell_{[q_n]}(c)=0$, according to Lemma \ref{confonduborne},  
we have $\lim_\omega d_n(I_n,J_n)=0$. 
For all $n\in\NN$, let  $x_n\in J_n$ be the last intersection point between $F_n(\widetilde{c})$ and the oriented geodesic $\Ax_n(\gamma)$. Let 
$(z_n)_{n\in\NN}$ be a sequence of points such that, for all $n$, $z_n$ belongs to $[x_n,\Ax_n(\gamma)(+\infty)[$, and let $z_{n\perp}$ be the orthogonal projection of
$z_n$ onto 
$F_n(\widetilde{c})$. Let $\widetilde{c}_n$ be the boundary component of $F_n(\widetilde{c})$ containing $z_{n\perp}$ and let $+$ be the side of $\widetilde{c}_n$
containing $z_n$. 
The interior of $F_n(\widetilde{c})$, if not empty, is not contained in the side $+$ of $\widetilde{c}_n$. 
Hence, according to Lemma \ref{Rafi}, the total curvature of a fundamental domain of $\widetilde{c}_n$ for the action of $\gamma_{\widetilde{c}}^{\ZZ}$, at the side $+$,
is at most $-\pi$. 
Let us prove that $d_n(x_n,z_{n\perp})\leqslant\ell_{[q_n]}(c)$. Otherwise, $z_{n\perp}$ is not equal to $x_n$ and the open segment $]x_n,z_{n\perp}[$ 
fully contains an entire fundamental domain of 
$\widetilde{c}_n$ for the action of $\gamma_{\widetilde{c}}^{\ZZ}$. Let $a_n$ be the last intersection point between $[z_n,x_n]$ and $[z_n,z_{n\perp}]$. 
The curve $[x_n,a_n]\cdot[a_n,z_{n\perp}]\cdot[z_{n\perp},x_n]$ bounds a topological disk $T$.

\begin{center}
	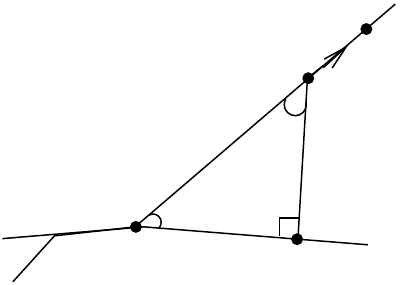
\end{center}

The interior angle
at $z_{n\perp}$ is at least $\frac{\pi}{2}$. Let $\theta_1$ and $\theta_2$ be the two other interior angles of the boundary of $T$ which are less than $\pi$. By
the Gauss-Bonnet formula applied to $T$, we should have
$$2\pi\chi(T)\leqslant \pi-\frac{\pi}{2}+\pi-\theta_1+\pi-\theta_2-\pi\leqslant\frac{3}{2}\pi,$$
\noindent
a contradiction since $\chi(T)=1$. Hence, we have $d_n(x_n,z_{n\perp})\leqslant\ell_{[{q}_n]}({c})$, and $d_n(z_n,{F}_n(\widetilde{c}))=d_n(z_n,z_{n\perp})\geqslant d_n(z_n,x_n)
-\ell_{[{q}_n]}({c})$. Since $\lim_\omega\ell_{[{q}_n]}({c})=0$, we deduce that $[z_n]_{n\in\NN}$ belongs to $[F_n(\widetilde{c})]_{n\in\NN}$
if and only if $[z_n]_{n\in\NN}=[x_n]_{n\in\NN}$. Similarly, for all $n\in\NN$, let $y_n$ be the first intersection point between $\Ax_n(\gamma)$ and $F_n(\widetilde{c})$, 
and let $z_n$ belongs to $[y_n,\Ax_n(\gamma)(-\infty)[$, then $[z_n]_{n\in\NN}$  belongs to $[F_n(\widetilde{c})]_{n\in\NN}$ if and only if $[z_n]_{n\in\NN}=[y_n]_{n\in\NN}$.
Moreover, since the distance between $I_n$ and $J_n$ $\omega$-converges to zero, the intersection between 
$[\Ax_n(\gamma)]_{n\in\NN}$ and $[F_n(\widetilde{c})]_{n\in\NN}$ is a segment orthogonal to the boundary components of $[F_n(\widetilde{c})]_{n\in\NN}$, possibly reduced to a point if $[F_n(\widetilde{c})]_{n\in\NN}$
is reduced to a single geodesic.  
\begin{center}
	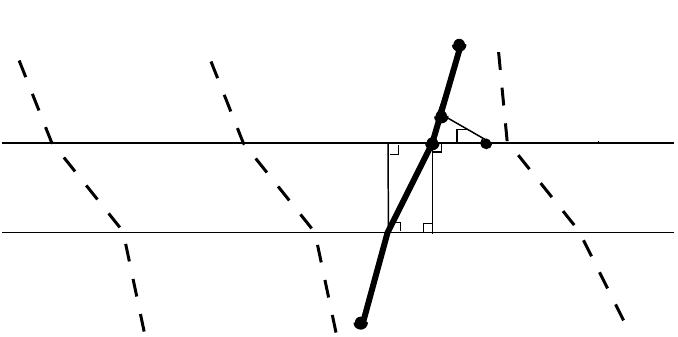
\end{center}

Let $(\widetilde{\beta}_n)_{n\in\NN}$ be a sequence of translation axes of another element of $\Gamma-\{e\}$ which is interlaced with $\widetilde{c}$,
or another sequence of translation axes of $\gamma$, and let $n\in\NN$. Since the action of $\gamma_{\widetilde{c}}^{\ZZ}$ on 
$\partial_\infty\widetilde{\Sigma}$ is a South-North dynamic with fixed points $\widetilde{c}(+\infty)$ and $\widetilde{c}(-\infty)$, there exists $k\in\NN$ such that 
$\Ax_n(\gamma)$ and $\widetilde{\beta}_n$ are caught between $\gamma_{\widetilde{c}}^k(\widetilde{\beta}_n)$ and $\gamma_{\widetilde{c}}^{-k}(\widetilde{\beta}_n)$. Moreover, the integer $k$
only depends on $E(\Ax_n(\gamma))$ and $E(\widetilde{\beta}_n)$, hence does not depend on $n\in\NN$. Hence,
we have $\lim_\omega d_n(\Ax_n(\gamma)\cap F_n(\widetilde{c}),\widetilde{\beta}_n\cap F_n(\widetilde{c}))\leqslant\lim_\omega
2k\,\ell_{[{q}_n]}({c})=0$, and $[\Ax_n(\gamma)]_{n\in\NN}\cap[F_n(\widetilde{c})]_{n\in\NN}=[\widetilde{\beta}_n]_{n\in\NN}\cap[F_n(\widetilde{c})]_{n\in\NN}$. 

According to Remark \ref{convergencedistancetranslation}, any translation axis or any set of fixed points of an element of $\Gamma-\{e\}$ in $\lim_\omega\revetn$ is the ultralimit of a sequence
of translation axes of $\gamma$ in $\revetn$. Hence, the element $\gamma$ has no other translation axis or set of fixed points than $\Ax_\omega(\gamma)$, and if 
$\gamma'\in\Gamma-\{e\}$ is interlaced with $\widetilde{c}$, we have $\Ax_\omega(\gamma)\cap[F_n(\widetilde{c})]_{n\in\NN}=\Ax_\omega(\gamma')\cap[F_n(\widetilde{c})]_{n\in\NN}$.\cqfd


\medskip

Let $\Gamma_{\widetilde{W}}$ be the stabilizer of $\widetilde{W}$ in $\Gamma$, let $\gamma_1,\gamma_2$ be two (not necessarily distinct) elements of
$\Gamma_{\widetilde{W}}$ 
that do not preserve any 
boundary component of $\widetilde{W}$
and let $\widetilde{c}$ be a boundary component
of $\widetilde{W}$.
Let $(\Ax_n(\gamma_1))_{n\in\NN}$ and $(\Ax_n(\gamma_2))_{n\in\NN}$ be two sequences of translation axes of $\gamma_1$ and $\gamma_2$, and let 
$\Ax_\omega(\gamma_1)=[\Ax_n(\gamma_1)]_{n\in\NN}$ and $\Ax_\omega(\gamma_2)=[\Ax_n(\gamma_2)]_{n\in\NN}$. Recall that $\gamma_{\widetilde{c}}$ is  a primitive generator 
of $\Gamma_{\widetilde{c}}$. For all $n\in\NN$, let $\widetilde{c}_n$ be the boundary component of 
$\widetilde{W}_n$ corresponding to $\widetilde{c}$.

\blemmA\label{compriss} If $[x_{1},x_{2}]$ is a geodesic segment joining $\Ax_\omega(\gamma_1)$ to $\Ax_\omega(\gamma_2)$ in $[\widetilde{W}_n]_{n\in\NN}$,
then the segment $[x_{1},x_{2}]$ can intersect $[\widetilde{c}_n]_{n\in\NN}$ at most at $\star_{\widetilde{c}}$. Notably, $\Ax_\omega(\gamma_1)$ and $\Ax_\omega(\gamma_2)$
can intersect $[\widetilde{c}_n]_{n\in\NN}$ at most at $\star_{\widetilde{c}}$. 
\elemmA

\demA There exists a sequence
$([x_{1,n},x_{2,n}])_{n\in\NN}$ of geodesic segments joining $\Ax_n(\gamma_1)$ to $\Ax_n(\gamma_2)$ such that $[x_1,x_2]=[[x_{1,n},x_{2,n}]]_{n\in\NN}$. 
Let $\widetilde{\ell}_0$ be a geodesic of $(\widetilde{\Sigma},[\widetilde{q}_0])$ that is 
interlaced with $\widetilde{c}$, and for all $n\in\NN$, let $\widetilde{\ell}_n$ be a geodesic of $\revetn$ having the same pair of points at infinity than 
$\widetilde{\ell}_0$. 
There exists $k\in\NN$ such that $\Ax_n(\gamma_1)$ and $\Ax_n(\gamma_2)$  are contained between $\gamma_{\widetilde{c}}^{-k}(\widetilde{\ell}_n)$ and $\gamma_{\widetilde{c}}^{k}(\widetilde{\ell}_n)$,
with $k$ independent of $n$ as in the proof of Lemma 
\ref{star}.

\begin{center}
	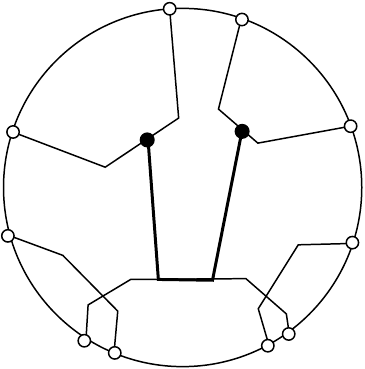
\end{center}

For all $n\in\NN$, the union $\Ax_n(\gamma_1)\cup[x_{1,n},x_{2,n}]\cup\Ax_n(\gamma_2)$ is contained in a connected generalized subsurface of $\widetilde{\Sigma}$
whose boundary is 
contained in 
$\gamma_{\widetilde{c}}^{k}(\widetilde{\ell}_n)\cup\gamma_{\widetilde{c}}^{-k}(\widetilde{\ell}_n)\cup\widetilde{c}_n$. Since we have
$[\gamma_{\widetilde{c}}^{k}(\widetilde{\ell}_n)]_{n\in\NN}\cap[\widetilde{c}_n]_{n\in\NN}=[\gamma_{\widetilde{c}}^{-k}(\widetilde{\ell}_n)]_{n\in\NN}\cap[\widetilde{c}_n]_{n\in\NN}
=\star_{\widetilde{c}}$,
we see that $[[x_{1,n},x_{2,n}]]_{n\in\NN}$ can intersect $[\widetilde{c}_n]_{n\in\NN}$ at most at $\star_{\widetilde{c}}$.\cqfd

\medskip

Let $x$ be a point of $\Ax_\omega(\gamma_1)$.

\blemmA\label{ilsysonttous} The geodesic segment $[x,\gamma_{\widetilde{c}}x]$ intersects $[\widetilde{c}_n]_{n\in\NN}$ at $\star_{\widetilde{c}}$.
\elemmA

\demA This follows from Lemma \ref{compriss} using $\gamma_2=\gamma_{\widetilde{c}}\gamma_1\gamma_{\widetilde{c}}^{-1}$ and from Remark \ref{convergencedistancetranslation}.\cqfd

\bigskip


In the remainder of this Section \ref{areazero}, we assume moreover that $\lim_\omega\operatorname{Diam}_{[\widehat{q}_n]}(\widehat{W}_n)=0$.

\blemmA\label{degenerationtypique} The ultralimit $[\widetilde{W}_n]_{n\in\NN}$ is an $\RR$-tree $T_{\widetilde{W}}$.\label{arbrereel}
\elemmA

\demA Since $\widetilde{W}_n$ is convex for all $n\in\NN$, the ultralimit $[\widetilde{W}_n]_{n\in\NN}$ is convex.
Since $\lim_\omega\operatorname{Diam}_{[\widehat{q}_n]}(\widehat{W}_n)$, for all $n\in\NN$, the union of $\partial \widetilde{W}_n$ with the set of singular points of $[\widetilde{q}_n]$
is $\varepsilon_n$-dense in $\widetilde{W}_n$, with $\lim_\omega\varepsilon_n=0$. According to Lemma \ref{hyperbolicity}, the convex set $[\widetilde{W}_n]_{n\in\NN}$
is $0$-hyperbolic, hence it is an $\RR$-tree.\cqfd

\medskip

Let $\Gamma_{\widetilde{W}}$ be the stabilizer of $\widetilde{W}$ in $\Gamma$. If $\Gamma_{\widetilde{W}}$ does not have a global fixed point in $T_{\widetilde{W}}$,
let $T_{\min}$ be the {\it minimal subtree} for the isometric action of $\Gamma_{\widetilde{W}}$ on $T_{\widetilde{W}}$, i.e. the smallest 
non empty subtree
of $T_{\widetilde{W}}$ that is invariant by $\Gamma_{\widetilde{W}}$. It is the union of the translation axes of the hyperbolic elements of $\Gamma_{\widetilde{W}}$. 
The isometric action of $\Gamma_{\widetilde{W}}$ on $T_{\min}$
is said to have {\it small edge stabilizers} 
if $\Gamma_{\widetilde{W}}$ has no global fixed point,
all the elements defined by the boundary components 
of $M$ are elliptic elements, and the stabilizers of the non trivial segments are trivial or cyclic. 

\blemmA\label{stabilizersdaretes}  Either $\Gamma_{\widetilde{W}}$ has a global fixed point in $T_{\widetilde{W}}$ or the isometric action of 
$\Gamma_{\widetilde{W}}$ on $T_{\min}$ has small edge stabilizers. 
\elemmA

\demA Assume that $\Gamma_{\widetilde{W}}$ has no global fixed point in $T_{\widetilde{W}}$ so that $T_{\min}$ is well defined. Let $\gamma\in\Gamma_{\widetilde{W}}$ be an element defined by a boundary component $b$ of 
$\widehat{W}$.
We have $\lim_\omega\ell_{[\widehat{q}_n]}(b)=0$, and according to Remark \ref{convergencedistancetranslation}, we have $\ell_{T_{\min}}(\gamma)=
\lim_\omega\ell_{\widetilde{W}_n}(\gamma)=0$, hence $\gamma$ is elliptic in $T_{\min}$. 

\medskip
Let $\alpha,\beta$ be two elements of
$\Gamma_{\widetilde{W}}-\{e\}$ which are not powers
of a common element. Let us prove that $\alpha$ and $\beta$ have at most one common fixed point. For all $n\in\NN$, let 
$\Ax_n(\alpha)$ and $\Ax_n(\beta)$ be two translation axes of $\alpha$ and $\beta$ in
$\widetilde{W}_n$, that minimize the distance between the translation axes of $\alpha$ and $\beta$. Let $[x_n]_{n\in\NN}$ and $[y_n]_{n\in\NN}$ be fixed points
of $\alpha$ and $\beta$ in $T_{\min}$.
For all $n\in\NN$, let $F_n(\alpha)$ be the flat strip union of all the translation axes of $\alpha$. 
According to Remark \ref{convergencedistancetranslation}, we have $d_n(x_n,F_n(\alpha))\leqslant d_n(x_n,\alpha x_n)$, hence $d_n(x_n,\Ax_n(\alpha))\leqslant 
d_n(x_n,\alpha x_n)+h_n(\alpha)$, where $h_n(\alpha)$ is the height of $F_n(\alpha)$.
Since ${h}_n(\alpha)\leqslant 2\operatorname{Diam}_{[\widehat{q}_n]}(\widehat{W}_n)$ $\omega$-converges to zero, we have
$\lim_\omega d_n(x_n,\Ax_n(\alpha))\leqslant \lim_\omega d_n(x_n,\alpha x_n)+\lim_\omega{h}_n(\alpha)=0$. 
Hence, we can always assume that $x_n$ 
belongs to $\Ax_{n}(\alpha)$ for all $n$. Similarly, we can always assume that $y_n$ 
belongs to $\Ax_{n}(\beta)$ for all $n$. 


Assume that $[x_n]_{n\in\NN}=[y_n]_{n\in\NN}$. Let $(z_n)_{n\in\NN}$ and $(w_n)_{n\in\NN}$ be two 
sequences of points such that $[z_n]_{n\in\NN}$ is a fixed point of $\alpha$ and $[w_n]_{n\in\NN}$ is a fixed point of $\beta$.
As above, we can assume that $z_n$ belongs to $\Ax_n(\alpha)$ and $w_n$ belongs to $\Ax_n(\beta)$ for all $n\in\NN$. 
Since $\lim_\omega h_n(\alpha)=\lim_\omega h_n(\beta)=0$, up to replacing $\Ax_n(\alpha)$ and $\Ax_n(\beta)$ by others
translation axes of $\alpha$ and $\beta$, 
we can always assume that the segment $[w_n,z_n]$ does not intersect a translation axis of $\alpha$ or $\beta$ other than $\Ax_n(\alpha)$ and 
$\Ax_n(\beta)$.
Since $\lim_\omega d_n(x_n,y_n)=0$ and $\lim_\omega\max\{\ell_{\widetilde{W}_n}(\alpha),\ell_{\widetilde{W}_n}(\beta)\}=0$, according to Lemma
\ref{doublepointfixes}, we have 
$\lim_\omega d_n(w_n,z_n)\geqslant\lim_\omega d_n(w_n,x_n)+\lim_\omega d_n(z_n, y_n)$. Hence $[w_n]_{n\in\NN}$ is equal to $[z_n]_{n\in\NN}$ if 
and only if
$[w_n]_{n\in\NN}=[z_n]_{n\in\NN}=[y_n]_{n\in\NN}=[x_n]_{n\in\NN}$. This proves the result.\cqfd

\medskip

Let $m$ be a  hyperbolic metric on $\Sigma$, and let $\widetilde{m}$ be its pullback on $\widetilde{\Sigma}$.
If $(\Lambda,\mu)$ is a measured hyperbolic lamination on $\srfcem$, it lifts to a measured hyperbolic lamination $(\widetilde{\Lambda},\widetilde{\mu})$ on 
$\revetm$. If $\widetilde{\lambda}$ is an isolated leaf of $\widetilde{\Lambda}$, the measure on any arc transverse to $\widetilde{\lambda}$, that is disjoint
of the rest of the lamination, is a Dirac measure at the intersection point between the arc and the leaf, of mass $\delta_{\widetilde{\lambda}}>0$.
We replace the leaf $\widetilde{\lambda}$ by a flat strip of width $\delta_{\widetilde{\lambda}}$ foliated by parallel leaves. Thus, we get a surface $\widetilde{\Sigma}'$
and a measured
lamination, with a transverse measure, such that the measures on the arcs transverse to the lamination have no atom. The map $d':\widetilde{\Sigma}'\times
\widetilde{\Sigma}'\to\RR^+$ defined by $d'(x,y)=
\underset{\widetilde{c}}\inf\,||\widetilde{\mu}_{\widetilde{c}}||$, where $\widetilde{c}$ is an arc transverse to $\widetilde{\Lambda}'$ between $x$ and $y$ and
$||\widetilde{\mu}_{\widetilde{c}}||$ is the total mass of $\widetilde{\mu}_{\widetilde{c}}$, is a pseudo-distance on $\widetilde{\Sigma}'$.
The quotient space $T_{(\Lambda,\mu)}=\widetilde{\Sigma}'/\sim$, with $x\sim y$ if $d'(x,y)=0$, endowed with the distance $d_\Lambda$ defined
by $d'$, is an $\RR$-tree,
called the {\it  $\RR$-tree dual to $(\Lambda,\mu)$} 
(see for example \cite[§~1]{MorSha91}), and the action of $\Gamma$ on 
$\widetilde{\Sigma}'$ defines an isometric action of $\Gamma$ on $T_{(\Lambda,\mu)}$.
For every element $\gamma\in\Gamma$, we denote by $\alpha_\gamma$ the associated free homotopy class of closed curves of $\Sigma$.
Recall that for every $\gamma\in\Gamma$, we have $\ell_{T_{(\Lambda,\mu)}}(\gamma)=i(\alpha_\gamma,(\Lambda,\mu))$ (see \cite[§~1]{MorSha91}).

\blemmA\label{dualite} If $\Gamma_{\widetilde{W}}$ has no global fixed point in $T_{\widetilde{W}}$, the minimal $\RR$-tree $T_{\min}$ endowed with the isometric action of 
$\Gamma_{\widetilde{W}}$ is dual to a measured hyperbolic lamination $(\Lambda,\mu)$ on $W$ (for any hyperbolic metric). 
\elemmA

\demA It is a consequence of Lemma \ref{stabilizersdaretes} and of \cite{Skora96}.\cqfd

\subsection{Decomposition of $\Sigma$ into pieces.}\label{decomposition}

In this short subsection, we introduce the notation which will be used in the next three subsections, that aim at decomposing the surface $\Sigma$ into revelant pieces.
For all $\alpha\in\C(\Sigma)$, let $\epsilon_\alpha=\lim_\omega\ell_{[q_n]}(\alpha)$. Recall that by the definition of $(\lambda_n)_{n\in\NN}$ and
$([q_n])_{n\in\NN}$ (see the beginning of Section \ref{sequenceofhalftranslations}), we have $\epsilon_\alpha<+\infty$ for every $\alpha\in\C(\Sigma)$.   Let $\SSSS_0=\{\alpha\in\SSSS(\Sigma)\;:\;
\epsilon_{\alpha_0}=0\}$.
If $\SSSS_0\not=\emptyset$, we denote by $\Sigma_0$ the tight subsurface filled up by $\SSSS_0$ (defined up to isotopy, see Lemma \ref{isotopy}).

Let $W$ be a tight connected subsurface of $\Sigma$. Let $\widehat{p}:\Sigma_W\to\Sigma$ be a $W$-cover of $\Sigma$, let 
$\widehat{W}$ be the connected component of
the preimage of $W$ in $\Sigma_W$ that is 
not simply connected, and let $\widetilde{W}$ be a connected component of the preimage of $W$ in $\widetilde{\Sigma}$.
For all 
$n\in\NN$, let $[\widehat{q}_n]$ and $[\widetilde{q}_n]$ be the pullbacks of $[q_n]$ on $\Sigma_W$ and on $\widetilde{\Sigma}$, and let
$\widehat{W}_n$ and $\widetilde{W}_n$ be the $[\widehat{q}_n]$ and $[\widetilde{q}_n]$-geometric realizations of $\widehat{W}$ and $\widetilde{W}$.
Finally, let $\Gamma_{\widetilde{W}}$ be the stabilizer of $\widetilde{W}$ in $\Gamma$.

Let us introduce the following example of typical degenerating sequence of half-translation structures on a surface.

\begin{center}\label{example}
 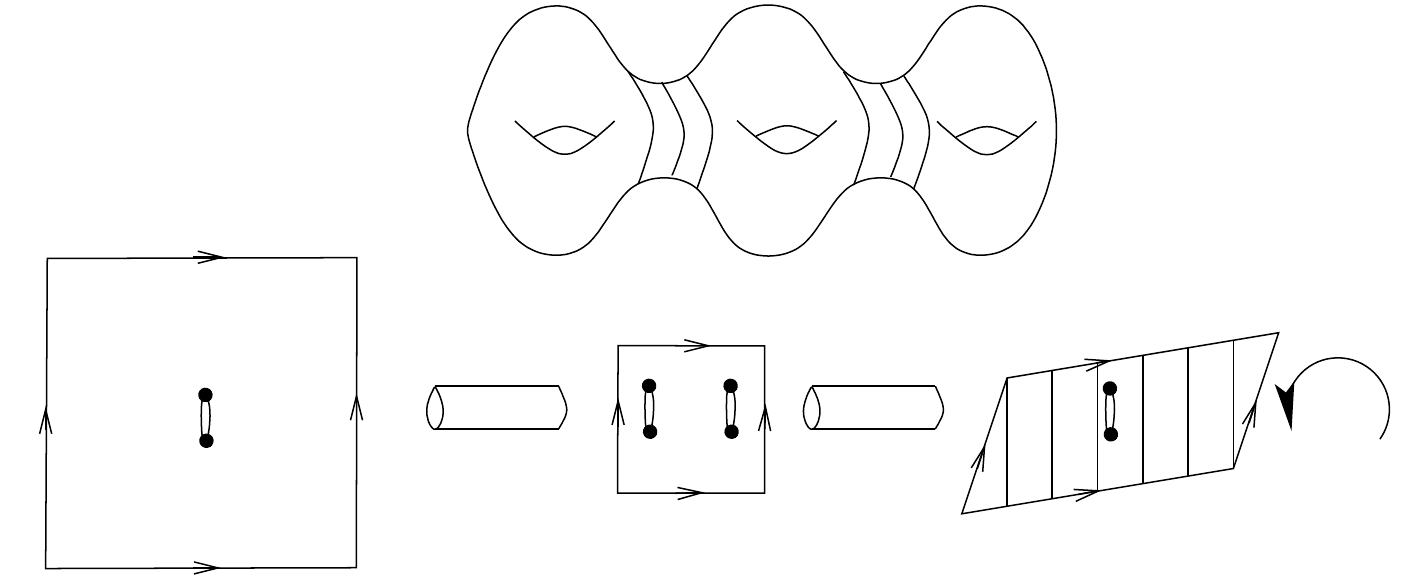
\end{center}

Let $n\in\NN$. The half-translation structure $[q_n]$ on $\Sigma$ is defined piece by piece. It is obtained by gluing by two flat cylinders of height 
$1$ and girth $\frac{1}{n}$ along the slits, the unit square flat torus with a vertical slit of perimeter $\frac{1}{n}$, a square flat torus with side lengths $\frac{3}{n}$ and two
vertical slits of perimeters $\frac{1}{n}$, and a square torus whose both vertical and horizontal measured foliations are minimal, whose transverse
measure to the horizontal foliation is multiplied by $\frac{1}{n}$, and whose transverse measure to the vertical foliation is constant.   
The singularities are at the $8$ marked points, and are of angle $3\pi$. 

Recall that for all $n\in\NN$, if $(\F_n^v,\mu_n^v)$ and $(\F_n^h,\mu_n^h)$ are the vertical and the horizontal foliations of a representative $q_n\in\Q(\Sigma)$ of $[q_n]$,
for all $\alpha\in\C(\Sigma)$, we have $\max\{i((\F_n^v,\mu_n^v),\alpha),i((\F_n^h,\mu_n^h),\alpha)\}\leqslant\ell_{[q_n]}(\alpha)\leqslant
i((\F_n^v,\mu_n^v),\alpha)+i((\F_n^h,\mu_n^h),\alpha)$. Hence, 

\medskip
\begin{itemize}
 \item[$\bullet$] $\lim_\omega\ell_{[q_n]}(\beta_1)=\lim_\omega\ell_{[q_n]}(\beta_2)=0$;
 \item[$\bullet$] for every $\alpha\in\SSSS(W_0)$, we have $\lim_\omega\ell_{[q_n]}(\alpha)=0$;
 \item[$\bullet$] for every $\alpha\in\SSSS(W_1)$, we have  $\lim_\omega\ell_{[q_n]}(\alpha)=\ell_{[q]}(\alpha)$, with $[q]$ the metric on the square flat torus of area $1$, minus a point;
 \item[$\bullet$] for every $\alpha\in\SSSS(W_2)$, we have $\lim_\omega\ell_{[q_n]}(\alpha)=i(\alpha,(\F^v,\mu^v))$, where $(\F^v,\mu^v)$ is the vertical foliation of 
 $q_n$ for all $n\in\NN$.
\end{itemize}

In this example, we have $\Sigma_0=W_0$. 

\subsection{The subsurface $\Sigma_0$.}\label{Sigma0}

If $\Sigma_0$ is not a union of cylinders, assume in this subsection \ref{Sigma0} that $W$ is a connected component of $\Sigma_0$ that is not a cylinder 
(the cylinder components will be considered in Subsection \ref{cylindresetpantalons}). 


There exists a subset $S_0(W)$ of $\SSSS_0$ with at most $K\in\NN$ elements 
(see Lemma \ref{FFini}) that fills up $W$.  
We can always assume that $S_0(W)$ does not contain the isotopy class
of any boundary component of $W$. Let $S_0(\widehat{W})$ be the set of lifts of the elements of $S_0(W)$ in $\widehat{W}$. 

\blemmA\label{crunch} We have $\lim_\omega\operatorname{Diam}_{[\widehat{q}_n]}(\widehat{W}_n)=0$, $\lim_\omega\operatorname{Area}_{[\widehat{q}_n]}(\widehat{W}_n)=0$ and
for every boundary component $\beta$ of $\widehat{W}$, we have $\lim_\omega\ell_{[\widehat{q}_n]}(\beta)=0$. Hence, the ultralimit 
$[\widetilde{W}_n]_{n\in\NN}$ is an $\RR$-tree $T_{\widetilde{W}}$.\elemmA

\demA The set $S_0(\widehat{W})$ has at most $K<+\infty$ elements and, for every $\alpha\in S_0(\widehat{W})$, we have
$\lim_\omega\ell_{[\widehat{q}_n]}(\alpha)=0$. Hence, for every $\varepsilon>0$, there exists
$I\in\omega$ such that for all $n\in I$ and $\alpha\in S_0(\widehat{W})$, we have
$\ell_{[\widehat{q}_n]}(\alpha)<\varepsilon$. Then, according to Lemma \ref{diameterandarea}, we
have $\operatorname{Diam}_{[\widehat{q}_n]}(\widehat{W}_n)\leqslant 11K\varepsilon$ 
and $\Area_{[\widehat{q}_n]}(\widehat{W}_n)\leqslant\frac{1}{\pi}(1+2b)^2(K\varepsilon)^2$ (where $b$ is the number of boundary components of $\widehat{W}$).
Hence $\lim_\omega\operatorname{Diam}_{[\widehat{q}_n]}(\widehat{W}_n)=0$ and $\lim_\omega\operatorname{Area}_{[\widehat{q}_n]}(\widehat{W}_n)=0$.
Hence, according to Lemma \ref{degenerationtypique}, the ultralimit $[\widetilde{W}_n]_{n\in\NN}$ is an $\RR$-tree.\cqfd

%

\blemmA\label{elliptic} The action of $\Gamma_{\widetilde{W}}$ on $T_{\widetilde{W}}$ has a unique 
global fixed point $\star_{\widetilde{W}}$, and for every boundary component $\widetilde{c}$ of $\widetilde{W}$, the point $\star_{\widetilde{c}}$ (see Lemma \ref{star}) is 
equal to $\star_{\widetilde{W}}$.
\elemmA

We will need the following well-known lemma, obtained by smoothing the Dirac masses of negative curvature at the singularities, see for instance \cite{Minsky92} for a
similar result.

\blemmA\label{smoothing} Let $[\widehat{q}]$ be a half-translation structure on  $\Sigma_W$, and let $E$ be a finite set of homotopy classes of
essential closed curves in $\widehat{W}$. For every $\varepsilon>0$, there  
exists a Riemannian metric $\widehat{g}$ with nonpositive curvature on the smooth surface $\Sigma_W$ such that $|\ell_{[\widehat{q}]}(\alpha)-\ell_{\widehat{g}}(\alpha)|<\varepsilon$ for all $\alpha\in E$.
\elemmA

{\noindent{\bf Proof of Lemma \ref{elliptic}.}} Let $\gamma$ be an element of $\Gamma_{\widetilde{W}}-\{e\}$, and let $\widehat{\alpha}_\gamma$ be the free homotopy
class of closed curves defined by $\gamma$ in $\widehat{W}$. For every $\varepsilon>0$, 
there exists $I\in\omega$ such that for all $n\in I$ and $s_0\in S_0(\widehat{W})$, we have $\ell_{[\widehat{q}_n]}(s_0)\leqslant\varepsilon$.
Let $E=S_0(\widehat{W})\cup\{\widehat{\alpha}_\gamma\}$.
According to Lemma \ref{smoothing}, there exists a Riemannian metric $\widehat{g}$ with nonpositive curvature on $\Sigma_{W}$ such that for every $\alpha\in E$, we have 
$|\ell_{[\widehat{q}_n]}(\alpha)-\ell_{\widehat{g}}(\alpha)|<\varepsilon$. Let $\widehat{W}_{\widehat{g}}$ be the $\widehat{g}$-geometric realization of $\widehat{W}$.
Since $S_0(\widehat{W})$ fills up $\widehat{W}$,
the union $G$ of some $\widehat{g}$-geodesic representatives of the
elements of $S_0(\widehat{W})$ (one by element of $S_0(\widehat{W})$) is a graph whose complementary components in $\widehat{W}_{\widehat{g}}$ are disks and half-open
cylinders that can be homotoped to some boundary components of $\widehat{W}_{\widehat{g}}$.
Moreover, since the lengths of these $\widehat{g}$-geodesics are at most $2\varepsilon$, and $S_0(\widehat{W})$ has at most $K$ elements, the sum of
the lengths of the edges of $G$ is at most $2K\varepsilon$. Hence, the perimeters of the complementary disks are at most $4K\varepsilon$ and the circonferences of the 
complementary cylinders are at most $2K\varepsilon$.
Since $\widehat{g}$ is a Riemannian metric with nonpositive curvature, the length of an intersection
segment of a $\widehat{g}$-geodesic representative of $\widehat{\alpha}_\gamma$ with a complementary component of $G$ is at most $4K\varepsilon$,
and since the intersection number between two $\widehat{g}$-geodesics is minimal, we have $\ell_{\widehat{g}}(\widehat{\alpha}_\gamma)\leqslant 2K\varepsilon
\sum_{s_0\in S_0(\widehat{W})}i(\widehat{\alpha}_{\gamma},s_0)$, and $\ell_{[\widehat{q}_n]}(\widehat{\alpha}_\gamma)\leqslant  2K\varepsilon 
\sum_{s_0\in S_0(\widehat{W})}i(\widehat{c}_\gamma,s_0)+\varepsilon$.
Hence, we have $\ell_{T_{\widetilde{W}}}(\gamma)=\lim_\omega\ell_{[\widehat{q}_n]}(\widehat{c}_\gamma)=0$, so that $\gamma$ has a fixed point in $T_{\widetilde{W}}$. 
Hence, all the elements of $\Gamma_{\widetilde{W}}$ are elliptic in $T_{\widetilde{W}}$, and according to a lemma of Serre (see \cite[p.~271]{Shalen87}), $\Gamma_{\widetilde{W}}$
has a global fixed point. Moreover, according to Lemma \ref{doublepointfixes} used as in the proof of Lemma \ref{stabilizersdaretes},
two elements of $\Gamma_{\widetilde{W}}$ which are not powers of a common element have 
at most one common fixed point. Hence, the global fixed point is unique . 

\medskip

Moreover, if $\gamma_1\in\Gamma_{\widetilde{W}}-\{e\}$ preserves a boundary component $\widetilde{c}$ of $\widetilde{W}$ and $\gamma_2\in\Gamma_{\widetilde{W}}$ does not preserve any boundary component
of $\widetilde{W}$, and if $(\Ax_n(\gamma_1))_{n\in\NN}$ and $(\Ax_n(\gamma_2))_{n\in\NN}$ are two sequences of translation axes of $\gamma_1$ and $\gamma_2$,
according to Lemma \ref{compriss}, the geodesic
$[\Ax_n(\gamma_1)]_{n\in\NN}$ can intersect $[\Ax_n(\gamma_2)]_{n\in\NN}$
at most at $\star_{\widetilde{c}}$. Since the sets $[\Ax_n(\gamma_1)]_{n\in\NN}$ and $[\Ax_n(\gamma_2)]_{n\in\NN}$ are the sets of fixed points of $\gamma_1$ and $\gamma_2$
in $T_{\widetilde{W}}$, the unique common point of $\gamma_1$ and $\gamma_2$ in $T_{\widetilde{W}}$ is $\star_{\widetilde{c}}$.
Hence, for every boundary component $\widetilde{c}$ of $\widetilde{W}$, the point $\star_{\widetilde{c}}$ is equal to $\star_{\widetilde{W}}$.\cqfd

\medskip

\remA According to Lemmas \ref{crunch} and \ref{elliptic}, in both cases ($\lim_\omega\lambda_n=+\infty$ or $\lim_\omega\lambda_n<+\infty$), the subsurface $\Sigma_0$
cannot be equal to $\Sigma$, since then $\Gamma$ would have a global fixed point in $\lim_\omega\revetn$.

\subsection{Complementary connected components of $\Sigma_0$ in $\Sigma$.}\label{composantescomplementaire}

In this Subsection \ref{composantescomplementaire}, we assume that $W$ is the closure of a connected component of $\Sigma-\Sigma_0$.
Then $W$ is non trivial, compact, connected and $\pi_1$-injective.
We assume first that $W$ is neither a cylinder nor a pair of pants. 


\blemmA\label{minorationuniforme} Either there exists $\varepsilon>0$ and $I\in\omega$ such that for all $n\in I$ and $\alpha\in\SSSS(W)$, we have
$\ell_{[q_n]}(\alpha)\geqslant\varepsilon$,
or there exists a sequence $(\alpha_k)_{k\in\NN}$ in $\SSSS(W)$ such that $(\epsilon_{\alpha_k})_{k\in\NN}$ (strictly) decreases to zero. 
\elemmA

\remA In the first case, for all $n\in I$ and $\alpha\in\C(W)$, we have $\ell_{[q_n]}(\alpha)\geqslant\varepsilon$. We will see that this case cannot happen if 
$\lim_\omega\lambda_n=+\infty$.

\medskip

We will need the following two lemmas. Let $m$ be a hyperbolic metric on $\Sigma$ and let $W_m$ be the $m$-geodesic realization of $W$.
Let $(\Lambda,\mu)$ be a measured hyperbolic lamination of $W_m$.

\blemmA\label{pasdeminoration} There exists a possibly constant sequence $(\alpha_k)_{k\in\NN}$ in $\SSSS(W)$ such that the sequence 
$(i((\Lambda,\mu),\alpha_k))_{k\in\NN}$
decreases to
zero. 
\elemmA

\demA If $\Lambda$ contains a simple closed geodesic $\beta$, we take $\alpha_k$ equal to the homotopy class of $\beta$ for all $k\in\NN$. Otherwise, 
let $L_1$ be a minimal component of
$\Lambda$. For every $k\in\NN$, let $c$ be a transverse arc to $L_1$, that meets $L_1$ and does not meet any other minimal component of $\Lambda$,
such that the total $\mu$-mass of $c$ is at most $\frac{1}{k}$.
Since the intersection of any leaf of $L_1$ with $c$ is dense in the Cantor set $L_1\cap c$,
there exist two distinct points in $c$ belonging to the same leaf. Then, if $\alpha_k'$ is the concatenation of the segment of this leaf and of the segment of $c$ 
between 
these two points, we have $i(\alpha_k',(\Lambda,\mu))\leqslant\frac{1}{k}$. Moreover, the image of $\alpha_k'$ contains the image of an essential simple
closed curve $\alpha_k\in\SSSS(W)$, and $i(\alpha_k,(\Lambda,\mu))\leqslant i(\alpha_k',(\Lambda,\mu))\leqslant\frac{1}{k}$.\cqfd

\medskip


We recall that $(\Lambda,\mu)$ {\it fills up} $W$ if $i(\alpha,(\Lambda,\mu))>0$ for all 
$\alpha\in\C(\Sigma)$ that topologically cuts $W$ (or equivalently for all $\alpha\in\C(W)$, since $W$ is not a cylinder). If $(\Lambda',\mu')$ is another measured hyperbolic lamination 
of $W_m$, recall that $(\Lambda,\mu)$  
and $(\Lambda',\mu')$  {\it jointly fill up} $W$ if $i(\alpha,(\Lambda,\mu))+i(\alpha,(\Lambda',\mu'))>0$ for all $\alpha\in\C(W)$, 
and $i((\Lambda,\mu),(\Lambda',\mu'))>0$.

\blemmA\label{jointlyfill} Assume that $(\Lambda,\mu)$ and $(\Lambda',\mu')$ jointly fill up $W$, and let $(\Lambda_n,\mu_n)_{n\in\NN}$ and 
$(\Lambda_n',\mu_n')_{n\in\NN}$ 
be two sequences of measured hyperbolic laminations of $\Sigma$ that respectively $\omega$-converge to two measured hyperbolic laminations of $\Sigma$ whose intersections 
with $W_m$ are equal to $(\Lambda,\mu)$ and $(\Lambda',\mu')$.
There exist $I\in\omega$ and 
$\varepsilon>0$ such that $i(\alpha,(\Lambda_n,\mu_n))+i(\alpha,(\Lambda_n',\mu_n'))\geqslant\varepsilon$ for all $n\in I$ and $\alpha\in\SSSS(W)$. 
\elemmA

\demA Assume for a contradiction that there exists a sequence $(\alpha_n)_{n\in\NN}$ in $\SSSS(W)$ such that $\lim_\omega i(\alpha_n,(\Lambda_n,\mu_n))=\lim_\omega i(\alpha_n,(\Lambda'_n,\mu'_n))=0$.
Then, there exists a positive real sequence $(t_n)_{n\in\NN}$ such that $(t_n\alpha_n)_{n\in\NN}$
$\omega$-converges to a measured hyperbolic lamination $(L,\nu)$ of $W_m$. First, $L$ is not a boundary component of $W_m$ since then $\alpha_n$ 
would spiral around this boundary component, for $n$ large enough, which is impossible for an essential simple closed geodesic of $W_m$. Moreover, we have
$\lim_\omega t_n<+\infty$,
hence $$i((L,\nu),(\Lambda,\mu))+
i((L,\nu),(\Lambda',\mu'))=0$$ If $L$ contains a simple closed geodesic $\alpha$, then $i(\alpha,(\Lambda,\mu))+i(\alpha,(\Lambda',\mu'))=0$, which is impossible since 
$(\Lambda,\mu)$ and $(\Lambda',\mu')$ jointly fill up $W$. Let 
$L_1$ be a minimal component of $L$, and let $\Sigma(L_1)\subseteq W_m$ be the connected subsurface filled  up by $L_1$ (see Definition \ref{remplieparunelamination}).
If 
$\Sigma(L_1)$ were not equal to $W_m$ up to isotopy, there would exist a boundary component $\alpha$ of $\Sigma(L_1)$ which would be essential in $W$, and then
$i(\alpha,(\Lambda,\mu))+i(\alpha,(\Lambda',\mu'))=0$, a contradiction. 
Hence, the lamination $L_1$ fills up $W$. By the condition $i((L,\nu),(\Lambda,\mu))+
i((L,\nu),(\Lambda',\mu'))=0$, we have  $\Lambda\subseteq L_1$ and $\Lambda'\subseteq L_1$, and since $L_1$ is minimal, we have $L_1=\Lambda=\Lambda'$, a contradiction
since $i((\Lambda,\mu),(\Lambda',\mu'))>0$.\cqfd

\medskip
\noindent
{\bf Proof of Lemma \ref{minorationuniforme}.} For all $n\in\NN$, we choose a quadratic differential $q_n$ representing $[q_n]$,
and we denote by $(\F_n^h,\mu_n^h)$ and $(\F_n^v,
\mu_n^v)$ (the equivalence classes of) its horizontal and vertical measured foliations. There exist two positive real sequences $(x_n)_{n\in\NN}$ and
$(y_n)_{n\in\NN}$ such that $x_n(\F_n^h,\mu_n^h)_{n\in\NN}$ and $y_n(\F_n^v,\mu_n^v)_{n\in\NN}$ $\omega$-converge 
to two (equivalence classes of) measured foliations $(\F^h,\mu^h)$ and $(\F^v,\mu^v)$. For all $n\in\NN$ and 
$\alpha\in\C(\Sigma)$, we have
\begin{equation}\label{label_equation}
	\max\{i(\alpha,(\F_n^h,\mu_n^h)),i(\alpha,(\F_n^v,\mu_n^v))\}\leqslant\ell_{[q_n]}(\alpha)\leqslant
	i(\alpha,(\F_n^h,\mu_n^h))+i(\alpha,(\F_n^v,\mu_n^v)).
\end{equation}

Assume for a contradiction that $\lim_\omega x_n=0$. Since $\mu^h$ is nonzero, there exists $\alpha\in\SSSS(\Sigma)$ such that $\lim_\omega 
i(\alpha,(\F_n^h,x_n\mu_n^h))=i(\alpha,(\F^h,\mu^h))>0$. However, $\lim_\omega i(\alpha,(\F_n^h,x_n\mu_n^h))\leqslant \lim_\omega x_n\ell_{[q_n]}(\alpha)=0$, since 
$\lim_\omega\ell_{[q_n]}(\alpha)<+\infty$. Hence, $\lim_\omega x_n$ is nonzero, and similarly $\lim_\omega y_n$ is nonzero.

Assume for a contradiction that $\lim_\omega x_n=\lim_\omega y_n=+\infty$. Then, for all $\alpha\in\SSSS(W)$, we have 
\begin{align*}
	\lim_\omega\ell_{[q_n]}(\alpha)&\leqslant\lim_\omega (i(\alpha,(\F_n^h,\mu_n^h))+i(\alpha,(\F_n^v,\mu_n^v)))\\
	&\leqslant\lim_\omega (\frac{1}{x_n}i(\alpha,(\F_n^h,x_n\mu_n^h))+\frac{1}{y_n}i(\alpha,(\F_n^v,y_n\mu_n^v)))\\
	&=0,
\end{align*}

\noindent
a contradiction since $\alpha$ is essential in $W$ and hence does not belong to $\SSSS_0$. Hence, up to permuting $(\F_n^h,\mu_n^h)_{n\in\NN}$ and $(\F_n^v,\mu_n^v)_{n\in\NN}$, we can assume that $\lim_\omega x_n<+\infty$.
This means that the sequence $(\F_n^h,\mu_n^h)_{n\in\NN}$, without renormalization, 
$\omega$-converges to a (nonzero equivalence class of)
measured foliation, still denoted by $(\F^h,\mu^h)$.

\medskip

Assume first that $\lim_\omega y_n<+\infty$. Then, the sequence $(\F_n^v,\mu_n^v)_{n\in\NN}$ also $\omega$-converges, without renormalization, 
to a (nonzero, equivalence class of) measured foliation, still denoted by $(\F^v,\mu^v)$. If $\alpha$ is a boundary component of $W$, it is 
also a boundary component of $\Sigma_0$
and we have 
$\lim_\omega\ell_{[q_n]}(\alpha)=0$. Hence, by Equation \eqref{label_equation} $i(\alpha,(\F^h,\mu^h))=i(\alpha,(\F^v,\mu^v))=0$. Consequently, if $(\Lambda^h,\nu^h)$ and $(\Lambda^v,\nu^v)$ are 
the measured
hyperbolic laminations of $\srfcem$ associated with $(\F^h,\mu^h)$ and $(\F^v,\mu^v)$ (see for instance \cite{Levitt81}), then no leaf of $\Lambda^v$
or $\Lambda^h$ intersects the boundary of $W_m$, and the intersections $\Lambda^h_W=\Lambda^h\cap W_m$  and $\Lambda^v_W=\Lambda^v\cap W_m$,
endowed with the induced transverse measures 
(if not empty), are measured hyperbolic laminations.

If $(\Lambda^h_W,\mu^h)$ and $(\Lambda_W^v,\mu^v)$ jointly fill up $W$, according to Lemma \ref{jointlyfill},
there exist $I\in\omega$ and $\varepsilon>0$ such that, for all $n\in I$ and $\alpha\in\SSSS(W)$, we have $i(\alpha,(\F_n^h,\mu_n^h))+i(\alpha,(\F_n^v,\mu_n^v))
\geqslant\varepsilon$, hence $\ell_{[q_n]}(\alpha)\geqslant\max\{i(\alpha,(\F_n^h,\mu_n^h)),i(\alpha,(\F_n^v,\mu_n^v))\}\geqslant\frac{\varepsilon}{2}$, so that the first 
conclusion of Lemma \ref{minorationuniforme} holds.

Recall that there exists no $\alpha\in\SSSS(W)$ such 
that $i(\alpha,(\Lambda^h,\mu^h))+i(\alpha,(\Lambda^v,\mu^v))=0$,  since then 
$\alpha$ would belong
to $\SSSS_0$ by Equation \eqref{label_equation}, and would have a representative contained in $\Sigma_0$.
Hence, if $(\Lambda^h_W,\mu^h)$ and $(\Lambda_W^v,\mu^v)$ do not jointly fill up $W$, then 
$i((\Lambda^h_W,\mu^h),(\Lambda_W^v,\mu^v))=0$. Assume first that $i((\Lambda^h_W,\mu^h),(\Lambda_W^v,\mu^v))=0$ and that the two sets $\Lambda^h_W$ and $\Lambda_W^v$ are
not empty. Let 
$\Sigma(\Lambda^h_W)\subseteq W_m$ be the subsurface with $m$-geodesic boundary filled up by $\Lambda^h_W$.
Assume for a contradiction that $\Sigma(\Lambda^h_W)$ is strictly contained in $W_m$. Let $\beta$ be a boundary component of $\Sigma(\Lambda^h_W)$. 
Since $\Lambda^h_W$ fills up $\Sigma(\Lambda^h_W)$ and $i((\Lambda^h_W,\mu^h),(\Lambda_W^v,\mu^v))=0$, we have $i(\beta,(\Lambda^v_W,\mu^v_W))=0$ and 
$i(\beta,(\Lambda^h_W,\mu^h_W))=0$ since $\beta$ is a boundary component of  $\Sigma(\Lambda^h_W)$. Hence $i(\beta,(\Lambda^h,\mu^h))+i(\beta,(\Lambda^v,\mu^v))=0$,
a contradiction since $\beta$ is essential in $W$. Hence, 
$\Lambda_W^h$ fills $W$ and similarly $\Lambda_W^v$ fills $W$. Since $i((\Lambda^h,\mu^h),(\Lambda^v,\mu^v))=0$, one of these hyperbolic laminations is contained in the
otherone. Moreover, the two laminations are minimal since $W$ is connected, and hence they are equal.  
Hence, either one of the sets $\Lambda^h_W,\Lambda^v_W$ is empty, or the laminations  
$\Lambda^h_W$ and $ \Lambda^v_W$ are equal. In both cases,  according to Lemma \ref{pasdeminoration}, there exists a sequence $(\alpha_k)_{k\in\NN}$ of 
$\SSSS(W)$ such that $(i(\alpha_k,(\Lambda^h_W,\mu^h))+i(\alpha_k,(\Lambda_W^v,\mu^v)))_{k\in\NN}$ strictly decreases to zero, hence the second conclusion of Lemma 
\ref{minorationuniforme} holds by Equation \eqref{label_equation}.

\medskip 

Assume  next that $\lim_\omega y_n=+\infty$. If $\beta$ is a boundary component of $W$, as above, by Equation \eqref{label_equation}, we have 
$i(\beta,(\F^h,\mu^h))\leqslant\lim_\omega
\ell_{[q_n]}(\beta)=0$.    
Consequently, if we denote by $(\Lambda^h,\nu^h)$ the measured hyperbolic
lamination of $\srfcem$ associated with $(\F^h,\mu^h)$, then no leaf of $\Lambda^h$
intersects the boundary of $W_m$. Moreover, the intersection $\Lambda^h_W=\Lambda^h\cap W_m$ is not empty, otherwise all the elements of $\SSSS(W)$ would belong to
$\SSSS_0$ by Equation \eqref{label_equation}. Hence, the set $\Lambda_W^h$, endowed with the induced transverse measure, is a measured hyperbolic lamination. 
For all $\alpha\in\SSSS(W)$ and  $n\in\NN$,
we have 
$$\max\{i(\alpha,(\F_n^h,\mu_n^h)),\frac{1}{y_n}i(\alpha,(\F_n^v,y_n\mu_n^v))\}\leqslant\ell_{[q_n]}(\alpha)\leqslant
i(\alpha,(\F_n^h,\mu_n^h))+\frac{1}{y_n}i(\alpha,(\F_n^v,y_n\mu_n^v)).$$ 
Hence, $\lim_\omega\ell_{[q_n]}(\alpha)=i(\alpha,(\Lambda^h,\nu^h))$. As above, there exists a sequence $(\alpha_k)_{k\in\NN}$ of 
$\SSSS(W)$ such that $i(\alpha_k,(\Lambda^h,\nu^h))$ strictly decreases to zero, hence the second conclusion of Lemma 
\ref{minorationuniforme} holds by Equation \eqref{label_equation}.\cqfd

\subsection{Case of degeneration.}\label{degenration}

In this Subsection \ref{degenration}, we assume as in the previous one that $W$ is the closure of a complementary connected component of $\Sigma_0$ which is neither a cylinder nor a 
pair of pants.
We will study the asymptotic behavior of the restriction
to $W$ of the sequence of half-translation structures $([{q}_n])_{n\in\NN}$ under the hypothesis that the second conclusion of Lemma \ref{minorationuniforme} holds, i.e.
there exists a sequence $(\alpha_k)_{k\in\NN}$ of $\SSSS(W)$ such 
that $(\epsilon_{\alpha_k})_{k\in\NN}$ strictly decreases to zero. 

\rem This case corresponds to the subsurface $W_2$ of the example of Subsection \ref{decomposition}, since for every 
$\alpha\in\SSSS(W_2)$ we have $\lim_\omega\ell_{[q_n]}(\alpha)=i(\alpha,(\F^v,\mu^v))$, and $(\F^v,\mu^v)$ is filling.  

\blemmA\label{caseodegeneration}
If there exists a sequence $(\alpha_k)_{k\in\NN}$ such that $(\epsilon_{\alpha_k})_{k\in\NN}$ strictly decreases to zero, then
we have $\lim_\omega\operatorname{Diam}_{[\widehat{q}_n]}(\widehat{W}_n)=0$, $\lim_\omega\operatorname{Area}_{[\widehat{q}_n]}(\widehat{W}_n)=0$ and
for every boundary component $\beta$ of $\widehat{W}$, we have $\lim_\omega\ell_{[\widehat{q}_n]}(\beta)=0$. And the ultralimit $[\widetilde{W}_n]_{n\in\NN}$ is 
an $\RR$-tree, on which the isometric action of $\Gamma_{\widetilde{W}}$ has no global fixed point. Moreover, the minimal subtree of $[\widetilde{W}_n]_{n\in\NN}$ 
for the action of $\Gamma_{\widetilde{W}}$ is dual to a measured hyperbolic lamination $(\Lambda,\mu)$
(for any hyperbolic metric on $\Sigma$), which fills up $W$.
\elemmA

\remA Since $W$ is connected, if $\Lambda$ is filling, it is also minimal.

\medskip

\demA For all $p\in\NN$, we denote by 
${W}_p$ the subsurface of $\Sigma$ filled up by $\{\alpha_k\,:\,k\geqslant p\}$ (defined up to isotopy). According to Lemma \ref{subfamily}, if $p_2\geqslant p_1$, up to 
isotopy, the 
subsurface
${W}_{p_2}$ is contained in the interior of $W_{p_1}$. Moreover, if ${W}_{p_2}$ is not isotopic to ${W}_{p_1}$, then at least one connected component of ${W}_{p_1}-W_{p_2}$
has a negative Euler characteristic,
hence $\chi(W_{p_2})>\chi(W_{p_1})$. Since the Euler characteristic of ${W}_p$ is non positive for all $p\in\NN$, there exists $P\in\NN$ such that ${W}_p$ is isotopic to ${W}_P$
for all $p\geqslant P$. We denote by ${W}'_P$ a (isotopy class of a) connected component of ${W}_P$. Since  ${W}'_P$ is filled up by some essential simple closed curves,
it is not a pair of pants.
Assume for a contradiction that
it is a
cylinder. Let $\alpha$ be a isotopy class of simple closed curves contained in ${W}'_P$. For every $\varepsilon>0$, there exists $p\geqslant P$ such that $\epsilon_{\alpha_p}
\leqslant\varepsilon$, and $\alpha_p$ fills up ${W}'_P$. Since $W_P'$ is assumed to be a cylinder, we have $\alpha=\alpha_p^{\pm 1}$, and $\epsilon_\alpha=\epsilon_{\alpha_p}<\varepsilon$. 
Hence $\alpha$ belongs to $\SSSS_0$, which is impossible since $\alpha$ is essential in $W$. Hence, the  subsurface ${W}_P'$ is not a cylinder.

Let $\varepsilon>0$. There exists $p_0>P$ such that $\epsilon_{\alpha_k}<\frac{\varepsilon}{2}$ for all $k>p_0$, and according to Lemma \ref{FFini}, there exists a subset
of $\{\alpha_k\}_{k>p_0}$, with at most $K\in\NN$ elements (where $K$ only depends on the topology of ${W}$), that fills up ${W}'_P$. Since $K<+\infty$, there exists $I\in\omega$
such that $\ell_{[q_n]}(\alpha)\leqslant\varepsilon$ for all $\alpha$ in this subset and $n\in I$. Let $\Sigma_{W'_P}$ be a $W'_P$-cover of 
$\Sigma$, let $\widehat{W}'_P$ be the connected component of the preimage of $W'_P$ in $\Sigma_{W_P'}$ which is not simply connected, and for all $n\in I$, let 
$[\widehat{q}_n]$ be the pullback of $[q_n]$ to 
$\Sigma_{W'_P}$. We denote by $\widehat{W}'_{P,n}$ the $[\widehat{q}_n]$-geometric realization of $\widehat{W}'_P$. According to Lemma
\ref{diameterandarea}, the length of any boundary component of $\widehat{W}'_{P,n}$ is at most $K\varepsilon$, its
diameter is at most $11K\varepsilon$ and its area is at most $\frac{1}{\pi}(1+2b)^2(K\varepsilon)^2$ (where $b$ is the number of boundary components
of $\widehat{W}'_P$). Hence, if $\beta$ is (the free homotopy class of) a boundary component of $W'_P$, then
$\epsilon_\beta=0$ and $\beta$ belongs to $\SSSS_0$. Therefore $\beta$ is a boundary component of $\Sigma_0$ and of $W$. Hence, $W'_P=W_P=W$
(we will then replace $W'_P$ and $W_P$ by $W$). Since the diameter of $\widehat{W}_n$
$\omega$-converges to zero, according to Lemma \ref{degenerationtypique}, the ultralimit $[\widetilde{W}_n]_{n\in\NN}$ is an $\RR$-tree. Recall that by the definition of 
$\Sigma_0$, for all $\alpha\in\SSSS(W)$, we have 
$\epsilon_\alpha>0$, hence the action of $\Gamma_{\widetilde{W}}$ on $[\widetilde{W}_n]_{n\in\NN}$ has no elliptic element, except the stabilizers of the boundary 
components
of $\widetilde{W}$. Hence, according to Lemmas \ref{stabilizersdaretes} and \ref{dualite}, the minimal subtree for the
action of $\Gamma_{\widetilde{W}}$ is dual to a measured hyperbolic lamination $(\Lambda,\mu)$ (for any hyperbolic metric $m$ on $\Sigma$), of (the $m$-geometric realization
of) $W$. Moreover, for every
$\alpha\in\C(W)$, we have $i(\alpha,(\Lambda,\mu))=\epsilon_\alpha>0$, hence $(\Lambda,\mu)$ fills up $W$.\cqfd

\medskip

Let $T_{\min}$ be the minimal subtree of $[\widetilde{W}_n]_{n\in\NN}$ for the action of $\Gamma_{\widetilde{W}}$, and let $(\Lambda,\mu)$ be the measured hyperbolic 
lamination of (the $m$-geometric realization of) $W$ dual to $T_{\min}$.

By the compactness of the space of projective measured hyperbolic laminations, there exists a positive real sequence $(t_k)_{k\in\NN}$ such that the sequence $(t_k\alpha_k)_{k\in\NN}$
$\omega$-converges to a measured hyperbolic lamination
$(\Lambda_\alpha,\mu_\alpha)$ of $W$.

\blemmA The laminations $\Lambda$ and $\Lambda_\alpha$ are equal.
\elemmA

\demA Since $\lim_\omega t_k<+\infty$, we have \begin{align*}
	i((\Lambda,\mu),(\Lambda_\alpha,\mu_\alpha))&=\lim_\omega\; t_k\, i((\Lambda,\mu),\alpha_k)\\
	&=\lim_\omega t_k\epsilon_{\alpha_k}=0.
\end{align*}
Since $\Lambda$ is minimal and fills up $W$, we have $\Lambda=\Lambda_\alpha$.\cqfd

\blemmA\label{pointfixearbreminimal} Let $\gamma\in\Gamma_{\widetilde{W}}-\{e\}$. Then $\gamma$ has a fixed point in $T_{\min}$ if and only if $\gamma$ preserves a boundary component $\widetilde{c}$ of 
$\widetilde{W}$. This fixed point is unique, equal to $\star_{\widetilde{c}}$.
\elemmA

\demA For every $\alpha\in\C(W)$, we have $\lim_\omega\ell_{[q_n]}(\alpha)=\epsilon_\alpha>0$, and if $\alpha$ is the homotopy class of a boundary component of $W$, then
$\lim_\omega\ell_{[q_n]}(\alpha)=0$. Hence, an element $\gamma\in
\Gamma_{\widetilde{W}}-\{e\}$ has a fixed point in $[\widetilde{W}_n]_{n\in\NN}$ if and only if $\gamma$ preserves a boundary component $\widetilde{c}$ of 
$\widetilde{W}$. For all $n\in\NN$, let $\widetilde{c}_n$ be the boundary component of $\widetilde{W}_n$ corresponding to $\widetilde{c}$. Then, according to Remark
\ref{convergencedistancetranslation}, the set of fixed points of $\gamma$ in $T_{\widetilde{W}}$ is $[\widetilde{c}_n]_{n\in\NN}$. Moreover, the minimal subtree $T_{\min}$ is the union of the 
translation axes of the hyperbolic elements of $\Gamma_{\widetilde{W}}$ in $T_{\widetilde{W}}$, 
and according to Lemma \ref{compriss}, this union
intersects $[\widetilde{c}_n]_{n\in\NN}$ only at $\star_{\widetilde{c}}$.\cqfd

\subsection{Case of non degeneration.}\label{nondegeneration}

In this Subsection \ref{nondegeneration}, we assume as in the previous one that $W$ is the closure of a complementary connected component of $\Sigma_0$ which is neither a cylinder nor a pair of pants. 
%
%
We  consider now the case where the first conclusion of Lemma \ref{minorationuniforme} holds, i.e. there exists $\varepsilon>0$ such that for all
$\alpha\in\SSSS(W)$ and $\omega$-almost all $n\in\NN$,s
we have $\ell_{[q_n]}(\alpha)\geqslant\varepsilon$. The boundary components of $W$ (if any) are also some boundary components of $\Sigma_0$, hence their $[q_n]$-lengths 
$\omega$-converge to zero. Since $W$ has finitely many boundary components, there exists $I\in\omega$ such that for every boundary component $c$ of $W$ and all $n\in I$,
we have $\ell_{[q_n]}(c)\leqslant\frac{\varepsilon}{5}$. 
We replace the sequence $([q_n])_{n\in\NN}$ by $([q_n])_{n\in I}$.  Hence, in the  remainder of this Subsection \ref{nondegeneration}, we can assume that $\varepsilon>0$ is 
given so that:

\medskip
\noindent
$\bullet$ for all $\alpha\in\C(W)$, we have $\lim_\omega\ell_{[q_n]}(\alpha)<+\infty$;

\medskip
\noindent
$\bullet$ for all $n\in\NN$ and  $\alpha\in\C(W)$, we have $\ell_{[q_n]}(\alpha)\geqslant\varepsilon$;

\medskip
\noindent
$\bullet$  for every boundary component $c$ of $W$ and $n\in\NN$, we have $\ell_{[q_n]}(c)\leqslant\frac{\varepsilon}{5}$, and $\lim_\omega\ell_{[q_n]}(c)=0$.

\medskip

\rem This case corresponds to the subsurface $W_1$ of the example of Subsection \ref{decomposition}, since for every 
$\alpha\in\SSSS(W_1)$ we have $\lim_\omega\ell_{[q_n]}(\alpha)=\ell_{[q]}(\alpha)\geqslant 1$. 




\blemmA\label{vraisurface} For all $n\in\NN$, the boundary components of $\widehat{W}_n$ (if any) are simple and pairwise disjoint. Moreover, there exists no essential 
arc of $\widehat{W}_n$ between two (possibly equal) boundary components of $\widehat{W}_n$, of length less than $\frac{\varepsilon}{4}$.
\elemmA

\demA For all $n\in\NN$, in the construction of $\widehat{W}_n$ (see \cite{Raf05}), the boundary components of $\widehat{W}_n$ are simple. 
Assume for a contradiction that there exist two distinct boundary components that intersect each other, or that there exists  
an essential arc of length less than $\frac{\varepsilon}{4}$ between two (possibly equal) boundary components of $\widehat{W}_n$. Since $W$ is not a cylinder, 
the union of the image(s) 
of this (or these) boundary component(s)
with (possibly) the essential arc, contains a closed curve of length less than $\frac{2}{5}\varepsilon+\frac{2}{4}\varepsilon<\varepsilon$, which is not freely 
homotopic to 
a point nor to a boundary component of $\widehat{W}$, and hence is (freely homotopic to) an essential closed curve of $\widehat{W}$, a contradiction.\cqfd 

\medskip 

\remA This does not necessarily mean that the geodesic representatives of the boundary components of $W$ in $(\Sigma,[q_n])$ are simple and pairwise disjoint
for all $n\in\NN$.

\medskip

According to Lemma \ref{vraisurface}, the boundary components of $\widetilde{W}_n$ are pairwise at distance at least 
$\frac{\varepsilon}{4}$.
Let $H_{\widetilde{W}}$ be the subset of non trivial and non peripheral elements of $\Gamma_{\widetilde{W}}$: their action on $[\widetilde{W}_n]_{n\in\NN}$
is hyperbolic. For every $\gamma\in H_{\widetilde{W}}$,
let $F(\gamma)$ be the (possibly degenerated) flat strip,
union of all the
translation axes of $\gamma$ in $[\widetilde{W}_n]_{n\in\NN}$. We will see that the convex hull of $\bigcup_{\gamma\in H_{\widetilde{W}}}F(\gamma)$ contains 
all the points $\star_{\widetilde{c}}$, with $\widetilde{c}$ a boundary component of $\widetilde{W}$ (see Lemma \ref{star} for the definition). Let $\widetilde{W}_\omega$
be the convex hull of $\bigcup_{\gamma\in H_{\widetilde{W}}}F(\gamma)$ minus the points $\star_{\widetilde{c}}$, with $\widetilde{c}$ a boundary
component of $\widetilde{W}$ (if any). 

\blemmA\label{enveloppeconvexe}
The set $\widetilde{W}_\omega$ is a surface without boundary, endowed with a $\Gamma_{\widetilde{W}}$-invariant
half-translation structure.
%
\elemmA

\noindent
{\bf Proof of Lemma \ref{enveloppeconvexe}.} We will need the following lemmas.

\blemmA\label{22} The set $[\widetilde{W}_n]_{n\in\NN}-[\partial\widetilde{W}_n]_{n\in\NN}$ is a surface endowed with a $\Gamma_{\widetilde{W}}$-invariant Euclidean
metric with conical singularities
of
angles $k\pi$ with $k\in\NN$, $k\geqslant 3$.
\elemmA

\demA Let $[x_n]_{n\in\NN}$ be a point of $[\widetilde{W}_n]_{n\in\NN}-[\partial\widetilde{W}_n]_{n\in\NN}$.

\blemmA\label{intersectionorbites} There exist $I\in\omega$ and $\delta>0$ such that for all $n\in I$, the open ball $B_n(x_n,\delta)$ of radius $\delta$ contains at
most one point of any $\Gamma_{\widetilde{W}}$-orbit. 
\elemmA

\demA Since $[x_n]_{n\in\NN}$ does not belong to $[\partial\widetilde{W}_n]_{n\in\NN}$, there exist $I\in\omega$ and $\eta>0$ such that for all $n\in I$, the open ball 
$B_n(x_n,\eta)$ is contained in the interior of $\widetilde{W}_n$. Let $\delta=\min\{\frac{\varepsilon}{2},\frac{\eta}{3}\}$. Let $(z_n)_{n\in\NN}$ be a sequence in
$\widetilde{\Sigma}$ such that $z_n$ belongs to $B_n(x_n,\delta)$ for all $n\in\NN$, and let $\gamma\in\Gamma_{\widetilde{W}}-\{e\}$. Either $\gamma$ does not stabilize
any boundary
component of $\widetilde{W}_n$, then $d_n(z_n,\gamma z_n)\geqslant\varepsilon\geqslant 2\delta$ and $\gamma z_n$ does not belong to $B_n(x_n,\delta)$, 
or $\gamma$ preserves a boundary component $\widetilde{c}_n$ of $\widetilde{W}_n$, and since $d_n(z_n,\widetilde{c}_n)\geqslant\eta-\delta\geqslant 2\delta$,
according to Lemma \ref{pointfixe}, we have $d_n(z_n,\gamma z_n)\geqslant 2\delta$ and
$\gamma z_n$ does not belong to $B_n(x_n,\delta)$. Hence, for all $n\in I$, the ball $B_n(x_n,\delta)$ contains at most one point
of any $\Gamma_{\widetilde{W}}$-orbit in $\widetilde{\Sigma}$.\cqfd

\medskip

We fix $\delta$ as in Lemma \ref{intersectionorbites}. According to the Gauss-Bonnet formula, and since the angles at the singular points are at least $3\pi$,
the number of singular points in $\widehat{W}_n$, and their total
angles are uniformely bounded. Consequently, 
the surface 
$\widehat{W}_n$ has a constant number of singular points of constant angles, for 
$\omega$-almost all $n\in\NN$. Hence, there exists $k\in\NN$ such that, for $\omega$-almost all $n\in\NN$, 
the ball $B_n(x_n,\delta)$ contains $k$ singular points, denoted by $y_{1,n},\dots,y_{k,n}$ with angles $\theta_1,\dots,\theta_k$ respectively. 

Let $[z_n]_{n\in\NN}$ be a point of $B_\omega([x_n]_{n\in\NN},\frac{\delta}{2})$. Assume first that there exists $\eta>0$ such that $B_n(z_n,\eta)$ does not contain
any singular 
point for $\omega$-almost all $n\in\NN$. Then, the ball $B_{\omega}([z_n]_{n\in\NN},\frac{\eta}{2})$ is a Euclidean disk.

If there does not exist such a $\eta$, there exists
a sequence of singular points $(y_n)_{n\in\NN}$ with $\lim_\omega d_n(z_n,y_n)=0$, and according to Lemma \ref{intersectionorbites}, there exists 
$\I\subseteq\{1,\dots,k\}$ such that for any sequence of singular points $(y_n)_{n\in\NN}$, we have $\lim_\omega d_n(y_n,z_n)=0$ if and only if there exists $i\in\I$
and $I\in\omega$ such that $y_n=y_{i,n}$ for all $n\in I$.

\blemmA\label{convergencedecourbure}  The point $[z_n]_{n\in\NN}$ is then a singular point of angle $2\pi+\sum_{i\in\I}(\theta_i-2\pi)$.
\elemmA

\demA Let $I\in\omega$ be such that for any sequence of singular points $(y_n)_{n\in\NN}$, 
we have $\lim_\omega d_n(y_n,z_n)=0$ if and only if there exist $i\in\I$
such that $y_n=y_{i,n}$ for all $n\in I$. Since, for all $n\in I$, the ball $B_n(x_n,\delta)$ contains $k<+\infty$ singular points, there exist $\eta'>\eta>0$ small enough
such that $B_n(z_n,\eta')$ is contained in $B_n(x_n,\delta)$ and for $\omega$-almost all $n\in I$, the only singular points contained in $B_n(z_n,\eta')$ belongs to $B_n(z_n,\frac{\eta}{2})$.

Let $n\in I$. The boundary $\C_n$ of 
$B_n(z_n,\eta)$ is naturally endowed with a cyclic order. Let $x_0$ be a point of $\C_n$ and let $(x_i)_{i\in\{0,\dots,p_n\}}$, with $p_n\in\NN$,
be a finite  sequence,  increasing for the cyclic
order, defined by $d_n(x_i,x_{i+1})=r$, with $r\in[\frac{\eta}{4},\frac{\eta}{3}]$, for all $i\in\{0,\dots,p_n-1\}$, and such that $x_{p_n}$ is contained between 
$x_0$ and $x_1$. 
There exists $r\in[\frac{\eta}{4},\frac{\eta}{3}]$ such that $x_{k_n}=x_0$. Moreover, the ball
$B_n(z_n,\eta)$ contains at most one point of any $\Gamma_{\widetilde{W}}$-orbit, hence it embeds into $\widehat{W}_n$, and according to 
\cite[Lem.~4.1]{Minsky92},  for all $n\in I$, the length of $\C_n$ is at most $L\eta$, where $L$ depends on $\chi(\widehat{W})$. Hence, the integer $p_n$
is bounded by ${L\eta}/{\frac{\eta}{3}}=3L$, and there exists $p\in\NN^*$ such that it
is equal to $p\in\NN$ for $\omega$-almost all $n\in I$. Up to changing $I$, we can assume that for all $n\in I$, there exists a topological disk $P_n$
whose boundary
is a finite union of $p$ Euclidean segments of equal length between $\frac{\eta}{4}$ and $\frac{\eta}{3}$, whose endpoints belong to $\C_n$. Moreover, for all
$n\in I$, the boundary $\C_n$ is contained in $\overline{B}_n(z_n,\eta)-B_n(z_n,\frac{2\eta}{3})$, and by the choice of $\eta$, the distance $d_\omega$ 
is locally Euclidean on $[\overline{B}_n(z_n,\eta)-B_n(z_n,\frac{2\eta}{3})]_{n\in\NN}$, hence $[P_n]_{n\in\NN}$ is a topological disk whose boundary is a union of $p$
Euclidean segments. 

\begin{center}
	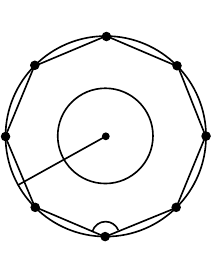
\end{center}

The interior Euclidean angle $\theta_n(x_i)$ at a vertex $x_i$ is determined by the distance between the midpoints of the segments adjacent to $x_i$,
and it $\omega$-converges to an angle $\theta_\omega(x_i)$ determined by the distance between the midpoints of the segments adjacent to the corresponding vertex. 
According to the Gauss-Bonnet formula applied to $P_n$, for all $n\in I$, we have $2\,\pi=\sum_{i\in\I}(2\pi-\theta_i)+\sum_{i=1}^p(\pi-\theta_n(x_i))$. Hence 
$\sum_{i=1}^p(\pi-\theta_n(x_i))=2\,\pi-\sum_{i\in\I}(2\pi-\theta_i)$. The sequence $(\sum_{i=1}^p(\pi-\theta_n(x_i)))_{n\in\NN}$ $\omega$-converges to
$\sum_{i=1}^p(\pi-\theta_\omega(x_i))$. Hence, $\sum_{i=1}^p(\pi-\theta_\omega(x_i))=2\,\pi-\sum_{i\in\I}(2\pi-\theta_i)$. Moreover, according to the remark preceding Lemma \ref{convergencedecourbure}, the distance $d_\omega$ is locally Euclidean on 
$[B_n(z_n,\eta)]_{n\in\NN}$, except at 
$[z_n]_{n\in\NN}$. Hence, the point $[z_n]_{n\in\NN}$ is a conical singularity of angle $\theta\geqslant0$, and by the Gauss-Bonnet formula applied to $[P_n]_{n\in\NN}$, we have 
$2\,\pi=2\,\pi-\theta-\sum_{i=1}^p(\pi-\theta_\omega(x_i))$, hence  $\theta=2\pi+\sum_{i\in\I}(\theta_i-2\pi)$.\cqfd

\medskip

Hence, the set $[\widetilde{W}_n]_{n\in\NN}-[\partial\widetilde{W}_n]_{n\in\NN}$ is a surface endowed with a locally Euclidean metric with conical singularities
of angles $k\pi$, $k\in\NN$ and $k\geqslant 3$, which is $\Gamma_{\widetilde{W}}$-invariant by naturality. This ends the proof of Lemma \ref{22}.\cqfd

\blemmA\label{pasdepointstordusdanslenveloppeconvexe} The convex hull of $\bigcup_{\gamma\in H_{\widetilde{W}}}F(\gamma)$, minus the points $\star_{\widetilde{c}}$, where 
$\widetilde{c}$ is a boundary component of $\widetilde{W}$ (if any), is a surface without boundary, endowed with a locally Euclidean metric with conical singularities
of angles of the form $k\pi$, $k\in\NN$ and $k\geqslant 3$, which is $\CAT(0)$ and $\Gamma_{\widetilde{W}}$-equivariant. Moreover,
the intersection of any translation axis of any
hyperbolic element of $\Gamma$ with $[\widetilde{W}_n]_{n\in\NN}$ is 
contained in the convex hull of $\bigcup_{\gamma\in H_{\widetilde{W}}}F(\gamma)$, if it is not empty.
\elemmA

\demA There exists a finite set of free homotopy classes (relative to the boundary of $W$) of simple arcs of $W$ joining some boundary components of $W$,
having pairwise disjoint representatives,
such that
no arc is freely homotopic to an arc contained
in a boundary component of $W$, and the union of some simple and pairwise disjoint representatives of these homotopy classes cuts $W$ into a finite set of disks.
Then, the union of all the lifts of these arcs in $\widetilde{W}$ cuts $\widetilde{W}$ into disks. We extend them into a set $\widetilde{S}$ of pairwise disjoint proper
biinfinite simple paths of $\widetilde{\Sigma}$ (possibly not in 
a $\Gamma$-equivariant way), that intersect the boundary of $\widetilde{W}$ in exactly two points, such that every path cuts $\widetilde{\Sigma}$ into two connected 
components.


Let $m$ be a hyperbolic metric on $\Sigma$ and let $\widetilde{m}$ be its pullback on $\widetilde{\Sigma}$.
We consider the $\widetilde{m}$-geodesic representatives of the elements of $\widetilde{S}$,
and we still denote by $\widetilde{S}$ the corresponding set of pairwise disjoint biinfinite $\widetilde{m}$-geodesics.
The elements of $\widetilde{S}$ are pairwise disjoint, each geodesic $\widetilde{s}\in\widetilde{S}$ divides $\widetilde{W}_{\widetilde{m}}$ into two connected
components, and the intersection of the union of the geodesics of $\widetilde{S}$
with $\widetilde{W}_{\widetilde{m}}$ cuts $\widetilde{W}_{\widetilde{m}}$ into disks. Since the geodesics are pairwise disjoint, the disks can be labelled by the
(unordered) sets of boundary components of
$\widetilde{W}_{\widetilde{m}}$ that meet theirs 
boundaries. We denote by $(D_k)_{k\in\NN}$ the closures of the connected complementary components of $(\bigcup\widetilde{S})\cap\widetilde{W}_{\widetilde{m}}$ in 
$\widetilde{W}_{\widetilde{m}}$, that are
(topological) disks. Let $n\in\NN$. For every 
$\widetilde{s}\in\widetilde{S}$, let $\widetilde{s}_n$ be a $[\widetilde{q}_n]$-geodesic having the same (ordered) pair of points at infinity than
$\widetilde{s}$. Since the geodesics $\widetilde{s}_n$ are properly homotopic to the pairwise disjoint elements of $\widetilde{S}$, 
to each disk $D_k$, with $k\in\NN$, corresponds a unique convex generalized disk $D_k(n)$ (i.e. a (possibly empty) disk with a finite number of spikes of finite lengths
(possibly zero)), whose boundary is contained in 
$(\bigcup_{\widetilde{s}\in\widetilde{S}}\widetilde{s}_n)\cup\partial\widetilde{W}_n$. It may happen that $D_k(n)$ is a graph. 

Let $k_1, k_2\in\NN$ and let $(x_{1,n})_{n\in\NN}$ and $(x_{2,n})_{n\in\NN}$ be two sequences of points such that for all $n\in\NN$, $x_{1,n}$ belongs to $D_{k_1}(n)$
and
$x_{2,n}$ belongs to $D_{k_2}(n)$. Since every geodesic of $\widetilde{S}$ cuts $\widetilde{W}_{\widetilde{m}}$ into two connected components,  
if the point $x_{1},x_{2}$ belong respectively to $D_{k_1}$ and $D_{k_2}$, there exists a finite sequence $(\widetilde{s}_i)_{i\in I}$ in
$\widetilde{S}$, such that the $\widetilde{m}$-geodesic segment joining $x_1$ to $x_2$ cuts the geodesics  $(\widetilde{s}_i)_{i\in I}$. Then, for
all $n\in\NN$, the geodesics $(\widetilde{s}_{i,n})_{i\in I}$ cut the geodesic segment joining $x_{1,n}$ to $x_{2,n}$ in $\widetilde{W}_n$ into a bounded number of
geodesic
segments of finite lengths (several
geodesics can be cut at the same time). 
Hence, the geodesic segment $[[x_{1,n},x_{2,n}]]_{n\in\NN}$ in $[\widetilde{W}_n]_{n\in\NN}$ is the concatenation of a finite 
number of geodesic segments whose endpoints belong to $\bigcup_{i\in I}([\widetilde{s}_{i,n}]_{n\in\NN})$, hence it is contained in 
$\bigcup_{k\in\NN}[D_k(n)]_{n\in\NN}$. Therefore $\bigcup_{k\in\NN}[D_k(n)]_{n\in\NN}$ is convex.


Let us show that the union $\bigcup_{k\in\NN}[D_k(n)]_{n\in\NN}$ intersects $[\partial\widetilde{W}_n]_{n\in\NN}$ only at the points $\star_{\widetilde{c}}$,
where ${\widetilde{c}}$ is a boundary component of $\widetilde{W}$. For all $n\in\NN$, if $\widetilde{s}_n$ is a geodesic corresponding to an element $\widetilde{s}\in\widetilde{S}$,
it is interlaced with two boundary components 
$\widetilde{c}$ and $\widetilde{c}'$ of $\widetilde{W}$, and its pair of points at infinity does not depend on $n$. 
Hence, there exists $\gamma\in\Gamma-\{e\}$ and $k\in\ZZ$ such that for all $n\in\NN$, a translation axis $\Ax_n(\gamma)$
of $\gamma$ is interlaced with $\widetilde{c}_n$, and $\widetilde{s}_n$ is caught between $\gamma_{\widetilde{c}}^{-k}(\Ax_n(\gamma))$ and 
$\gamma_{\widetilde{c}}^k(\Ax_n(\gamma))$, where $\gamma_{\widetilde{c}}\in\Gamma_{\widetilde{W}}-\{e\}$ preserves $\widetilde{c}$. 
Similarly, there exists $\gamma'\in\Gamma-\{e\}$ and $k'\in\ZZ$ such that for all $n\in\NN$, a translation axis $\Ax_n(\gamma')$
of $\gamma'$ is interlaced with $\widetilde{c}'_n$, and $\widetilde{s}_n$ is caught between $\gamma_{\widetilde{c}'}^{-k'}(\Ax_n(\gamma'))$ and
$\gamma_{\widetilde{c}'}^{k'}(\Ax_n(\gamma'))$, where $\gamma_{\widetilde{c}'}\in\Gamma_{\widetilde{W}}-\{e\}$ preserves $\widetilde{c}'$.
According to Lemma \ref{star}, 
the geodesic $[\widetilde{s}_n]_{n\in\NN}$ can intersect $[\partial\widetilde{W}_n]_{n\in\NN}$ only at the points $\star_{\widetilde{c}}$ and $\star_{\widetilde{c}'}$.
Hence, the convex union $\bigcup_{k\in\NN}[D_k(n)]_{n\in\NN}$,
minus 
the points $\star_{\widetilde{c}}$, is contained in $[\widetilde{W}_n]_{n\in\NN}-[\partial\widetilde{W}_n]_{n\in\NN}$, and according to Lemma \ref{22}, it is a 
surface without boundary. Since it is convex, it is $\CAT(0)$ as $[\widetilde{W}_n]_{n\in\NN}$.
Finally, it is naturally $\Gamma_{\widetilde{W}}$-invariant.

\medskip

Let $\gamma\in H_{\widetilde{W}}$. By the definition of $\widetilde{S}$, there exists a sequence $(\widetilde{s}_{i})_{i\in\ZZ}$ of $\widetilde{S}$
that cuts the translation axis $\Ax_{\widetilde{m}}(\gamma)$ of $\gamma$ in $\widetilde{W}_{\widetilde{m}}$
into bounded intervals. Then, for all $n\in\NN$, if $\Ax_{n}(\gamma)$ is a translation axis of $\gamma$ in $\widetilde{W}_n$, it is cut by the sequence of 
geodesics $(\widetilde{s}_{i,n})_{i\in\NN}$ into bounded intervals, and the order of the 
$(\widetilde{s}_i)_{i\in\NN}$ does not depend on $n\in\NN$ (it may happen that several geodesics are intersected at the same time).
Hence, the geodesic $[\Ax_n(\gamma)]_{n\in\NN}$ is contained in $\bigcup_{k\in\NN}[D_k(n)]_{n\in\NN}$.
Hence, all the (possibly degenerated) flat strips $F(\gamma)$, with $\gamma\in H_{\widetilde{W}}$, are contained in $\bigcup_{k\in\NN}[D_k(n)]_{n\in\NN}$, and since
$\bigcup_{k\in\NN}[D_k(n)]_{n\in\NN}$ is convex, the convex hull of $\bigcup_{\gamma\in H_{\widetilde{W}}}F(\gamma)$ is contained in $\bigcup_{k\in\NN}[D_k(n)]_{n\in\NN}$. 
Moreover, since the convex hull of $\bigcup_{\gamma\in H_{\widetilde{W}}}F(\gamma)$ is $\Gamma_{\widetilde{W}}$-equivariant, according to Lemma
\ref{pointfixe}, it contains all the points $\star_{\widetilde{c}}$, where $\widetilde{c}$ is a boundary component of $\widetilde{W}$, hence it contains all the geodesic
segments joining them, and notably the segments $[\widetilde{s}_n]_{n\in\NN}\cap[\widetilde{W}_n]_{n\in\NN}$, with
$\widetilde{s}\in\widetilde{S}$.
Hence, the convex hull of $\bigcup_{\widetilde{s}\in\widetilde{S}}[\widetilde{s}_n]_{n\in\NN}\cap[\widetilde{W}_n]_{n\in\NN}$, which is equal to 
$\bigcup_{k\in\NN}[D_k(n)]_{n\in\NN}$, is equal to the convex hull of $\bigcup_{\gamma\in H_{\widetilde{W}}}F(\gamma)$.

\medskip

Finally, if $\gamma$ is a hyperbolic element of $\Gamma-\Gamma_{\widetilde{W}}$ whose translation axis $\Ax_\omega(\gamma)$ in $\lim_\omega\revetn$ intersects $[\widetilde{W}_n]_{n\in\NN}$, according to Lemma \ref{star}, 
there exists two boundary components
$\widetilde{c}$ and $\widetilde{c}'$ of $\widetilde{W}$ such that the intersection $\Ax_\omega(\gamma)\cap[\widetilde{W}_n]_{n\in\NN}$ is the segment $[\star_{\widetilde{c}},\star_{\widetilde{c}'}]$, hence
is contained in
the convex hull of $\bigcup_{\gamma\in H_{\widetilde{W}}}F(\gamma)$.\cqfd

\blemmA\label{extension} The locally Euclidean metric with conical singularities on $\widetilde{W}_\omega$ is a half-translation structure $[\widetilde{q}_\omega]$ on $\widetilde{W}_\omega$, that
is the pullback of a half-translation structure 
$[\widehat{q}_\omega]$ on a finite type surface $\widehat{W}_\omega$ homeomorphic to $\widehat{W}-\partial\widehat{W}$. Moreover, if $W\not=\Sigma$,
$[\widehat{q}_\omega]$ 
can be extended to the compact surface obtained from $\widehat{W}_\omega$ by filling in the punctures, with possibly some singular points of angle $\pi$ at the added points.
\elemmA

\demA Let $x=[x_n]_{n\in\NN}$ be a point of $\widetilde{W}_\omega$ such that there exists $r>0$ small enough, such that the (Euclidean) disk $D_\omega(x,2r)$
does not contain any singular point. Then, there exists $I$ in $\omega$ such that for all $n\in I$, the disk $D_n(x_n,r)$ in 
$\widetilde{W}_n$ does not contain any
singular point of $[\widetilde{q}_n]$. Let $x_1=[x_{1,n}]_{n\in\NN}$ and $x_2=[x_{2,n}]_{n\in\NN}$ be two points of $D_\omega(x,r)$, and let $r_1,\,r_2>0$ be
such that the disks
$D_\omega(x_1,r_1)$ and $D_\omega(x_2,r_2)$ are contained in $D_\omega(x,r)$. We can always assume that the distance $d_n(x_{1,n},x_{2,n})$ is constant for 
all $n\in\NN$.
For all $n\in I$, let $\psi_{i,n}:D_n(x_{i,n},r_i)\to\DD(0,r_i)$ 
(with $i=1,2$, where $\DD(0,r_i)=\{z\in\CC\,:\,|z|<r_i\}$) be the inverse of an exponential map at $x_{i,n}$.    
Then, if $r_1+r_2<d_\omega(x_1,x_2)$, for $\omega$-almost all $n\in\NN$, there exists $c_n\in\CC$ such that the restriction of 
$\psi_{2,n}\circ\psi_{1,n}^{-1}$ to $\psi_{1,n}(D_n(x_{1,n},r_1)\cap
D_n(x_{2,n},r_2))$ is either $z\mapsto z+c_n$ or $z\mapsto -z+c_n$. Moreover, we have $|c_n|=d_n(x_{1,n},x_{2,n})$, hence the sequence 
$(c_n)_{n\in\NN}$ $\omega$-converges to $c\in\CC$.
Let us define the isometry $\psi_{i,\omega}:D_\omega(x_i,r_i)\to\DD(0,r_i)$  
by 
$\psi_{i,\omega}([w_n]_{n\in\NN})=\lim_\omega\psi_{i,n}(w_n)$ ($i\in\{1,2\}$). Then, the restriction of
$\psi_{2,\omega}\circ\psi_{1,\omega}^{-1}$ to $\psi_{1,\omega}(D_\omega(x_1,r_1)\cap
D_\omega(x_2,r_2))$ is either $z\mapsto z+c$ or $z\mapsto -z+c$. 

Hence, the set of maps thus defined spans a maximal atlas of charts on the complement of the singular points in $\widetilde{W}_\omega$, such that the exchange maps are
of the form $z\mapsto\pm z+c$, with $c\in\CC$. Since the singular points 
are of angles $k\pi$, with $k\in\NN$ and $k\geqslant3$, it defines a half-translation structure $[\widetilde{q}_{\omega}]$ on $\widetilde{W}_\omega$.
This half-translation structure is 
naturally $\Gamma_{\widetilde{W}}$-invariant, and since $\widetilde{W}_\omega$ is simply connected and the action of $\Gamma_{\widetilde{W}}$ is 
properly discontinuous, it defines a half-translation structure
$[\widehat{q}_{\omega}]$ on the surface 
$\widehat{W}_\omega=\Gamma_{\widetilde{W}}\backslash\widetilde{W}_\omega$. Moreover, since $\Gamma_{\widetilde{W}}$ is finitely generated, the surface $\widehat{W}_\omega$
is of finite type. Let us prove that $\widehat{W}_\omega$ is homeomorphic to $\widehat{W}-\partial\widehat{W}$. Since $\Gamma_{\widetilde{W}}$ has $b$ conjugation classes of parabolic elements ($b$ is
the number of boundary components of $\widehat{W}$), the surface $\widehat{W}_\omega$ has $b$ punctures. Moreover, the surface $\widetilde{W}_\omega$ is orientable,
and 
two translation axes of two elements $\gamma_1$ and $\gamma_2$ of
$H_{\widetilde{W}}$ in 
$\widetilde{W}_\omega$ are properly homotopic to two disjoint biinfinite paths if and only if  any two translation axes of $\gamma_1$ and $\gamma_2$ in
$\widetilde{W}_n$ are properly homotopic to two disjoint biinfinite paths, for any $n\in\NN$ (see the proof of \cite[Lem.~6.26]{Morzy4}). 
Hence, the cardinality of a maximal set of isotopy classes of  essential simple closed curves of $\widehat{W}_\omega$ having pairwise disjoint representatives 
is equal to the cardinality of a
maximal set of isotopy classes of simple closed curves having pairwise disjoint representatives of $\widehat{W}$.  
Hence, the Euler characteristics of $\widehat{W}_\omega$ and of $\widehat{W}$ are equal, and since $\widehat{W}_\omega$ has $b$ punctures,
the surfaces $\widehat{W}_\omega$ and $\widehat{W}$ have the same genus.
Hence, $\widehat{W}_\omega$ is homeomorphic to $\widehat{W}-\partial\widehat{W}$.

The complex structure defined by $[\widehat{q}_\omega]$ on $\widehat{W}_\omega$ can be extended to the compact surface without puncture obtained from $\widehat{W}_\omega$
by filling in the punctures. Let $a$ be a puncture of $\widehat{W}_\omega$, and let $(U,z)$ be a chart of the extended complex structure, with $a\in U$. 
Let $\widehat{q}_\omega\in\Q(\widehat{W}_\omega)$ be a representative of $[\widehat{q}_\omega]$. We define the holomorphic
map $\varphi$ on $U-\{a\}$ by $\varphi(x)=\widehat{q}_{\omega,1}(x)(\frac{dz_1}{dz})^2$, for any chart $(U_1,z_1)$ of $\widehat{W}_\omega$ such that $x$ belongs to 
$U_1$, where $\widehat{q}_{\omega,1}$ is the representative of $\widehat{q}_\omega$ in $(U_1,z_1)$.
Let $\widetilde{a}$ be the fixed point of an elliptic element $\gamma_{\widetilde{a}}\in\Gamma_{\widetilde{W}}-\{e\}$ in the completion of $\widetilde{W}_\omega$, such that the 
free homotopy class of closed curves defined by $\gamma_{\widetilde{a}}$ in $\widehat{W}_\omega$ as some representatives contained in any neighborhood of $a$.
Let $x$ be a point of $\widetilde{W}_\omega$, close enough to ${\widetilde{a}}$
(for the distance $d_\omega$). Then, by the isoperimetric inequality (see for example \cite[Thm.~2.17~p.~426]{BriHae99}), the angle sector,
bounded by $[{\widetilde{a}},x]$, and $[{\widetilde{a}},\gamma_{\widetilde{a}}x]$ has a finite area.
Moreover, it projects to a closed neighborhood of $a$ in $\widehat{W}_\omega$, whose $[\widehat{q}_\omega]$-area is finite 
since the cover projection does not increase the area.
Hence, the integral of $|\varphi|$ on a closed neighborhood of $a$ is finite, and $\varphi$ can be extended to a meromorphic function, with at most a simple pole
in $a$.
Hence, the half-translation structure $[\widehat{q}_\omega]$ on the compact surface with punctures $\widehat{W}_\omega$ can be extended to a half-translation structure on
the compact surface, possibly with singularities of angle $\pi$ at the punctures.  \cqfd

\subsection{Complementary components of $\Sigma_0$ which are cylinders or pair of pants.}\label{cylindresetpantalons}

It remains to study the ultralimits of the geometric realizations of the connected components of the preimages in $\widetilde{\Sigma}$ of the  connected components of
$\Sigma_0$ and of
$\Sigma-\Sigma_0$ which are
pair of pants or cylinders, and of the cylinders that can be homotoped to a boundary component of $\Sigma_0$. 

\medskip 

Let $\alpha$ be the free homotopy class of closed curves of a boundary component of $\Sigma_0$, and let $\widetilde{\alpha}$ be a lift of $\alpha$ in
$\widetilde{\Sigma}$. For all $n\in\NN$, we denote by $F_n(\widetilde{\alpha})$ the
(possibly degenerated) flat strip, union of all the geodesics of $\revetn$ having the same pair of points at infinity than $\widetilde{\alpha}$.
The ultralimit $[F_n(\widetilde{\alpha})]_{n\in\NN}$ exists and can be neither a plane nor a half-plane, since the height of $F_n(\widetilde{\alpha})$ is bounded, hence it is either a flat strip or a 
single geodesic. According to Remark \ref{convergencedistancetranslation}, it is the set of fixed points in 
$\lim_\omega\revetn_{n\in\NN}$ of the elements of the stabilizer of $\widetilde{\alpha}$ in $\Gamma$, and according to
Lemma \ref{star}, the intersection of any translation axis of any hyperbolic element of $\Gamma$ with $[F_n(\widetilde{\alpha})]_{n\in\NN}$ is either empty, or a point, or a geodesic segment 
orthogonal to the boundary components of the flat strip $[F_n(\widetilde{\alpha})]_{n\in\NN}$.

\medskip

Assume until the end of this section that $W$ is the closure of a connected component of $\Sigma-\Sigma_0$ which is a pair of pants.

\blemmA\label{pantalon}

The ultralimit $[\widetilde{W}_n]_{n\in\NN}$ is an 
$\RR$-tree and the action of $\Gamma_{\widetilde{W}}$ on $[\widetilde{W}_n]_{n\in\NN}$ has a global fixed point $\star_{\widetilde{W}}$
which is equal to $\star_{\widetilde{c}}$ for every boundary component $\widetilde{c}$ of $\widetilde{W}$ (see the lines following the statement of Lemma \ref{star} 
for the definition).
\elemmA

\demA Let $\widehat{c}_1$, $\widehat{c}_2$ and $\widehat{c}_3$ be the homotopy classes of the boundary components of $\widehat{W}$. Let $n\in\NN$.
Assume for a contradiction
that the boundary components of $\widehat{W}_n$ are pairwise disjoint. According to \cite[Rem.~3.2]{Raf05}, the interior curvature of any boundary component of $\widehat{W}_n$ 
is 
at most $-\pi$. Hence, according to the Gauss-Bonnet formula, we should have $2\pi\chi(\widehat{W})\leqslant-3\pi$, which is impossible since $\chi(\widehat{W})=-1$.
Hence, up to changing the notation, for $\omega$-almost all $n\in\NN$, we can assume that the boundary components $\widehat{c}_{1,n}$ and $\widehat{c}_{2,n}$ 
of $\widehat{W}_n$, corresponding to $\widehat{c}_1$ and $\widehat{c}_2$, intersect each other.
The geodesic $\widehat{c}_{3,n}$ either intersects the union of the other ones, or is disjoint from the union of the other ones, and the interior of $\widehat{W}_n$
is a cylinder. In the second case, the union $\widehat{c}_{1,n}\cup\widehat{c}_{2,n}$ contained a simple closed curve isotopic to $\widehat{c}_{3,n}$, and 
according to Lemma \ref{cylinder}, there exists a point in $\widehat{c}_{1,n}\cup\widehat{c}_{2,n}$ at distance of $\widehat{c}_{3,n}$
at most $(\ell_{[\widehat{q}_n]}(\widehat{c}_{1,n})+\ell_{[\widehat{q}_n]}(\widehat{c}_{2,n}))$, and $\widehat{W}_n$ can be cut into a disk of perimeter at most 
$2(\ell_{[\widehat{q}_n]}(\widehat{c}_{1,n})+\ell_{[\widehat{q}_n]}(\widehat{c}_{2,n}))
+(\ell_{[\widehat{q}_n]}(\widehat{c}_{1,n})+\ell_{[\widehat{q}_n]}(\widehat{c}_{2,n}))+\ell_{[\widehat{q}_n]}(\widehat{c}_{3,n})\leqslant 
4(\ell_{[\widehat{q}_n]}(\widehat{c}_{1,n})+\ell_{[\widehat{q}_n]}(\widehat{c}_{2,n}))$. 
In both cases, the diameter of $\widehat{W}_n$ is at most $4(\ell_{[\widehat{q}_n]}(\widehat{c}_{1,n})+\ell_{[\widehat{q}_n]}(\widehat{c}_{2,n}))$.
According to Lemma \ref{arbrereel}, the 
ultralimit $[\widetilde{W}_n]_{n\in\NN}$ is an $\RR$-tree. Moreover, there exist some homotopy classes of arcs joining the boundary components of $\widehat{W}$, whose
$[\widehat{q}_n]$-lengths $\omega$-converge to zero. As in the proof of Lemma \ref{elliptic}, the action of $\Gamma_{\widetilde{W}}$
has a global fixed point $\star_{\widetilde{W}}$ which is equal to $\star_{\widetilde{c}}$ for every boundary component $\widetilde{c}$ of $\widetilde{W}$.\cqfd

\section{Mixed structures.}\label{mixedstructure}

Let $\Sigma$ be a compact, connected, orientable surface, such that $\chi(\Sigma)<0$, and let $p:\widetilde{\Sigma}\to\Sigma$ be a universal
cover with covering group $\Gamma$. Let $\Sigma_0$ be a tight proper subsurface of $\Sigma$, and let $\widetilde{\Sigma}_0$ be its preimage in $\widetilde{\Sigma}$. 
If $\widetilde{W}$ is a connected component of $\widetilde{\Sigma}_0$ or of $\widetilde{\Sigma}-\widetilde{\Sigma}_0$, 
we denote by $\Gamma_{\widetilde{W}}$ the stabilizer of $\widetilde{W}$ in $\Gamma$. 

\medskip

Let $W$ be  a connected component of $\Sigma-\Sigma_0$ and let $\widetilde{W}$ be a connected component of the preimage of
$W$ in $\widetilde{\Sigma}$. Assume first that $W$ is neither a cylinder nor a pair of pants. Let $T_{\Lambda}$ be the $\RR$-tree dual  to a 
measured hyperbolic lamination $(\Lambda,\mu)$ on $W$ (for any complete hyperbolic metric).

\blemmA\label{pointfixelamination}(see \cite[Lem.~7.1]{Morzy4}) If $\Lambda$ is filling, the stabilizers of the boundary components of $\widetilde{W}$
have a unique fixed point in $T_\Lambda$. The other elements have no fixed point.
\elemmA

The surface $W$ is a finite type surface. Let $[q_W]$ be a half-translation structure on $W$, that can be extended
to the compact surface obtained from $W$ by filling in the punctures, with possibly some singularities of angle $\pi$ at the added point.
Let $p_W:(\widetilde{W},[\widetilde{q}_{W}])\to(W,[q_W])$ be a universal cover. Then
$(\widetilde{W},[\widetilde{q}_W])$ is not complete.

\blemmA\label{pointfixestructureplate} The completion $\widetilde{W}^c$ of $(\widetilde{W},[\widetilde{q}_W])$ is the union of $\widetilde{W}$ and of countably many isolated points which
are exactly the 
fixed points of the stabilizers of the boundary components of $\widetilde{W}$ in $\Gamma_{\widetilde{W}}$.
\elemmA


Let $\widetilde{\Sigma}_0$ be as above. Let $\widetilde{W}$ be a connected component of $\widetilde{\Sigma}_0$ or of $\widetilde{\Sigma}-\widetilde{\Sigma}_0$. 
Let $X_{\widetilde{W}}$ be a complete geodesic metric space endowed with an isometric action of the stabilizer $\Gamma_{\widetilde{W}}$ of $\widetilde{W}$ in $\Gamma$,
such that: 

\medskip
\noindent
$\bullet$~ if $\widetilde{W}$ is a strip, $X_{\widetilde{W}}$ is empty;

\medskip
\noindent
$\bullet$~ if $\widetilde{W}$ is a connected component of $\widetilde{\Sigma}_0$, which is not a strip, $X_{\widetilde{W}}$ is a point, and the action of
$\Gamma_{\widetilde{W}}$ on $X_{\widetilde{W}}$ is trivial;  

\medskip
\noindent
$\bullet$~ if $\widetilde{W}$ is a connected component of $\widetilde{\Sigma}-\widetilde{\Sigma}_0$ such that $W=p(\widetilde{W})$ is neither a 
cylinder nor a pair of pants, $X_{\widetilde{W}}$ is

\medskip

-either the $\RR$-tree dual to a filling measured
hyperbolic
lamination $(\Lambda,\mu)$ of $W$ (for any complete hyperbolic metric), and then the action of $\Gamma_{\widetilde{W}}$ on $X_{\widetilde{W}}$ is dual to $(\Lambda,\mu)$,

\medskip

-or $X_{\widetilde{W}}$ is the completion of a universal cover of a half-translation structure on $W$, that can be extended 
to a half-translation structure on the compact surface obtained from $W$ by filling in the punctures, with possibly some
singularities of angle $\pi$ at the added points, and then the action of $\Gamma_{\widetilde{W}}$ on $X_{\widetilde{W}}$ is the covering action, extended at the added points;

\medskip
\noindent
$\bullet$~ if $\widetilde{W}$ is a connected component of $\widetilde{\Sigma}-\widetilde{\Sigma}_0$ such that $W=p(\widetilde{W})$ is 
a pair of pants, $X_{\widetilde{W}}$ is a point, and the action of $\Gamma_{\widetilde{W}}$ on $X_{\widetilde{W}}$ is trivial.

\medskip

Finally, for every proper homotopy class $\widetilde{c}$ of boundary components of $\widetilde{\Sigma}_0$ (two boundary components of $\widetilde{\Sigma}_0$ can be 
properly homotopic because of the strips), let $X_{\widetilde{c}}$ be an edge 
(i.e. a compact interval of $\RR$, possibly reduced to a point). We assume that if $\widetilde{c}$ bounds a connected component of
$\widetilde{\Sigma}-\widetilde{\Sigma}_0$ which
is a strip, the length of this edge is nonzero and if $\widetilde{W}$ is a connected component of $\widetilde{\Sigma}-\widetilde{\Sigma}_0$ whose image in $\Sigma$ is a 
pair of pants, then there exists at least one boundary component $\widetilde{c}$ of $\widetilde{W}$ such that the length of $X_{\widetilde{c}}$ is nonzero.
We endow $X_{\widetilde{c}}$ with the trivial action of the stabilizer $\Gamma_{\widetilde{c}}$ of $\widetilde{c}$ in $\Gamma$. 
Moreover, if $\widetilde{W}$ and $\widetilde{W}'$ are the connected components of $\widetilde{\Sigma}_0$  or of $\widetilde{\Sigma}-\widetilde{\Sigma}_0$ bounded by 
$\widetilde{c}$, we identify the fixed points of $\Gamma_{\widetilde{c}}$ in $X_{\widetilde{W}}$ and $X_{\widetilde{W}'}$
(which exist and are unique, according to Lemmas \ref{pointfixelamination} and
\ref{pointfixestructureplate}) with the (possibly non distinct) endpoints of $X_{\widetilde{c}}$.

\begin{center}
 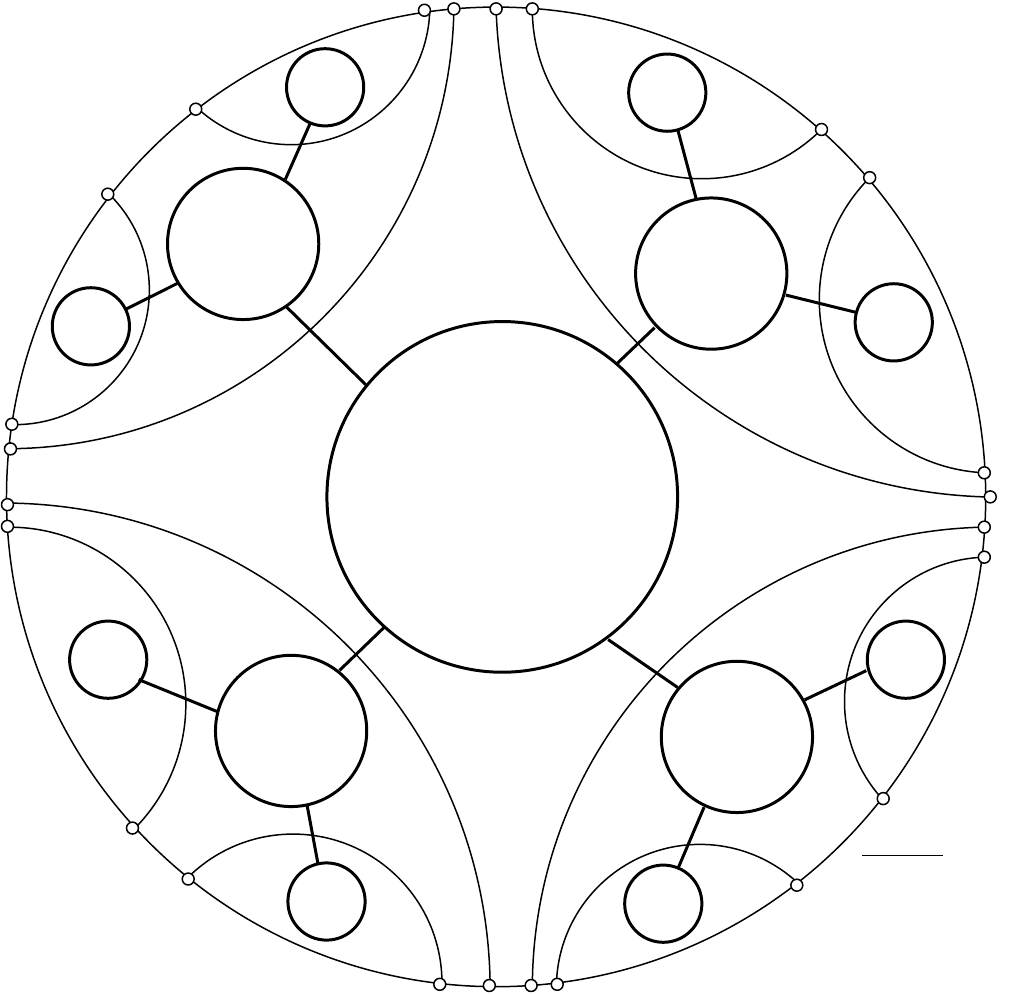
\end{center}

Let $(X,d)$ be the topological space obtained by this gluing, endowed with
the length distance induced by the distances on the different metric spaces. We identify the metric spaces $X_{\widetilde{W}}$,
and $X_{\widetilde{c}}$ with their images in $X$. The actions of the stabilizers $\Gamma_{\widetilde{W}}$  above extend uniquely in an isometric
action of $\Gamma$ by setting $\gamma(X_\star)=X_{\gamma\star}$ for every $\gamma\in\Gamma$, where $\star$ is a connected component of
$\widetilde{\Sigma}_0$ or of $\widetilde{\Sigma}-\widetilde{\Sigma}_0$, or a boundary component of $\widetilde{\Sigma}_0$. If $\Gamma$ acts on a $\CAT(0)$ geodesic space $X$, and if $\Gamma'$ is
a subgroup of $\Gamma$, if there exists a unique minimal non empty closed convex subset of $X$ which is preserved by $\Gamma'$, it is called the
{\it convex core} of $\Gamma'$. By definition of the action of $\Gamma$ on $X$, we see that if $\widetilde{W}$ is a connected component of $\widetilde{\Sigma}_0$ or of 
$\widetilde{\Sigma}-\widetilde{\Sigma}_0$, then the convex core of $\Gamma_{\widetilde{W}}$ in $X$ is $X_{\widetilde{W}}$, and if $\widetilde{c}$ is a boundary component
of $\widetilde{\Sigma}_0$, the set of fixed points of $\Gamma_{\widetilde{c}}$ in $X$ is $X_{\widetilde{c}}$. The elements of $\Gamma$ that have a fixed point in $X$ are 
those that preserve a connected component of $\widetilde{\Sigma}_0$ or a connected component of $\widetilde{\Sigma}-\widetilde{\Sigma}_0$ whose image in $\Sigma$ is a 
pair of pants. The other elements are hyperbolic. Moreover, if $\gamma\in\Gamma$ is a hyperbolic element in $X$, it can have more than one 
translation axis only if it preserves a piece which is the completion of a surface endowed with a half-translation structure. 


\blemmA\label{biendefini} The space $(X,d)$ just defined is a $\CAT(0)$ tree-graded space, whose set of pieces are the spaces $X_{\widetilde{W}}$ and $X_{\widetilde{c}}$ 
as above, which is uniquely determined by the half-translation structures and by 
the measured hyperbolic laminations on the connected components
of $\Sigma-\Sigma_0$, and by the lengths of the edges corresponding to the boundary components of $\widetilde{\Sigma}_0$, up to $\Gamma$-equivariant isometry.
Moreover, the action of $\Gamma$ on $X$ determines uniquely
the subsurface $\Sigma_0$ (up to isotopy), the measured hyperbolic laminations (up to isotopy) and the half-translation structures (up to isometry), on the connected
components of $\Sigma-\Sigma_0$.  
\elemmA

\remA Recall that $\Sigma_0$ is proper in $\Sigma$, and by construction, the space $X$ is not reduced to a point,
and the action of $\Gamma$ on $X$ has no global fixed point.

\medskip

\demA According to the construction, two distinct pieces $Y$ and $Y'$ have at most one common point $x$ which is then a fixed point of the stabilizer
$\Gamma_{\widetilde{c}}$ of a 
boundary component $\widetilde{c}$ of $\widetilde{\Sigma}_0$ in $\Gamma$. Moreover, the boundary component $\widetilde{c}$ cuts $\widetilde{\Sigma}$ into two
connected components, and if $\widetilde{W}$ and $\widetilde{W}'$ are two connected components of $\widetilde{\Sigma}_0$ or of $\widetilde{\Sigma}-\widetilde{\Sigma}_0$ that 
are not contained in the same connected component of $\widetilde{\Sigma}-\widetilde{c}$, then there exists no common point between the pieces $X_{\widetilde{W}}$ and 
$X_{\widetilde{W}'}$, except possibly $x$. Hence, any simple loop of $X$ is contained in a unique piece, and since all the pieces are simply connected, the space $X$ is
simply connected. Moreover, it is clearly locally $\CAT(0)$ and hence globally $\CAT(0)$. And any simple geodesic triangle is contained in a unique piece,
hence $(X,d)$ is a tree-graded space.
%

\medskip

We now show that the action of $\Gamma$ on $(X,d)$ determines $\Sigma_0$, up to isotopy. Let $Z$ be a non trivial $\pi_1$-injective closed connected subsurface of
$\Sigma$ and let $\widetilde{Z}$ be a connected component of
the preimage of $Z$ in $\widetilde{\Sigma}$. Assume that
the stabilizer $\Gamma_{\widetilde{Z}}$ of $\widetilde{Z}$ in $\Gamma$ is neither trivial nor cyclic, and preserves a unique 
point $\star_{\widetilde{Z}}$ in $X$. We can always assume that $Z$ is not a proper subsurface of a bigger subsurface having this property.   
Then, all the elements of $\Gamma_{\widetilde{Z}}$ preserve a connected components of $\widetilde{\Sigma}_0$ or of $\widetilde{\Sigma}-\widetilde{\Sigma}_0$ whose image in
$\Sigma$ is a pair of pants. Hence ${Z}$ is either a connected component of $\Sigma_0$ or the closure of a connected component of $\Sigma-\Sigma_0$ which is a pair of pants. 
Since no connected component of $\Sigma_0$ is a pair of pants, in the first case, the group $\Gamma_{\widetilde{Z}}$ is either a free group of rank at least $3$,
or if $Z$ is a 
torus minus a disk, $\Gamma_{\widetilde{Z}}$ is a free group of rank $2$. In that case, all the boundary components of $\widetilde{Z}$ belong to the same 
$\Gamma_{\widetilde{Z}}$-orbit. Moreover, if $\widetilde{c}$ is a boundary component of $\widetilde{\Sigma}_0$ and $A$ is a connected component of
$\widetilde{\Sigma}-\widetilde{c}$, then the image of the union of the pieces corresponding to the connected components of $\widetilde{\Sigma}_0$ or of 
$\widetilde{\Sigma}-\widetilde{\Sigma}_0$, contained in $A$, by $\gamma\in\Gamma_{\widetilde{Z}}$, is the union of the pieces corresponding to the connected components of $\widetilde{\Sigma}_0$ or of 
$\widetilde{\Sigma}-\widetilde{\Sigma}_0$, contained in $\gamma(A)$. Hence, in that case, the action of $\Gamma_{\widetilde{Z}}$ on the connected components of
$X-\{\star_{\widetilde{Z}}\}$ has a unique 
orbit. In the case where $Z$ is a pair of pants, the group $\Gamma_{\widetilde{Z}}$ is a free group of rank $2$ and there are three 
$\Gamma_{\widetilde{Z}}$-orbits
of boundary components of $\widetilde{Z}$. Hence, the action of $\Gamma_{\widetilde{Z}}$ on the connected components of $X-\{\star_{\widetilde{Z}}\}$ has three
orbits. 

Conversely, if $Z$ is a connected component of $\Sigma_0$ which is not a cylinder, then $\Gamma_{\widetilde{Z}}$ is neither trivial nor cyclic, it fixes a point
$\star_{\widetilde{Z}}$ in $X$, either $\Gamma_{\widetilde{Z}}$
is a free group of rank at least $3$ or a free group of rank $2$ and the action of $\Gamma_{\widetilde{Z}}$ on $X-\{\star_{\widetilde{Z}}\}$ has one orbit, and
${Z}$ is maximal (for the inclusion) for these properties. Hence, the connected components of $\Sigma_0$ which are not cylinder are exactly the connected
closed subsurface of
$\Sigma$ whose stabilizers of the connected components of the preimage in $\widetilde{\Sigma}$ satisfy these properties, and the action of $\Gamma$ on $(X,d)$ determines
the connected components of $\Sigma_0$ which are not cylinders, up to isotopy.

Let us show that the action of $\Gamma$ on $X$ determines the cylinder components of $\Sigma_0$.  No cylinder connected component
of $\Sigma_0$ can be homotoped to a boundary component of another connected component of $\Sigma_0$. Hence, if $Z$ is a cylinder component of $\Sigma-\Sigma_0$,
it is not bounded by a cylinder component of $\Sigma_0$ but by two connected components of $\Sigma_0$ whose Euler characteristics are negative. And
$Z$ is a cylinder component of $\Sigma_0$ if and only if $\Gamma_{\widetilde{Z}}$ is cyclic and  pointwise preserves an edge (possibly reduced to a point) in $X$ and $Z$
is not a connected component of $\Sigma-\Sigma_0$. Hence, the action of $\Gamma$ on $(X,d)$ determines the cylinder components of $\Sigma_0$, up to isotopy.

\medskip

Let us show that the action of $\Gamma$ on $X$ determines the half-translation structures and the measured hyperbolic laminations, on the connected components of 
$\Sigma-\Sigma_0$ that are neither a cylinder nor a pair of pants. Let $\widetilde{W}$ be such a connected component of $\widetilde{\Sigma}-\widetilde{\Sigma}_0$.
If the piece $X_{\widetilde{W}}$ is the $\RR$-tree dual to a measured hyperbolic lamination $(\Lambda,\mu)$, then 
$(\Lambda,\mu)$ is determined (up to isotopy) by the set of translation distances of the elements of $\Gamma_{\widetilde{W}}$ (see \cite{Bonahon97}).
If $X_{\widetilde{W}}$ is the completion
of $(\widetilde{W},[\widetilde{q}_{\widetilde{W}}])$ as above, then $(W,[q_W])=(\widetilde{W},[\widetilde{q}_{\widetilde{W}}])/\Gamma_{\widetilde{W}}$
is determined by the action of $\Gamma_{\widetilde{W}}$ on $(\widetilde{W},[\widetilde{q}_{\widetilde{W}}])$.

\medskip

Let us show that the space $(X,d)$ and the action of $\Gamma$ are determined by the data of $\Sigma_0$, of the half-translation structures and the measured hyperbolic
laminations on the connected components of $\Sigma-\Sigma_0$ that are neither a cylinder nor a pair of pants, and by the lengths of the edges, up to $\Gamma$-equivariant
isometry. Assume that $(X,d)$ and $(X',d')$
correspond to two constructions as above, 
with the same subsurface $\Sigma_0$, the same  half-translation structures and measured 
hyperbolic laminations, and the same lengths of the edges.
Let $\widetilde{W}$ be a connected component
of $\widetilde{\Sigma}-\widetilde{\Sigma}_0$ that is not a pair of pants. Then, if the piece $X_{\widetilde{W}}$ is the completion of the universal cover 
of $(W,[q_W])$,
it is determined by $(W,[q_W])$ up to a $\Gamma_{\widetilde{W}}$-equivariant isometry. Similarly, if $X_{\widetilde{W}}$ is the $\RR$-tree dual to a measured hyperbolic lamination
$(\Lambda,\mu)$, it is determined by $(\Lambda,\mu)$ up to a $\Gamma_{\widetilde{W}}$-equivariant isometry. 
In both cases, the map that associates to the fixed point of a non trivial elliptic element $\gamma\in\Gamma_{\widetilde{W}}$ in the piece $X_{\widetilde{W}}$ of $X$,
the 
fixed point of $\gamma$ in the piece $X'_{\widetilde{W}'}$ of $X'$ (the fixed points exist and are unique according to Lemmas \ref{pointfixelamination} and
\ref{pointfixestructureplate}), extends in a unique way into a $\Gamma_{\widetilde{W}}$-equivariant isometry. Hence, the map
$\Phi$ which associates to the set of fixed point(s) 
of a non trivial elliptic element $\gamma\in\Gamma-\{e\}$ in $X$, the set of fixed point(s) 
of $\gamma$ in $X'$, extends in a unique way in an isometry, that defines an $\Gamma_{\widetilde{W}}$-equivariant isometry between $X_{\widetilde{W}}$ and
$X'_{\widetilde{W}}$, for every connected component $\widetilde{W}$ of $\widetilde{\Sigma}_0$ or of $\widetilde{\Sigma}-\widetilde{\Sigma}_0$, and since 
$\gamma(X_\star)=X_{\gamma\star}$, and $\gamma(X'_\star)=X'_{\gamma\star}$ for every $\gamma\in\Gamma$, the global isometry is $\Gamma$-equivariant.\cqfd 

\medskip

\bdefiA\label{defmixedstructure} A tree-graded metric space as defined above is called a mixed structure. 

\edefiA

Let $\gamma\in\Gamma-\{e\}$ and let $\alpha_{\gamma}$ be the free homotopy class of closed curves defined by $\gamma$ in $\Sigma$. Let $W$ be 
a non trivial $\pi_1$-injective connected (open) subsurface of
$\Sigma$. We identify $\alpha_\gamma$ and $W$ with their geodesic representative and geodesic realization, for any hyperbolic metric on $\Sigma$, and 
we denote by $\alpha_{\gamma,W}$ the closed curve 
or the union of essential arcs between two (possibly equal) punctures of $W$, which is the intersection of $\alpha_\gamma$ and $W$. We still denote by
$\alpha_{\gamma,W}$ the isotopy class (relative to the punctures) of $\alpha_{\gamma,W}$ in $W$ (that does not depend on the choice of the hyperbolic metric). 
If $W$ is endowed with a half-translation 
structure $[q_W]$ as above, let $\ell_{[q_W]}(\alpha_{\gamma,W})$ be the length of a $[q_W]$-geodesic representative
of $\alpha_{W}$.

Let $(X,\PPPP)$ be a mixed structure on $\Sigma$, and let $m$ be a hyperbolic metric on $\Sigma$. Let $(\Lambda,\mu)$ be the union of the 
measured hyperbolic laminations on $\Sigma$ which are dual to the $\RR$-tree pieces ($\Lambda$ has no closed leaf), let $c_1,\dots,c_n$ be the boundary
components of $\Sigma_0$, and let
$t_1,\dots,t_n\in [0;+\infty[$ be the lengths of the corresponding edges in $(X,d)$. 
For every $\gamma\in\Gamma-\{e\}$, let $\ell_{[q]}(\alpha_\gamma)$ be the sum of the lengths $\ell_{[q_W]}(\alpha_{\gamma,W})$
of the geodesic representatives of the arcs $\alpha_{\gamma,W}$, associated with the different pieces as above.    

\blemmA\label{MoonDuchin} For every $\gamma\in\Gamma-\{e\}$, 
the translation length of $\gamma$ in $(X,d)$ is $\ell_{X}(\gamma)=i((\Lambda,\mu),\alpha_\gamma)+\ell_{[q]}(\alpha_\gamma)+\sum_{i=1}^nt_ii(c_i,\alpha_\gamma)$.
\elemmA

\demA Let $\gamma\in\Gamma-\{e\}$. Assume first that $\gamma$ belongs to the stabilizer $\Gamma_{\widetilde{W}}$ of a connected component 
$\widetilde{W}$ of 
$\widetilde{\Sigma}_0$ or
of $\widetilde{\Sigma}-\widetilde{\Sigma}_0$. Let $W=p(\widetilde{W})$ and let $X_{\widetilde{W}}$ be the piece preserved by $\Gamma_{\widetilde{W}}$.
If $X_{\widetilde{W}}$ is a point, then $\ell_{X}(\gamma)=0$, and since $W$ is a connected component of $\Sigma_0$ or a pair of pants in $\Sigma-\Sigma_0$, we have 
$i((\Lambda,\mu),\alpha_\gamma)+\ell_{[q]}(\alpha_\gamma)+\sum_{i=1}^nt_ii(c_i,\alpha_\gamma)=0$. If $X_{\widetilde{W}}$ is the dual tree to a measured hyperbolic lamination 
$(\Lambda',\mu')$ on $W$, then $\ell_{X}(\gamma)=i((\Lambda',\mu'),\alpha_\gamma)=i((\Lambda,\mu),\alpha_\gamma)+\ell_{[q]}(\alpha_\gamma)+\sum_{i=1}^nt_i
i(c_i,\alpha_\gamma)$. Similarly, if $X_{\widetilde{W}}$ is the completion of a 
universal cover of $(W,[q_W])$  as above, then $\ell_{X}(\gamma)=\ell_{[q_W]}(\alpha_\gamma)=i((\Lambda,\mu),\alpha_\gamma)+\ell_{[q]}(\alpha_\gamma)+\sum_{i=1}^nt_i
i(c_i,\alpha_\gamma)$.


\medskip

Assume next that $\gamma$ does not preserve any piece. 
We replace $(X,d)$ by the space $(X',d')$ obtained by replacing the edges of length $t_i$ (possibly equal to $0$), with $i\in\{1,\dots,n\}$, by edges of length
$t_i+1$.  

Then, according to the picture of the introduction, the quotient of $(X',d')$ by the equivalence relation $x\sim y$ if and only if $x,y$ belong to the same piece which is not an
edge, is the (simplicial) tree $T$ dual to the multicurve $((t_i+1)c_i)_{i\in\{1,\dots,n\}}$.
Let $\Ax_{\widetilde{m}}(\gamma)$ be the translation axis of $\gamma$ in $(\widetilde{\Sigma},\widetilde{m})$ and let 
$\widetilde{W}_{0},\widetilde{W}_{1},\dots,\widetilde{W}_k$ ($k\in\NN$), be the connected components of $\widetilde{\Sigma}_0$ and of
$\widetilde{\Sigma}-\widetilde{\Sigma}_0$, that are successively intersected by a fundamental
domain of $\Ax_{\widetilde{m}}(\gamma)$, for the action of $\gamma^{\ZZ}$, starting at a boundary component $\widetilde{c}$ of $\widetilde{\Sigma}_0$.
Then, the translation axis of $\gamma$ in $T$ has a fundamental domain (for the action of $\gamma^{\ZZ}$) starting at an endpoint
of the edge pointwise preserved by the
stabilizer $\Gamma_{\widetilde{c}}$ of $\widetilde{c}$ in $\Gamma$, that 
successively meets the vertices preserved by $\Gamma_{\widetilde{W}_0},\Gamma_{\widetilde{W}_1},\dots,\Gamma_{\widetilde{W}_k}$. By the  definition of $T$, the translation axis of $\gamma$ in $X'$ (which is unique
since it is not contained in a piece) has a fundamental domain that  starts at an endpoint 
of the edge pointwise preserved by $\Gamma_{\widetilde{c}}$, and 
successively meets the pieces $X_{\widetilde{W}_0},X_{\widetilde{W}_1},\dots,X_{\widetilde{W}_k}$ preserved by 
$\Gamma_{\widetilde{W}_0},\Gamma_{\widetilde{W}_1},\dots,\Gamma_{\widetilde{W}_k}$. Hence, by addition, we have $\ell_{X'}(\gamma)=
i((\Lambda,\mu),\alpha_\gamma)+\ell_{[q]}(\alpha_\gamma)+\sum_{i=1}^n(t_i+1)i(c_i,\alpha_\gamma)$, and then $\ell_{X}(\gamma)=
i((\Lambda,\mu),\alpha_\gamma)+\ell_{[q]}(\alpha_\gamma)+\sum_{i=1}^nt_ii(c_i,\alpha_\gamma)$.\cqfd

\medskip

Let $\Mix(\Sigma)$ be the set of $\Gamma$-equivariant isometry classes of mixed structures on $\Sigma$. We endow $\Mix(\Sigma)$ with the equivariant Gromov topology
defined as follows. Let $\E$ be a set of metric
spaces endowed with an isometric action of $\Gamma$. For any $X\in\E$, for any finite subset $K$ of $X$, for any finite subset $P$ of $\Gamma$ and for any $\varepsilon>0$,
let $\V(X,K,P,\varepsilon)$ be the set of elements $X'\in\E$ such that there exist a finite subset $K'\subseteq X'$ and a relation $\R\subseteq K\times K'$,
whose projections on $K$ and $K'$ are surjective, such that 
$$\forall x,y\in K\;\forall x',y'\in K'\;\forall \gamma\in P, \mbox{ if } x\R x' \mbox{ and } y\R y',\mbox{ then } \; |d(x',\gamma y')-d(x,\gamma y)|<\varepsilon.$$ 
The sets $\V(X,K,P,\varepsilon)$ span a topology on $\E$ called the {\it equivariant Gromov  topology} (see for instance \cite{Paulin88,Paulin09a}). 
The equivariant Gromov topology naturally defines a topology on the set of $\Gamma$-equivariant isometry classes of metric
spaces endowed with an isometric action of $\Gamma$, still called the equivariant Gromov topology, and we endow $\Mix(\Sigma)$ with this topology.

\blemmA\label{plongementmix} 
The map $X\mapsto(\ell_{X}(\gamma))_{\gamma\in\Gamma}$ from $\Mix(\Sigma)$ to $\RR_+^\Gamma$ is continuous and injective.
\elemmA

\demA First, we prove that the map is injective. Let $(X,\PPPP)$ and $(X',\PPPP')$ be two mixed structures on $\Sigma$,
let $\Sigma_0$ and $\Sigma_0'$ be the associated tight subsurfaces
of $\Sigma$ 
and let $\ell_{[q]}$, $(\Lambda,\mu)$ 
and $(t_ic_i)_{i\in\{1,\dots,n\}}$ (resp. $\ell_{[q']}$, $(\Lambda',\mu')$ and $(t'_ic'_i)_{i\in\{1,\dots,n'\}}$)
be the length functions, measured hyperbolic laminations and multicurves
defined by $X$ and $X'$ as in Lemma \ref{MoonDuchin}.
Assume that for every $\gamma\in\Gamma$, we have $\ell_{X}(\gamma)=\ell_{X'}(\gamma)$.
Let $\widetilde{W}_0$ be a connected component of $\widetilde{\Sigma}_0$, which is not a strip. The stabilizer $\Gamma_{\widetilde{W}_0}$ of $\widetilde{W}_0$ 
fixes a point
$\star_{\widetilde{W}_0}$ in $X$.
Since $\ell_X(\gamma)=0$ if and only if $\ell_{X'}(\gamma)=0$, all
the elements of $\Gamma_{\widetilde{W}_0}$ have a fixed point in $X'$. Let us show that they have a common fixed point. As in the proof of Lemma \ref{MoonDuchin}, 
let $(X'_2,d'_2)$ be the space obtained by replacing the edges of length $t'_i$ (possibly equal to $0$), with $i\in\{1,\dots,n'\}$, by edges of length
$t'_i+1$, and let $T'_2$ be the quotient space $(X'_2,d'_2)/\sim$, by the equivalence relation generated by $x\sim y$ if $x$ and $y$ belong to a common piece which is not an
edge. The quotient metric space $T'_2$ is the simplicial tree dual to $\{(t'_i+1)c'_i\}_{i\in\{1,\dots,n'\}}$, and it is endowed by the action of $\Gamma$ defined by the
action of $\Gamma$ on $X'_2$. All the elements of $\Gamma_{\widetilde{W_0}}$ have a fixed 
point in $T'_2$. According to a Lemma of Serre (see \cite[p.~271]{Shalen87}), the subgroup $\Gamma_{\widetilde{W}_0}$ has a global fixed point in $T'_2$, and by definition of the action of $\Gamma$ on
$X'_2$, and hence on $T'_2$, this point is a vertex of $T'_2$. Hence, $\Gamma_{\widetilde{W}_0}$ preserve a piece of $X'_2$, and hence it preserves a piece of $X'$. Since
all the elements of $\Gamma_{\widetilde{W}_0}$ have a fixed point in $X'$, by the definition of the action of $\Gamma$ on $X'$, this piece is a point, and
$\Gamma_{\widetilde{W}_0}$ preserve a connected component $\widetilde{W}'_0$
of $\widetilde{\Sigma}'_0$ or of $\widetilde{\Sigma}-\widetilde{\Sigma}'_0$ whose image $W'_0$ in $\Sigma$ is a pair of pants, and $W_0$ 
is contained in $W'_0$, up to isotopy. However, $W_0$ is either a torus minus a disk or we 
have $\chi(W_0)\leqslant-4$, hence $W_0$ cannot be contained in a pair of pants. Hence, up to isotopy, $W_0$ is contained in a connected component ${W}'_0$ of 
$\Sigma_0'$. Similarly, up to isotopy, any cylinder component of $\widetilde{\Sigma}_0$ is contained in a connected component of $\Sigma_0'$ or of 
$\Sigma-\Sigma_0'$ which is a pair of pants. In both cases, it is isotopic to a cylinder contained in a connected component of $\Sigma'_0$.
Hence, up to isotopy, any connected component $W_0$ of $\Sigma_0$ is contained in a connected component $W_0'$ of $\Sigma_0'$.

Conversely, there exists a connected component 
${Z}_0$ of ${\Sigma}_0$ such that ${W}_0'$ is contained in ${Z}_0$, up to isotopy. Since the connected components of 
${\Sigma}_0$ are pairwise disjoint, we have ${W}_0={W}_0'$. Hence, the connected components of $\Sigma_0$ and of $\Sigma_0'$ are equal,
up to isotopy, and ${\Sigma}_0'={\Sigma}_0$, up to isotopy.
Notably, we have $n=n'$ and $\{c_1,\dots, c_n\}=\{c'_1,\dots,c'_{n'}\}$. 

\medskip
Let $\widetilde{W}$ be a connected component of 
$\widetilde{\Sigma}-\widetilde{\Sigma}_0$. Let $X_{\widetilde{W}}$ (resp. $X'_{\widetilde{W}}$) be the piece of $X$ (resp. $X'$) preserved by $\Gamma_{\widetilde{W}}$, 
and let $H_{\widetilde{W}}$ be the set of non trivial 
and non peripheral elements of $\Gamma_{\widetilde{W}}$. The pieces $X_{\widetilde{W}}$ and $X_{\widetilde{W}'}$ are reduced to a point if and only if $W=p(\widetilde{W})$
is a pair of pants.

If $X_{\widetilde{W}}$ is a tree dual to a measured hyperbolic lamination, according to Lemma
\ref{pasdeminoration}, for 
every $\varepsilon>0$, there exists $\gamma\in H_{\widetilde{W}}$ such that $\ell_{X}(\gamma)<\varepsilon$. If $X_{\widetilde{W}}$ is the completion of a
half-translation structure
on a surface, then there exists $\varepsilon>0$ such that $\ell_{X}(\gamma)\geqslant\varepsilon$ for every $\gamma\in H_{\widetilde{W}}$. 
And the same properties hold for $X'_{\widetilde{W}}$.
Hence, the piece $X_{\widetilde{W}}$ is the completion of a surface if and only if the piece $X'_{\widetilde{W}}$ is also the completion of a surface. 
Assume that $X_{\widetilde{W}}$ and $X'_{\widetilde{W}}$ are the completions of $(\widetilde{W},[\widetilde{q}_{\widetilde{W}}])$ and 
$(\widetilde{W},[\widetilde{q}'_{\widetilde{W}}])$.
Since $(\ell_{X}(\gamma))_{\gamma\in\Gamma_{\widetilde{W}}}=(\ell_{X'}(\gamma))_{\gamma\in\Gamma_{\widetilde{W}}}$, 
according to \cite[Thm.~1]{DucLeiRaf10}, the quotients $(\widetilde{W},[\widetilde{q}_{\widetilde{W}}])/\Gamma_{\widetilde{W}}$ and
$(\widetilde{W},[\widetilde{q}'_{\widetilde{W}}])/\Gamma_{\widetilde{W}}$ are isometric. Notably, the corresponding length functions  
$\ell_{X_{\widetilde{W}}}$ and $\ell_{X'_{\widetilde{W}}}$ are equal, and by
addition, we have $\ell_{[q]}=\ell_{[q']}$. Next, for every $\gamma\in\Gamma$, we have 
$$i((\Lambda,\mu),\alpha_\gamma)+\ell_{[q]}(\alpha_\gamma)+\sum_{i=1}^nt_ii(c_i,\alpha_\gamma)
=i((\Lambda',\mu'),\alpha_\gamma)+\ell_{[q']}(\alpha_\gamma)+\sum_{i=1}^{n}t'_i
i(c_i,\alpha_\gamma),$$
\noindent
and by soustraction $i((\Lambda,\mu),\alpha_\gamma)+\sum_{i=1}^nt_ii(c_i,\alpha_\gamma)
=i((\Lambda',\mu'),\alpha_\gamma)+\sum_{i=1}^{n}t'_i
i(c_i,\alpha_\gamma)$. Since $(\Lambda,\mu)\cup(t_i c_i)_{i\in\{1,\dots,n\}}$ and
$(\Lambda',\mu')\cup(t'_i c_i)_{i\in\{1,\dots,n\}}$ are some measured hyperbolic laminations and since the boundary of $\Sigma$ is empty, according to \cite[Exp.~7]{FatLauPoe79}, we have 
$(\Lambda,\mu)\cup(t_i c_i)_{i\in\{1,\dots,n\}}=(\Lambda',\mu')\cup(t'_i c_i)_{i\in\{1,\dots,n\}}$, and
since the laminations $\Lambda$ and $\Lambda'$ have no closed leaf, we have
$(\Lambda,\mu)=(\Lambda',\mu')$ and
$(t_i c_i)_{i\in\{1,\dots,n\}}=(t'_i c_i)_{i\in\{1,\dots,n\}}$. Finally, according to Lemma \ref{biendefini},
there exists a $\Gamma$-equivariant isometry between $(X,d)$ and $(X',d')$. Hence, the map is injective.

\medskip

Next, we prove that the map is continuous. 
It is sufficient to prove that for every $\gamma\in\Gamma-\{e\}$, the map $X\mapsto \ell_X(\gamma)$ from $\Mix(\Sigma)$ to $\RR^+$ is
continuous. Let $\gamma\in\Gamma$, 
let $X\in\Mix(\Sigma)$ and let $(X_n)_{n\in\NN}$ be a sequence of $\Mix(\Sigma)$ that converges to $X$.  

\medskip
\noindent
{\bf Case 1.}~  Assume that $\gamma$ is elliptic in $X$. Let $x\in X$ be a fixed point of $\gamma$. Let $P=\{\gamma\}$, $K=\{x\}$ and $\varepsilon>0$. Then, for $n$ large enough there
exists $x_n\in X_n$ such that $d(x_n,\gamma x_n)<\varepsilon$. Hence $\ell_{X_n}(\gamma)<\varepsilon$.

\medskip
\noindent
{\bf Case 2.}~ Assume that $\gamma$ is hyperbolic in $X$. Let $\Ax_X(\gamma)$ be a translation axis of $\gamma$ in $X$ and let $x\in\Ax_X(\gamma)$. Then, we have
$d(x,\gamma x)=\ell_X(\gamma)$ and $d(x,\gamma^2 x)-d(x,\gamma x)-d(\gamma x,\gamma^2 x)=0$. Let $K=\{x\}$ and 
$P=\{\gamma,\gamma^2\}$. Let $0<\varepsilon<\frac{1}{4}\ell_X(\gamma)$. If $n$ is large enough, there exists $x_n \in X_n$ such that 

$$|d(x_n,\gamma x_n)-d(x,\gamma x)|<\varepsilon \mbox{ and } |d(x_n,\gamma^2 x_n)-d(x,\gamma^2 x)|<\varepsilon$$
$$\mbox{so that }\,|d(x_n,\gamma^2 x_n)-d(x_n,\gamma x_n)-d(\gamma x_n,\gamma^2 x_n)|<3\varepsilon.$$

Let $\PPPP_n$ be the set of pieces covering $X_n$, let $\Sigma_{0n}$ be the tight subsurface associated with $X_n$ as in Lemma \ref{biendefini},
and let $\widetilde{\Sigma}_{0n}$ be its preimage in $\widetilde{\Sigma}$. 
Assume for a contradiction that $\gamma$ is elliptic in $X_n$. Then $\gamma$ belongs to the stabilizer 
$\Gamma_{\widetilde{W}}$ of a connected component $\widetilde{W}$ of $\widetilde{\Sigma}_{0n}$ or of $\widetilde{\Sigma}-\widetilde{\Sigma}_{0n}$ whose image in $\Sigma$ 
is a pair of pants. Either $\gamma$
does not preserve any boundary component of $\widetilde{W}$, and it fixes a unique point $\star_{\widetilde{W}}$ in $X_n$, or $\gamma$ preserves a boundary component
of $\widetilde{W}$. In the former case, it fixes an edge (possibly reduced to a point), and it preserves another piece $Y_n$ of $X_n$. According to Lemma 
\ref{pointfixelamination} and 
\ref{pointfixestructureplate}, any element of $\Gamma-\{e\}$ has at most one fixed point in any piece of $X_n$ that is not an edge. Hence, up to (possibly) replacing
the fixed point 
$\star_{\widetilde{W}}$ of 
$\gamma$ by the good endpoint of the edge, the point $\star_{\widetilde{W}}$ is the
unique common point between the segments $[x_n,\star_{\widetilde{W}}]$ and $[\star_{\widetilde{W}},\gamma x_n]$, and by construction of $X_n$, 
the union $[x_n,\star_{\widetilde{W}}]\cup[\star_{\widetilde{W}},\gamma x_n]$ is a geodesic segment, equal to $[x_n,\gamma x_n]$.
Similarly, the segment  $[x_n,\gamma^2 x_n]$ is the union $[x_n,\star_{\widetilde{W}}]\cup[\star_{\widetilde{W}},\gamma^2 x_n]$. 
Hence, we have $|d(x_n,\gamma^2 x_n)-d(x_n,\gamma x_n)-d(\gamma x_n,\gamma^2 x_n)|=2\, d(x_n,\star_{\widetilde{W}})$, hence 
$d(x_n,\star_{\widetilde{W}})<\frac{3}{2}\varepsilon$ and $d(x_n,\gamma x_n)<3\varepsilon$, a contradiction. Hence $\gamma$ is hyperbolic in $X_n$. 

\medskip 

If $x_n$ belongs to a translation axis of $\gamma$ in $X_n$, then $d(x_n,\gamma x_n)=\ell_{X_n}(\gamma)$, and $|\ell_X(\gamma)-\ell_{X_n}(\gamma)|<\varepsilon$. Otherwise, let 
$x_{n\perp}$ be the orthogonal projection  of $x_n$ onto the (possibly degenerated) flat strip union of all the translation axes of $\gamma$ in $X_n$, and let 
$\Ax_n(\gamma)$ be the translation axis containing $x_{n\perp}$. The segment $[x_n,x_{n\perp}]$ and $[x_n,\gamma x_n]$ may share a initial segment.
Let $a$ be the last intersection point between $[x_n,x_{n\perp}]$ and $[x_n,\gamma x_n]$. Similarly, let $b$ be the last intersection point between
$[\gamma x_n,\gamma x_{n\perp}]$ and $[\gamma x_n, x_n]$. Assume for a contradiction that the segment $[x_n,\gamma x_n]$ is disjoint from $\Ax_n(\gamma)$. Then, the curve 
$[a,x_{n\perp}]\cdot[x_{n\perp},\gamma x_{n\perp}]\cdot[\gamma x_{n\perp},b]\cdot[b,a]$ is a simple geodesic quadrangle, hence 
it is contained in a piece $Y_n\in\PPPP_n$ (see \cite[Lem.~2.5]{Drutu05}). 

Let $y_n$ be in $[x_n,x_{n\perp}]\cap Y_n$. Then $\gamma y_n$ belongs to $Y_n\cap[\gamma x_n,\gamma x_{n\perp}]$. Since $\gamma x_{n\perp}$ is not equal to
$x_{n\perp}$, the piece $Y_n$ is neither a point nor an edge. If $Y_n$ is the $\RR$-tree dual to a filling lamination or if $Y_n$ is the completion
of a surface endowed with a half-translation structure, according to Lemma \ref{pointfixe}, the segment $[y_n,\gamma y_n]$ intersects
$\Ax_n(\gamma)$, a contradiction. Hence, the segment $[x_n,\gamma x_n]$ meets $\Ax_n(\gamma)$, and similarly the segments 
$[\gamma x_n,\gamma^2 x_n]$ and $[x_n,\gamma^2 x_n]$ meet $\Ax_n(\gamma)$. 

\medskip

Moreover, since there exists a unique segment joining two points in a $\CAT(0)$ geodesic space, the first points of intersection between $[x_n,\gamma x_n]$ and 
$\Ax_n(\gamma)$,
and $[x_n,\gamma^2 x_n]$ and $\Ax_n(\gamma)$ are equal. Similarly, the last points of intersection between $[\gamma x_n,\gamma^2 x_n]$ and $\Ax_n(\gamma)$,
and $[x_n,\gamma^2 x_n]$ and $\Ax_n(\gamma)$ are equal. Let $w_n$ (resp. $z_n$) be the last (resp. first) point of intersection between $[x_n,\gamma x_n]$ and $\Ax_n(\gamma)$
(resp. $[\gamma x_n,\gamma^2 x_n]$ and $\Ax_n(\gamma)$).

\begin{center}
	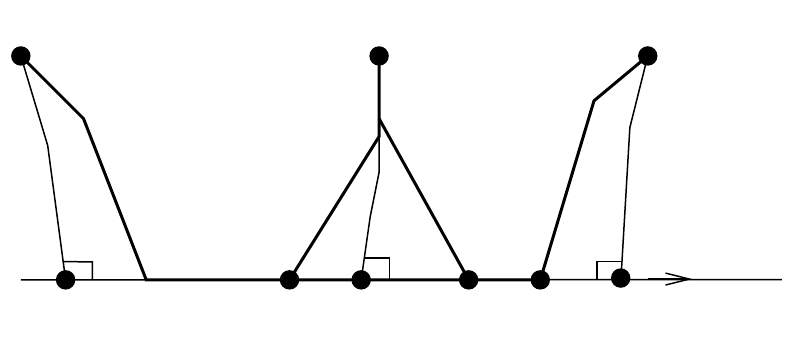
\end{center}

We have $|d(x_n,\gamma^2 x_n)-d(x_n,\gamma x_n)-d(\gamma x_n,\gamma^2 x_n)|=d(w_n,\gamma x_n)+d(z_n,\gamma x_n)-d(w_n,z_n)$. 
Moreover, since the space is $\CAT(0)$, we have 
$d(w_n,\gamma x_n)^2\geqslant d(w_n,\gamma x_{n\perp})^2+d(\gamma x_n,\gamma x_{n\perp})^2$ and $d(z_n,\gamma x_n)^2\geqslant d(z_n,\gamma x_{n\perp})^2+
d(\gamma x_n,\gamma x_{n\perp})^2$.
Hence \begin{align*}
	d(w_n,\gamma x_n)-d(w_n,\gamma x_{n\perp})&\geqslant\frac{d(\gamma x_n,\gamma x_{n\perp})^2}{d(w_n,\gamma x_n)+d(w_n,\gamma x_{n\perp})}\\
	&\geqslant \frac{d(x_n, x_{n\perp})^2}{2d(x_n,\gamma x_n)}
\end{align*}
Similarly, we have $d(z_n,\gamma x_n)-d(z_n,\gamma x_{n\perp})\geqslant\frac{d(x_n, x_{n\perp})^2}{2d(x_n,\gamma x_n)}$. Hence
$|d(x_n,\gamma^2 x_n)-d(x_n,\gamma x_n)-d(\gamma x_n,\gamma^2x_n)|\geqslant \frac{d(x_n,x_{n\perp})^2}{d(x_n,\gamma x_n)}$, and so 
$d(x_n,x_{n\perp})^2\leqslant 3\varepsilon d(x_n,\gamma x_n)\leqslant 3\varepsilon(\ell_X(\gamma)+\varepsilon).$
Hence $|d(x_n,\gamma x_n)- \ell_{X_n}(\gamma)|\leqslant 2\sqrt{3\varepsilon(\ell_X(\gamma)+\varepsilon)}$, so $|\ell_X(\gamma)-\ell_{X_n}(\gamma)|
\leqslant 2\sqrt{3\varepsilon(\ell_X(\gamma)+\varepsilon)} +\varepsilon$.\cqfd

\medskip

Let $W$ be a (open) non trivial, $\pi_1$-injective, connected subsurface of $\Sigma$, and let $\Flat(W)$ be the set of (isotopy classes)
of half-translation structures on $W$, that can be extended 
to half-translation structures on the compact surface obtained from $W$ by filling in the punctures, with possibly some singularities of angle $\pi$ at the added points.
We identify $\Flat(W)$ with the $\Gamma$-equivariant isometry classes of universal covers of half-translation structures on $W$ and we endow $\Flat(W)$ 
with the Gromov equivariant topology.
Let $\widetilde{W}\to W$ be a universal cover with covering group $\Gamma_{\widetilde{W}}$. 
For every element $[q_W]$ of $\Flat(W)$, we denote by 
$[\widetilde{q}_{\widetilde{W}}]$ its lift on $\widetilde{W}$. Let $H_{\widetilde{W}}$ be the set of non trivial and non peripheral elements of $\Gamma_{\widetilde{W}}$.

\blemmA\label{plongementflat} The map $\pi:(W,[{q}_W])\mapsto(\ell_{[\widetilde{q}_{\widetilde{W}}]}
(\gamma))_{\gamma\in\Gamma_{\widetilde{W}}}$
from $\Flat(W)$ to $(\RR^+)^{\Gamma_{\widetilde{W}}}-\{0\}$ is a homeomorphism onto its image.
\elemmA

\demA As in Lemma \ref{plongementmix}, this map is injective and continuous. Let $((x_{\gamma,n})_{\gamma\in\Gamma_{\widetilde{W}}})_{n\in\NN}$ be a sequence of 
the image of $\Flat(W)$ in $(\RR^{+})^{\Gamma_{\widetilde{W}}}-\{0\}$ that converges to $(x_\gamma)_{\gamma\in\Gamma_{\widetilde{W}}}$ in the image of $\Flat(W)$. Let
$({W},[{q}_n])_{n\in\NN}$ be the sequence of preimages of  $((x_{\gamma,n})_{\gamma\in\Gamma_{\widetilde{W}}})_{n\in\NN}$
in $\Flat(W)$. For all $n\in\NN$, let $\lambda_n$ be the minimal displacement of generator of $(\widetilde{W},
[\widetilde{q}_n])$, for a finite generating set of $\Gamma_{\widetilde{W}}$, as in Section \ref{sequenceofhalftranslations}. Assume for a contradiction that 
$\lim_\omega\lambda_n=+\infty$. Then, there exists a generator $s\in\Gamma$ such that $\lim_\omega\ell_{(\widetilde{W},[\widetilde{q}_n])}(s)=+\infty$, and 
$x_{s,n}$ does not $\omega$-converge to $x_s\in\RR^+$. Hence, we have $\lim_\omega\lambda_n<+\infty$, and
the asymptotic limit $\lim_\omega(\widetilde{W},[\widetilde{q}_n])_{n\in\NN}$, without renormalization, exists and is endowed with an isometric action of 
$\Gamma_{\widetilde{W}}$.
Since the point $(x_\gamma)_{\gamma\in\Gamma_{\widetilde{W}}}$ is in the image of $\Flat(W)$, there exists $\varepsilon'>0$ such that
for all $\gamma\in H_{\widetilde{W}}$, we have $x_\gamma>\varepsilon'$. Hence, according to Lemma \ref{minorationuniforme}
(that immediately extends to the universal cover of a finite type surface), there
exist $\varepsilon>0$ and $I\in\omega$ such that for all $\gamma\in H_{\widetilde{W}}$ and $n\in I$, 
we have $\ell_{(\widetilde{W},[\widetilde{q}_n])}(\gamma)>\varepsilon$. Hence, according to Section \ref{nondegeneration} (that immediately extends to the universal cover of a finite
type surface), the space $\lim_\omega(\widetilde{W},[\widetilde{q}_n])_{n\in\NN}$, minus the fixed points of the peripheral elements of $\Gamma_{\widetilde{W}}$,
is the universal cover of $W$ endowed with a half-translation structure, that extends to the punctures. Moreover, the ultralimit $\lim_\omega(\widetilde{W},[\widetilde{q}_n])_{n\in\NN}$ is the limit of a subsequence of 
$(\widetilde{W},[\widetilde{q}_n])_{n\in\NN}$ for the Gromov equivariant topology (see \cite{Paulin89a}). 
Hence, every subsequence of $(W,[q_n])_{n\in\NN}$ has a subsequence that converges to an 
element of $\Flat(W)$,
equal to $({W},[{q}_W])$ by injectivity of $\pi$. Hence the sequence $(W,[q_n])_{n\in\NN}$ converges to $({W},[{q}_W])$, and the inverse map is continuous.\cqfd
%


\blemmA\label{topologieseparable} The equivariant Gromov topology on $\Mix(\Sigma)$ is metrizable.
\elemmA

\demA According to Lemma \ref{plongementmix}, the space $\Mix(\Sigma)$ is Hausdorff and with countable bases of neighborhoods of points. Let us prove that it is separable. 
Let $\Sigma_0$ be a tight subsurface of $\Sigma$, let $\widetilde{\Sigma}_0$ be its preimage in $\widetilde{\Sigma}$, 
let $A$ be the set of connected components of $\widetilde{\Sigma}-\widetilde{\Sigma}_0$ whose image in $\Sigma$ is not a pair of pants, 
and let $B$ be the finite set of $\Gamma$-orbits of boundary components of $\widetilde{\Sigma}_0$. Let $E(\Sigma_0)$ 
be the subset of $\Mix(\Sigma)$ whose underlying tight subsurface is $\Sigma_0$. 
For every connected component $\widetilde{W}$ of $\widetilde{\Sigma}-\widetilde{\Sigma}_0$, let $W=p(\widetilde{W})$ and let $F_1(\widetilde{W})$ be the set of
completions of universal covers of (isotopy classes of) half-translation structures on $W$, and $F_2(\widetilde{W})$ be the set of $\RR$-trees dual
to a filling measured hyperbolic lamination on $W$. We endow $F_1(\widetilde{W})$ and $F_2(\widetilde{W})$ with the Gromov equivariant topology.
Let $\Phi$ be the map from $\underset{\widetilde{W}\in A}\prod(F_1(\widetilde{W})\coprod F_2(\widetilde{W}))\times(\RR^+)^B$, endowed with the product topology, 
to $\Mix(\Sigma)$, that associates
to $((X_{\widetilde{W}})_{\widetilde{W}\in A},(t_b)_{b\in B})$ the $\Gamma$-equivariant isometry classes of mixed structures such that the piece preserved by 
$\Gamma_{\widetilde{W}}$, with $\widetilde{W}\in A$, is $X_{\widetilde{W}}$, and the length of the edges pointwise preserved by the stabilizers of the elements of $b\in B$
is $t_b$ (that exists and is unique under the few conditions before Lemma \ref{biendefini}).

Moreover, for every element $\widetilde{W}\in A$, according to \cite[Thm.~5.2]{Paulin89a},  the space $F_2(\widetilde{W})$ is separable,
and according to Lemma \ref{plongementflat}, the space $F_1(\widetilde{W})$ is separable. Since $(\RR^+)^B$ is separable and $A$ is countable, the space 
$\underset{\widetilde{W}\in A}\prod(F_1(\widetilde{W})\coprod F_2(\widetilde{W}))\times(\RR^+)^B$ is separable, and since the map $\Phi$ is  
continuous and surjective, the space $E(\Sigma_0)$ is separable. Finally, the space $\Mix(\Sigma)$ is the countable union of the sets $E(\Sigma_0)$,  
with $\Sigma_0$ a isotopy class of tight subsurface of $\Sigma$, hence $\Mix(\Sigma)$ is separable.
Since $\Mix(\Sigma)$ has countable bases of neighborhoods of points, it has a countable basis. Moreover, the space $(\RR^+)^\Gamma-\{0\}$ is normal, hence
according to Lemma \ref{plongementmix}, the space $\Mix(\Sigma)$ is also normal. According to a theorem of Urysohn
(see for instance \cite[Ch.~9, 9.2]{Dugundji}), it is metrizable.\cqfd

\medskip

The group $\RR^{+*}$ acts on $\Mix(\Sigma)$ by multiplications of the distances. Let $\PP\Mix(\Sigma)$ and $\PP\Flat(\Sigma)$ be the quotients of
$\Mix(\Sigma)$ and of $\Flat(\Sigma)$ by these actions, endowed with the quotient topology of the equivariant Gromov topology.


\btheoA\label{theorem1} The space $\PP\Flat(\Sigma)$ is an open and dense subset of $\PP\Mix(\Sigma)$, which is compact.  
\etheoA
%

\demA  Let $(\widetilde{\Sigma},[\widetilde{q}'_n])_{n\in\NN}$ be a sequence of $\Flat(\Sigma)$.
Let $S$ be a finite generating set of $\Gamma$, and for all $n\in\NN$, 
let $\lambda_n$ be the minimal displacement of the generators on $(\widetilde{\Sigma},[\widetilde{q}'_n])$ and let $\star_n$ be a point such that 
$\max\{d'_n(\star_n, s\star_n)\,:\,s\in S\}$ is less than $\lambda_n+1$, as in Section \ref{sequenceofhalftranslations}. 
Let $[\widetilde{q}_n]=\frac{1}{\lambda_n}[\widetilde{q}'_n]$, and let $[q_n]$ be the half-translation structure on $\Sigma$ defined by $[\widetilde{q}_n]$.
Let $\omega$ be a non principal ultrafilter on $\NN$ as in Section \ref{ultralimite}.
Then, the ultralimit $\lim_\omega(\widetilde{\Sigma},[\widetilde{q}_n],\star_n)_{n\in\NN}$ is a $\CAT(0)$, geodesic, complete space, endowed with an isometric action of 
$\Gamma$. 

\medskip

Let us show that  $\lim_\omega(\widetilde{\Sigma},[\widetilde{q}_n],\star_n)_{n\in\NN}$ has a $\Gamma$-invariant subset which is 
a mixed structure.
Let $\SSSS_0=\{\alpha\in\SSSS(\Sigma)\,:\,\lim_\omega\ell_{[q_n]}(\alpha)=0\}$ and let $\Sigma_0$ be the tight subsurface of $\Sigma$ filled up
by $\SSSS_0$, and let $\widetilde{\Sigma}_0$ be its preimage in $\widetilde{\Sigma}$. No connected component of $\Sigma_0$ is a pair of pants.
For all connected component ${\widetilde{W}}$ of $\widetilde{\Sigma}_0$ or of $\widetilde{\Sigma}-\widetilde{\Sigma}_0$ and $n\in\NN$, we denote by $\widetilde{W}_n$ the 
$[\widetilde{q}_n]$-geometric realization of $\widetilde{W}$. If $\widetilde{W}$ is a connected component
of $\widetilde{\Sigma}-\widetilde{\Sigma}_0$, let $W$ be the 
finite type surface, image of $\widetilde{W}$ in $\Sigma$. According to Lemmas \ref{elliptic}, \ref{caseodegeneration}, \ref{enveloppeconvexe}, 
\ref{pantalon}, the convex core of $\Gamma_{\widetilde{W}}$ in $[\widetilde{W}_n]_{n\in\NN}$ is a unique point $\star_{\widetilde{W}}$ if
$\widetilde{W}$ is a connected component of $\widetilde{\Sigma}_0$ or if $W$ is a pair of pants in $\Sigma-\Sigma_0$, and it is either 
the completion of a universal cover of $W$ endowed with a half-translation structure that can be extended 
(as in Section \ref{nondegeneration}), or the 
$\RR$-tree dual to a filling measured hyperbolic lamination on $W$ (for any complete hyperbolic metric). In the two cases, we denote it by $X_{\widetilde{W}}$.
Finally, if $\Gamma_{\widetilde{c}}$ is the stabilizer in $\Gamma$ of a boundary component $\widetilde{c}$ of $\widetilde{\Sigma}_0$,
we have seen in  Section \ref{cylindresetpantalons}, that $\Gamma_{\widetilde{c}}$ pointwise preserves a (maximal) flat strip (possibly reduced to a geodesic), and
there exists a geodesic segment $e_{\widetilde{c}}$, possibly reduced to a point, orthogonal to the boundary components of the flat strip (if any), such that
the intersection of the translation axis 
of any hyperbolic element of $\Gamma$ which is interlaced with $\widetilde{c}$ (see Section \ref{sequenceofhalftranslations} for the definition),
is $e_{\widetilde{c}}$. Then, the union $X$ of the points $\star_{\widetilde{W}}$, of the convex subsets $X_{\widetilde{W}}$ and of the geodesic 
segments
$e_{\widetilde{c}}$, 
where $\widetilde{W}$ is a connected component of $\widetilde{\Sigma}_0$ or of $\widetilde{\Sigma}-\widetilde{\Sigma}_0$ and $\widetilde{c}$ is a
boundary component of $\widetilde{\Sigma}_0$, is a convex tree-graded subspace of $\lim_\omega(\widetilde{\Sigma},[\widetilde{q}_n],\star_n)_{n\in\NN}$,
notably it is $\CAT(0)$. Moreover, it is naturally $\Gamma$-invariant, and hence endowed with an isometric action of $\Gamma$.
Moreover, if $\widetilde{W}$ is a connected component of the preimage of a cylinder
component $W$ of $\Sigma-\Sigma_0$, and if for all $n\in\NN$ we denote by $h_n(\widetilde{W})$ the height of the 
$[\widetilde{q}_n]$-geometric realization of $\widetilde{W}$, then $\lim_\omega h_n(\widetilde{W})>0$, otherwise there would exist an element of $\SSSS_0$ 
that would cut $W$ and that would be contained in the union of $W$ and of the two connected components of $\Sigma_0$ that bound $W$, and $W$ would be contained in $\Sigma_0$. Similarly, 
if $\widetilde{W}$ is a connected component of $\widetilde{\Sigma}-\widetilde{\Sigma}_0$ whose image in $\Sigma$ is a pair of pants, there exists at least a 
boundary component $\widetilde{c}$ of $\widetilde{W}$ such that $\lim_\omega h_n(\widetilde{c})>0$, with $h_n(\widetilde{c})$ the height of the flat strip union of all the 
$[\widetilde{q}_n]$-geometric representatives of $\widetilde{c}$, , otherwise $W$ would be contained in 
$\Sigma_0$.
Hence, the subspace $X$, endowed with the action of $\Gamma$, satisfies all the properties of Definition \ref{defmixedstructure}. 

\medskip

Moreover, by the definition of the ultralimits, 
the ultralimit $\lim_{\omega}(\widetilde{\Sigma},[\widetilde{q}_n],\star_n)_{n\in\NN}$ is a limit of a subsequence of 
$(\widetilde{\Sigma},[\widetilde{q}_n])_{n\in\NN}$ for the  equivariant Gromov topology, in the space of $\Gamma$-equivariant isometry classes
of metric spaces endowed with an
isometric action of $\Gamma$. Since $X$ is $\Gamma$-invariant in $\lim_{\omega}(\widetilde{\Sigma},[\widetilde{q}_n],\star_n)_{n\in\NN}$, the space $X$ is also
a limit of a subsequence of $(\widetilde{\Sigma},[\widetilde{q}_n])_{n\in\NN}$ for the Gromov equivariant topology. If the sequence $(\widetilde{\Sigma},[\widetilde{q}_n])_{n\in\NN}$ converges in
$\Mix(\Sigma)$, since $\Mix(\Sigma)$ is separated (by Corollary \ref{topologieseparable}), the space $X$ 
is its unique limit in $\Mix(\Sigma)$. 

\medskip

Moreover, by the geometric construction of a sequence in $\Flat(\Sigma)$ at the end of the proof of 
\cite[Thm.~6 p.~27]{DucLeiRaf10}, we can obtain any mixed structure  as a $\Gamma$-equivariant subspace of the
ultralimit of a sequence of $\PP\Flat(\Sigma)$. Let us recall this construction of \cite{DucLeiRaf10}. Let $(X,d)$ be a mixed structure on $\Sigma$ and let $\Sigma_0$ be the associated tight subsurface.
We build a sequence in $\Flat(\Sigma)$ piece by piece as follows. Let $W$ be a connected component of $\Sigma-\Sigma_0$ which is neither a cylinder nor a pair of pants and let 
$\widetilde{W}$ be a connected component of the preimage
of $W$ in $\widetilde{\Sigma}$. Then, the piece $X_{\widetilde{W}}$ is non empty. 

If the piece $X_{\widetilde{W}}$ is the completion of the universal cover of a half-translation structure $[q_W]$ on $W$, let $q_W\in\Q(W)$ be a representative of $[q_W]$
and let $(W_n,q_{W,n})$ be the half-translation surface obtained by cutting vertical slits at the punctures of $(W,q_{W})$ of length $\frac{1}{n^2}$.

If $X_{\widetilde{W}}$ is the dual tree to a filling (and minimal) measured hyperbolic lamination $(\Lambda_W,\mu_W)$ on $W$, let $q_W$ be a quadratic differential on $W$ whose
vertical measured foliation is associated with $(\Lambda_W,\mu_W)$ (by the map of \cite{Levitt81}),  and let $(W_n,q_{W,n})$ be the surface obtained from the half-translation 
surface $(W,\begin{pmatrix} 1 & 0 \\ 0 &\frac{1}{n^2}\end{pmatrix}\cdot q_W)$ by cutting vertical slits at the punctures of length $\frac{1}{n^2}$. 

Finally, let $q_0$ be a quadratic differential on the interior of the connected components of $\Sigma_0$ which are not cylinders, 
that can be extended at the punctures (with possibly singularities of angle $\pi$), and whose vertical foliation is minimal 
(on every connected component, which is always possible) 
and let us cut vertical 
slits of 
lengths $\frac{1}{n}$ at the punctures. Thus, we get a surface $(\Sigma_0',q'_{0})$. We set $(\Sigma_{0,n},q_{0,n})=(\Sigma_0',\frac{1}{n}q'_0)$.

For every free homotopy class $c$ of boundary components of $\Sigma_0$, let $C_{c,n}$ be the flat cylinder of height $\ell_X(X_{\widetilde{c}})$, where $X_{\widetilde{c}}$ is 
the edge associated to a lift $\widetilde{c}$ of $c$, and of girth $\frac{1}{n^2}$. Any free homotopy class of boundary components $c$ of $\Sigma_0$ bounds two
connected components 
$W$ and $W'$ of $\Sigma_0$ or of $\Sigma-\Sigma_0$ which are not cylinders. If $W$ and $W'$ are not pair of pants, for all $n\in\NN$, we glue isometrically the boundary components of $C_{c,n}$ on the corresponding slits of $W_n$ and
of $W'_n$.


Finally, if $W$ is a pair of pants, it is bounded by three boundary components of   $\Sigma_0$. We cover it by gluing the three corresponding cylinders, and we glue their 
other boundary components on the corresponding slits of $\Sigma_{0,n}$. 
\begin{center}
\begin{picture}(0,0)%
\includegraphics{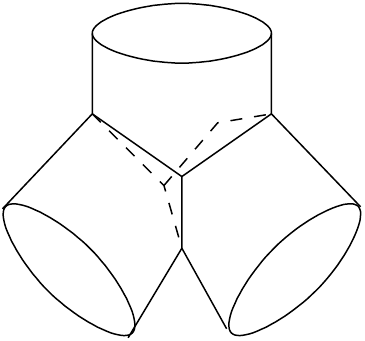}%
\end{picture}%
\setlength{\unitlength}{3771sp}%
\begingroup\makeatletter\ifx\SetFigFont\undefined%
\gdef\SetFigFont#1#2#3#4#5{%
  \reset@font\fontsize{#1}{#2pt}%
  \fontfamily{#3}\fontseries{#4}\fontshape{#5}%
  \selectfont}%
\fi\endgroup%
\begin{picture}(1823,1704)(888,-973)
\end{picture}%
 
\end{center}

We thus obtain a locally Euclidean metric on $\Sigma$ with conical singular points of angle $k\pi$, with $k\in\NN$ and $k\geqslant 3$. Moreover, all the slits on which are 
glued the flat cylinders are vertical. Hence, the surface has a singular foliation, whose leaves are vertical (with respect to the quadratic differentials 
on the pieces). Hence, this metric comes from  a
quadratic differential $q_n$ on $\Sigma$. Let $(\widetilde{\Sigma},[\widetilde{q}_n])_{n\in\NN}$ be the corresponding sequence of $\Flat(\Sigma)$. Then, 
the sequence converges to $X\in\Mix(\Sigma)$. Hence the closure of $\PP\Flat(\Sigma)$ for the  equivariant Gromov topology is $\PP\Mix(\Sigma)$.

\medskip

Let $d$ be a distance on $\PP\Mix(\Sigma)$ that induces the topology. Let $(X_n)_{n\in\NN}$ be a sequence in $\PP\Mix(\Sigma)$. For all $n\in\NN$, there exists a sequence  
$(X_{n,k})_{k\in\NN}$ in $\PP\Flat(\Sigma)$ such that $d(X_{n,k},X_n)<\frac{1}{k}$. We have seen that there exists $X\in\PP\Mix(\Sigma)$ such that 
the sequence $(X_{n,n})_{n\in\NN}$ $\omega$-converges to  $X$ in $\PP\Mix(\Sigma)$. 
Then, the sequence $(X_n)_{n\in\NN}$ $\omega$-converges to $X$. Hence, the space 
$\PP\Mix(\Sigma)$ is sequencially compact, and since it is Hausdorff, it is compact.

\medskip

Finally, let us prove that $\PP\Mix(\Sigma)-\PP\Flat(\Sigma)$ is closed. By the definition of the mixed structures, and according to Lemma \ref{pasdeminoration}, the set 
$\Mix(\Sigma)-\Flat(\Sigma)$ is the set of the elements $X\in\Mix(\Sigma)$ such that for all $\varepsilon>0$, there exists $\gamma\in\Gamma-\{e\}$ such that 
$\ell_X(\gamma)<\varepsilon$.
And a sequence $(Y_k)_{k\in\NN}$ of $\Flat(\Sigma)$ $\omega$-converges to an element of $\Flat(\Sigma)$ if and only if the $\omega$-limits of the translation distances of the 
non trivial elements of 
$\Gamma$ are finite and uniformely bounded below by the injectivity radius of the limit of the sequence in $\Flat(\Sigma)$. 
According to Lemma \ref{minorationuniforme}, this happens
if and only if there exists $\varepsilon>0$ such that for all $\gamma\in\Gamma-\{e\}$ and $\omega$-almost all $n\in\NN$, we have 
$\ell_{Y_k}(\gamma)\geqslant\varepsilon$ (and $\lim_\omega\ell_{Y_k}(\gamma)<+\infty$).

\medskip

Let $(X_n)_{n\in\NN}$ be a sequence of $\Mix(\Sigma)-\Flat(\Sigma)$ that converges to $X\in\Mix(\Sigma)$, and for all $n\in\NN$
let $(X_{n,k})_{k\in\NN}$ be a sequence in $\Flat(\Sigma)$ that converges to $X_n$. For all $n\in\NN$, there exists 
$\gamma_n\in\Gamma-\{e\}$ such that 
$\ell_{X_n}(\gamma_n)<\frac{1}{n}$. Since the map $X\mapsto\ell_X(\gamma)$ is continuous, 
there exists $k(n)\in\NN$ such that $\ell_{(X_{n,k(n)})}(\gamma_n)<\frac{1}{n}$. Since the sequence $(X_{n,k(n)})_{n\in\NN}$ converges to $X$ in
$\Mix(\Sigma)$, $X$ does not belong to $\Flat(\Sigma)$. Hence $\Mix(\Sigma)-\Flat(\Sigma)$ is closed, and so is 
$\PP\Mix(\Sigma)-\PP\Flat(\Sigma)$.\cqfd

\btheoA\label{theorem2} The map $[X]\mapsto[\ell_X(\gamma)]_{\gamma\in\Gamma}$ from $\PP\Mix(\Sigma)$ to $\PP\RR^\Gamma$ is a homeomorphism onto its image, equivariant
under the action of the mapping class group of $\Sigma$.
\etheoA

\demA It is a consequence of Lemma \ref{plongementmix} and Theorem \ref{theorem1}.\cqfd

\medskip

In \cite[Thm.~4]{DucLeiRaf10}, Duchin-Leininger-Rafi define an embedding $[q]\mapsto L_{[q]}$ of $\Flat(\Sigma)$ into the space $\CCC(\Sigma)$
of geodesic currents on $\Sigma$ (see for instance \cite{Bonahon88} for the definition of the geodesic currents), uniquely defined such that for every 
$\alpha\in\C(\Sigma),\; i(L_{[q]},\alpha)=\ell_{[q]}(\alpha)$. Moreover, in \cite[§.~5]{DucLeiRaf10}, they define 
the space $\M(\Sigma)$ of mixed structures on $\Sigma$ as the subset of $\CCC(\Sigma)$ of geodesic currents which are the sum of a geodesic current $L_{[q']}$,
defined by $i(L_{[q']},\alpha_\gamma)=\ell_{[q']}(\alpha_\gamma)$ for every $\gamma\in\Gamma$
where $[q']$ is a half-translation
structure on a subsurface $\Sigma'$ of $\Sigma$ (possibly equal to $\Sigma$), whose connected components are $\pi_1$-injective and have a negative Euler characteristic,
and of a measured hyperbolic laminations (for any hyperbolic metric on $\Sigma$) whose support is disjoint from $\Sigma'$ (up to isotopy). 
Let $\PP\CCC(\Sigma)=(\CCC(\Sigma)-\{0\})/\RR^{+*}$ and $\PP\M(\Sigma)$ be the image of $\M(\Sigma)-\{0\}$ in $\PP\CCC(\Sigma)$. They prove at \cite[Thm.~6]{DucLeiRaf10}
that the closure of the image of $\PP\Flat(\Sigma)$ in $\PP\CCC(\Sigma)$ is exactly $\PP\M(\Sigma)$. Let $\Phi$ be the map from $\Mix(\Sigma)$
to $\CCC(\Sigma)$ that associates to the lengths function $\ell_{[q]}$,
the measured hyperbolic lamination $(\Lambda,\mu)$ and the multicurve $(t_ic_i)_{i\in\{1,\dots,n\}}$ associated with 
$(X,\PPPP)\in\Mix(\Sigma)$, the geodesic current $L_{[q]}+(\Lambda,\mu)+(t_ic_i)_{i\in\{1,\dots,n\}}$
(where the measured hyperbolic laminations are seen as geodesic currents).

\blemmA\label{dualitee} The map $\Phi$ is a homeomorphism which is an extension to $\Mix(\Sigma)$ of the map $[q]\mapsto L_{[q]}$.  
\elemmA

\demA It is a consequence of the fact that the map $\mu\mapsto(i(\mu,\alpha))_{\alpha\in\C(\Sigma)}$ is an embedding of $\CCC(\Sigma)$ into its image in
$(\RR^{+})^{\C(\Sigma)}$ (see for instance \cite[Thm.~10,11]{DucLeiRaf10} and of Theorem \ref{theorem2}.\cqfd
\bibliographystyle{alphanum}
\bibliography{biblio}{}
\end{document}